\documentclass{article}
\usepackage{latexsym,amsmath,amssymb,amscd}
\usepackage[dvips]{color}
\numberwithin{equation}{section}

\begin{document}

\title{On the axisymmetric metric generated by a rotating perfect fluid with the vacuum boundary}
\author{Tetu Makino \footnote{Professor Emeritus at Yamaguchi University, Japan. E-mail: makino@yamaguchi-u.ac.jp}}
\date{\today}
\maketitle

\newtheorem{Lemma}{Lemma}
\newtheorem{Proposition}{Proposition}
\newtheorem{Theorem}{Theorem}
\newtheorem{Definition}{Definition}

\newtheorem{Remark}{Remark}
\newtheorem{Corollary}{Corollary}
\newtheorem{Notation}{Notation}

\begin{abstract}
We consider the equations for the coefficients of stationary rotating axisymmetric metrics governed by the Einstein-Euler equations, that is, the Einstein equations together with the energy-momentum tensor of a barotropic perfect fluid. 
Although the reduced system of equations for the potentials in the co-rotating co-ordinate system is known, we derive the system of equations for potentials in the so called zero angular momentum observer co-ordinate system. We newly give a proof of the equivalence between the reduced system and the full system of Einstein equations. It is done under the assumption that the angular velocity is constant on the support of the density. Also the consistency of the equations of the system is analyzed. On this basic theory we construct on the whole space the stationary asymptotically flat metric generated by a slowly rotating compactly supported perfect fluid with vacuum boundary. \\



\end{abstract}

\section{Introduction}

In this article we discuss the Einstein equations with the energy-momentum tensor of a rotating compactly supported perfect fluid. The fluid is supposed to be barotropic and the equation of state is approximately the $\gamma$-law near the vacuum. Mathematically rigorous treatment of this situation requires delicate analysis. Of course we seek an axially symmetric stationary metric, but moreover, we want to construct the metric globally so that it turns out to be asymptotically flat infinitely far away. 

For the special case of spherically symmetric metric for static (non-rotating) fluid source, the problem can be said to be already solved. In this case the density distribution and the metric on the support of density is given by solving the Tolman-Oppenheimer-Volkoff equation ( see T.M. \cite{TM1998}, A. D. Rendall and B. G. Schmidt \cite{RendallS} ) and the metric can be extended using the exterior Schwarzschild metric in $C^2$-fashion across the vacuum boundary ( see \cite[Supplementary Remark 4]{ssEE} ). Of course the Schwarzschild metric is asymptotically flat. 

However the treatment of the rotating case is harder. Historically speaking, the problem has been attacked in the direction to try to find an interior metric on the support of density under the situation that an axially symmetric vacuum metric on the exterior domain is given so that the interior metric be matched well with the given exterior one across the vacuum boundary. Naturally one supposes that the given exterior vacuum metric is of the Kerr metric, which is the rotating version of the Schwarzschild metric, although any uniqueness theorem analogous to the Birkhoff's theorem for the spherically symmetric static metric is not known. But this strategy, of e.g., the work by W. Roos \cite{Roos}, is involved in the difficulty of determination of the
location of the matching boundary, which should be a free boundary determined by the unknown interior solution. It must be avoided to go round and round in circle. Anyway this is the problem of so called `a source of the Kerr field'. Several scholars have attacked it but not yet solved it. The textbook \cite{PlebanskiK} by J. Pleba\'{n}ski and A. Krasi\'{n}ski, 2006, pp.499-495, says

\begin{quote}
The Kerr solution has been known for more than 40 years now, and from the very beginning
its existence provoked the simple question: what material body could generate such a vacuum field around it? Several authors have tried very hard to find a model of the source, but so far without success. The most promising positive result is that of Roos (1976), who investigated the Einstein equations with a perfect fluid with the boundary condition that the Kerr metric is matched to the solution. All attempts  so far to find an explicit example of a solution failed. The continuing lack of success prompted some authors to spread the suspicion that a perfect fluid source might not exist; rumours about this suspicion were then taken as a serious suggestion. The opinion of one of the present authors (A. K.) is that a bright new idea is needed, as opposed to routine standard tricks tested so far. 
\end{quote}

On the other hand the success of the solution of the problem for spherically symmetric case suggests the other strategy of the opposite direction to try to find the metric on the exterior vacuum region under that situation that an interior metric on the support of the density is already given. But mathematically rigorous construction of stationary rotating compactly supported compressible fluid mass distributions is far from being an easy task even for the non-relativistic situation, say, for the solutions governed by the Euler-Poisson equations. Standing upon the  results of Juhi Jang and T. M. \cite{JJTM1}, \cite{JJTM2}, the counter-parts on the non-relativistic settings, T. M. \cite{asEE} constructed axially symmetric metrics with compactly supported rigidly rotating density distributions on a bounded domain which includes the supports of the density distributions. In the non-relativistic problem, once solutions on a bounded domain which includes the support of the density $\rho$ are constructed, the problem is already solved, since we can assume arbitrary velocity filed on the whole vacuum region  merely by putting $\rho=0$ outside the bounded domain. However the situation is not so self-evident for the relativistic problem. 
The Lense-Thirring's dragging effect ( see e.g., \cite[ \S 12.18]{PlebanskiK} ) could be neglected nowhere in the vacuum region even far away from the rotating source.
It was difficult to find a global extension of the metric infinitely far away so that it be asymptotically flat. The difficulty comes from the usage of the co-rotating co-ordinate system in the work \cite{asEE}. The co-rotating co-ordinate system seems not to work well globally. Therefore in this article we try to use not the co-rotating co-ordinate system but the so called `ZAMO (= zero angular momentum observer )
co-ordinate system' On this ground we can construct an asymptotically flat axysymmetric metric generated by a rigidly and slowly rotating perfect fluid source with relatively small central density. In order to do it we should prepare a basic analysis of the system of equations for the potentials of the axisymmetric metric described in the ZAMO co-ordinate system. Namely we should perform the analysis of the equivalence of the reduced system to the full system of Einstein equations, and the consistency of the equations and so on. This basic analysis seems not to be found in the existing literatures.\\

Let us describe the settings precisely. We are considering the metric 
\begin{equation}
ds^2=g_{\mu\nu}dx^{\mu}dx^{\nu}
\end{equation}
satisfying the Einstein equations
\begin{equation}
R_{\mu\nu}-\frac{1}{2}g_{\mu\nu}R=
\frac{8\pi\mathsf{G}}{\mathsf{c}^4}T_{\mu\nu}.
\end{equation}
The positive constants $\mathsf{c}, \mathsf{G}$ are the speed of light and the constant of gravitation. 
The energy-momentum tensor $T^{\mu\nu}$ is supposed to be that of a perfect fluid
\begin{equation}
T^{\mu\nu}=(\epsilon +P)U^{\mu}U^{\nu}-Pg^{\mu\nu}
\quad
\mbox{with}\quad \epsilon =\mathsf{c}^2\rho.
\end{equation}

Here we put the following assumption:\\

{\bf (A)}: {\it The pressure $P$ is  a  smooth function of $\rho>0$ such that
 $0<P,
0<dP/d\rho <\mathsf{c}^2$ for $\rho >0$ and there are constants $\gamma, \mathsf{A}$ and a function $\Upsilon$ which is analytic near $0$ and satisfies $\Upsilon(0)=0$ such that
\begin{equation}
P=\mathsf{A}\rho^{\gamma}(1+\Upsilon(\mathsf{A}\rho^{\gamma-1}/\mathsf{c}^2)) \label{0104}
\end{equation}
for $\rho > 0$ and $1<\gamma < 2, \mathsf{A} >0$.
 }\\

The typical example which we keep in mind for the assumption {\bf (A)} is the equation of  neutron stars (\cite[Page 188, (6.8.4), (6.8.5)]{Zeldovich}):
\begin{equation}
P=B\mathsf{c}^5\int_0^Q\frac{q^4dq}{\sqrt{1+q^2}},\quad
\rho=3B\mathsf{c}^3\int_0^Q\sqrt{1+q^2}q^2dq,
\end{equation}
$Q$ being an auxiliary variable, $0\leq Q < +\infty$, $B$ being a positive constant. 
Actually this equation of state enjoys \eqref{0104}
with $\displaystyle \gamma=\frac{5}{3}, \mathsf{A}=\frac{1}{5B^{2/3}}$, and enjoys
$$\frac{dP}{d\rho}=\frac{\mathsf{c}^2}{3}\frac{Q^2}{1+Q^2} < \frac{\mathsf{c}^2}{3}\quad\mbox{for}\quad 0\leq Q <+\infty.$$
We avoid the simple $\gamma$-law $P=\mathsf{A}\rho^{\gamma}, \gamma >1$, since
the speed of sound $\displaystyle \sqrt{\frac{dP}{d\rho}}=\sqrt{\mathsf{A}\gamma}\rho^{\frac{\gamma-1}{2}}$ would exceed the speed of light $\mathsf{c}$ for large $\rho$.
In this sense, the remainder term $\Upsilon(\mathsf{A}\rho^{\gamma-1}/\mathsf{c}^2)$ is not introduced in order to consider coupling of two polytropes. One might assume the simple $\gamma$-law $P=\mathsf{A}\rho^{\gamma_0}, 1<\gamma_0<2$ without this remainder term near the vacuum $0<\rho \ll 1$ and smoothly match it with a equation of state like, e.g., $P=(\gamma_{\infty}-1)\mathsf{c}^2\rho$ with $1<\gamma_{\infty}<2$ in a region of large $\rho \gg 1$
in such a manner that $dP/d\rho <\mathsf{c}^2$ be satisfied throughout $0 \leq \rho <+\infty$. 
Under the assumption {\bf (A)}, particularly by \eqref{0104}, the equation of state tends to the simple $\gamma$-law $P=\mathsf{A}\rho^{\gamma}$ as $\mathsf{c} $ tends to $+\infty$, say, in the Newtonian limit. But this is supposed for the sake of simplicity of the discussion; namely the Newtonian limit could be supposed to be more generally admitting, e.g., the equation of state of the white dwarf (\cite[Chapter XI]{ChandraSS}): 
$$P=B_1[Q(2Q^2-3)\sqrt{1+Q^2}+3\sinh^{-1}Q],\quad \rho=B_2Q^3, $$
or $\rho =B_0 [u(1+u)]^{\frac{3}{2}}$ with $u=\int dP/\rho$,
$B_1,B_2, B_0$ being positive constants. 

In Section 2, where we shall discuss the reduction of the Einstein-Euler equations to a system of equations on the coefficients of the axisymmetric metric and state variables and its equivalence, the consistency, we need not the detailed behavior ruled by \eqref{0104} of $P$ as $\rho \rightarrow +0$; It is sufficient to suppose that $P$ is a $C^1$- function of $\rho \geq 0$. But, in Section 3, where we construct asymptotically flat metric generated by rotating fluid, the behavior ruled by \eqref{0104} will play a crucial r\^{o}le. Actually, to begin with, we need it to start with a rotating configuration in the Newtonian limit under $\frac{6}{5} < \gamma <2$, where $\gamma >\frac{6}{5}$ is necessary in order that the configuration have a finite diameter. Moreover during our post-Newtonian approximation procedure, the behavior ruled by \eqref{0104} will be often used. The restriction $\gamma <2$ is
mobilized in order that $\rho$ vanishes in $C^1$-way at the vacuum boundary
$\Sigma=\partial\{ \rho >0\}$, like
$\rho \sim C\Big[\mathrm{dist}(\cdot, \Sigma)\Big]^{\frac{1}{\gamma-1}}$,
with $\frac{1}{\gamma-1}>1$.
The analyticity near $0$ of the function $\Upsilon$ which vanishes at $0$ in the remainder term is supposed for the sake of simplicity of discussion. This restriction  might be loosen to the $C^{2,\alpha}$-class regularity, $\alpha$ being a positive number $<1$. But more rough regularities cause many troubles in the discussion in Section 3.
\\

In this article we consider axially symmetric metrics, that is, taking
the co-ordinates
\begin{equation}
x^0=\mathsf{c}t,\quad x^1=\varpi,\quad x^2=\phi,\quad x^3=z,
\end{equation}
we consider the metric $ds^2=g_{\mu\nu}dx^{\mu}dx^{\nu}$ in the following form (Lewis 1932 \cite{Lewis}, 
Papapetrou 1966 \cite{Papapetrou} ), say, ` Lanczos form' after \cite{Lanczos}: 
\begin{equation}
ds^2=e^{2F}
(\mathsf{c}dt+Ad\phi)^2-
e^{-2F}[e^{2K}
(d\varpi^2+dz^2)+\Pi^2d\phi^2], \label{2}
\end{equation}
that is,
\begin{align*}
& g_{00}=e^{2F}, \quad g_{02}=g_{20}=e^{2F}A, \quad
g_{11}=g_{33}=-e^{2(-F+K)}, \\
&g_{22}=e^{2F}A^2-e^{-2F}\Pi^2,\quad \mbox{other}\  g_{\mu\nu}=0,
\end{align*}
where the quantities $F, A, K$ and $\Pi$, which will be called `potentials',  are functions of only  $\varpi$ and $z$. \\

The 4-velocity vector filed $U^{\mu}$ is supposed to be of the form
\begin{equation}
U^{\mu}\frac{\partial}{\partial x^{\mu}}=e^{-G}\Big(\frac{1}{\mathsf{c}}\frac{\partial}{\partial t}+
\frac{\Omega}{\mathsf{c}}\frac{\partial}{\partial\phi}\Big),
\end{equation}
that is,
\begin{equation}
U^0=e^{-G},\quad U^1=U^3=0,\quad
U^2=e^{-G}\frac{\Omega}{\mathsf{c}},
\end{equation}
where the angular velocity $\Omega$ depends only on $\varpi, z$.

Since $U^{\mu}U_{\mu}=1$, the factor $e^{-G}=U^0$ is given by
$$ e^{-2G}\Big(g_{00}+2g_{02}\frac{\Omega}{\mathsf{c}}+g_{22}\frac{\Omega^2}{\mathsf{c}^2}\Big)=1, $$
that is,
\begin{equation}
e^{2G}=e^{2F}\Big(1+\frac{\Omega}{\mathsf{c}}A\Big)^2
-e^{-2F}\frac{\Omega^2}{\mathsf{c}^2}\Pi^2. \label{18}
\end{equation}

In this point of view, $F, A, \Pi$ to be considered should satisfy the following assumption 
{\bf (B) } with respect to the prescribed $\Omega$ :\\

{\bf (B)}: {\it It holds that}
\begin{equation}
e^{2F}\Big(1+\frac{\Omega}{\mathsf{c}}A\Big)^2-e^{-2F}
\frac{\Omega^2}{\mathsf{c}^2}\Pi^2 >0.
\end{equation}\\

We want to find potentials $F, A, \Pi, K$ which give an asymptotically flat metric, namely, we use the following

\begin{Definition}
The metric \eqref{2} on $r:=\sqrt{\varpi^2+z^2} \gg 1$ is said to be asymptotically flat if
\begin{subequations}
\begin{align}
&F=-\frac{\mathsf{G}M}{\mathsf{c}^2r}+O\Big(\frac{1}{r^2}\Big), \label{7.2a}\\
&A=\frac{2\mathsf{G}J\varpi^2}{\mathsf{c}^3r^3}
+O\Big(\frac{\varpi^2}{r^4}\Big),
\label{7.2b}\\
&\Pi=\varpi\Big(1+O\Big(\frac{1}{r^2}\Big)\Big), \label{7.2c}\\
&e^K=1+O\Big(\frac{1}{r^2}\Big) \label{7.2d}
\end{align}
\end{subequations}
as $r \rightarrow +\infty$. Here $M, J$ are constants.  
\end{Definition}
As for the physical reasoning of this asymptotic conditions, see \cite[p.19]{Islam} and \cite[Chapter 19]{MisnerTW}.

Actually the Kerr metric is asymptotically flat. The Kerr metric described by the Boyer-Lindquist Schwarzschild type co-ordinate system $(t, \underline{r}, \underline{\vartheta}, \phi)$ reads
\begin{align}
ds^2=&\mathsf{c}^2dt^2-\mathsf{\Sigma}\Big(\frac{d\underline{r}^2}{\mathsf{\Delta}}+d\underline{\vartheta}^2\Big)
-(\underline{r}^2+a^2)\sin^2\underline{\vartheta} d\phi^2 + \nonumber \\
&-\frac{2\mathsf{m} \underline{r}}{\mathsf{\Sigma}}(a\sin^2\underline{\vartheta}d\phi-\mathsf{c}dt)^2, \label{Kerrmetric}
\end{align}
where
\begin{equation}
\mathsf{\Delta}=\underline{r}^2-2\mathsf{m}\underline{r}+a^2,\qquad
\mathsf{\Sigma}=\underline{r}^2+a^2\cos^2\underline{\vartheta}
\end{equation}
with parameters $\mathsf{m}, a$ such that $0\leq |a|\leq \mathsf{m}$. See 
\cite[(2.8), (2.12), (2.13)]{BoyerL} or \cite[p.306, (217)]{ChandraBH}. We consider the metric in the exterior vacuum region $\{ \underline{r} >2\mathsf{m} \}$. (Note that $2\mathsf{m} \geq \underline{r}_+:=\mathsf{m}+\sqrt{\mathsf{m}^2-a^2}$ so that $\mathsf{\Delta} >0$ when $\underline{r} >2\mathsf{m}$. ) Then the change of co-ordinates
$(\underline{r}, \underline{\vartheta}) \mapsto (\varpi, z)$ defined by
\begin{equation}
\varpi=\sqrt{\mathsf{\Delta}}\sin\underline{\vartheta},\qquad z=(\underline{r}-\mathsf{m})\cos\underline{\vartheta} 
\end{equation}
brings the metric \eqref{Kerrmetric} to the Lanczos form \eqref{2} with
the potential $F, A, \Pi, K$ determined by
\begin{align*}
&e^{2F}=1-\frac{2\mathsf{m}\underline{r}}{\mathsf{\Sigma}}, \\
&e^{2F}A=\frac{2\mathsf{m}\underline{r}}{\mathsf{\Sigma}}a\sin^2\underline{\vartheta}, \\
&e^{-2F}\Pi^2=e^{2F}A^2+
(\underline{r}^2+a^2)\sin^2\underline{\vartheta}+\frac{2\mathsf{m}\underline{r}}{\mathsf{\Sigma}}a^2\sin^4\underline{\vartheta}, \\
&e^{-2F}e^{2K}=\frac{\mathsf{\Sigma}}{\mathsf{\Delta}+(\mathsf{m}^2-a^2)\sin^2\underline{\vartheta}}.
\end{align*}
Note that $\displaystyle 1-\frac{2\mathsf{m}\underline{r}}{\mathsf{\Sigma}}>0$ when $\underline{r} >2\mathsf{m}$. Since
$$r=\underline{r}\sqrt{\Big(1-\frac{\mathsf{m}}{\underline{r}}\Big)^2-\frac{\mathsf{m}^2-a^2}{\underline{r}^2}\sin^2\underline{\vartheta}}
=\underline{r}\Big(1-\frac{\mathsf{m}}{\underline{r}}+O\Big(\frac{1}{\underline{r}^2}\Big)\Big)
$$
so that $$\underline{r}=r\Big(1+\frac{\mathsf{m}}{r}+O\Big(\frac{1}{r^2}\Big)\Big) $$
as $r \rightarrow +\infty$, the conditions \eqref{7.2a} $\sim$ \eqref{7.2d} of the asymptotic flatness are satisfied with $M, J$ determined by
 $\displaystyle \mathsf{m}=\frac{\mathsf{G}M}{\mathsf{c}^2}, a=\frac{J}{\mathsf{c}M}$. \\

In Section 2, we drive the system of equations for the potentials $F, A, \Pi, K$ with the density distribution $\rho$ and the angular velocity distribution $\Omega$. 

Although the system for the adjoint potential $F', A', \Pi'(=\Pi), K'$ described in a co-rotating co-ordinate system with constant $\Omega$ is given in \cite{Meinel} as 
\cite[(1.34a), (1.34b), (1.34c), (1.35a), (1.35b)]{Meinel} without proof, we here want to derive the system in the ZAMO co-ordinate system. The reason is that it seems difficult to discuss the metric on the whole space on the setting of \cite{Meinel} as discussed in \cite[VII Discussion]{asEE}. 

The derivation can be done along the calculation way described in \cite{Islam} on the vacuum, say, $\rho=0$. However the recoverability of the full set of Einstein equations is not self-evident. Actually proof of the recoverability requires a tedious and tricky calculation. It is done under the assumption that either $\Omega$ is constant or $\rho=0$ in the considered domain. The reversibility is doubtful if differential rotations with variable $\Omega$ on the support of $\rho$ are concerned.  Anyway we have the equivalence result, Theorem \ref{Th.1}.

The consistency of the first order system of equations for $K$ also needs a careful treatment. Actually it happens to be said that we need not write down the second equation, say for 
$\partial K/\partial z$, because it provides no new information than the first equation, say, for $\partial K/\partial \varpi$ . See \cite[p.202, passage after (4d)]{ButterworthI}.
One considers that there are no problem, since only simple quadrature is sufficient. This is a hasty conclusion. It is not the case when the vacuum is not presupposed. Through a tedious calculation we prove a useful condition for the consistency, Theorem \ref{ThN2}.

Section 3 is devoted to construction of an asymptotically flat metric generated by a slowly rotating compactly supported perfect fluid. As in \cite{asEE} we work on the post-Newtonian approximation using the model given by the distorted Lane-Emden function established in
\cite{JJTM1}, \cite{JJTM2}. The angular velocity is supposed to be constant on a bounded domain including the support of the density as the study developed in \cite{asEE}, but, in order to get a global potentials, we have to cut off the angular velocity outside the bounded domain. This is another  reason why we derive the system of equations in the ZAMO co-ordinate system with variable $\Omega$ in Section 2. The co-rotating co-ordinate system of \cite{Meinel},\cite{asEE} may not work well in the infinitely far away region.  On the bounded domain an almost same discussion as that of \cite{asEE} can be mobilized, but we should prepare a set of functional spaces of functions defined on the whole space, say, $r=\sqrt{\varpi^2+z^2}<\infty$. Of course on the exterior domain the vacuum is expected. This framework is made by using the Kelvin transformation well known in the classical potential theory. The results are Theorems \ref{Th.WYX}, \ref{Th.K}.
Summing up, we declare the final goal Theorem \ref{Th.5}.

In Appendix there will be listed up constants, parameters and variables appearing frequently in this article as a glossary. \\

Here let us describe comments on the significance of the results of the present research in context of some preceding research achievements.  

(1) The recent important achievement in the matter-vacuum matching problem is \cite{Hernandez2017}
by J. L. Hernandez-Pastora and L. Herrera published in 2017. This paper provides axially symmetric metric smoothly matched on the boundary surface of the fluid distribution to the Kerr metric. But the 
energy-momentum tensor of the interior solution constructed there is anisotropic, that, is, the pressure tensor has the non-zero anisotropic components $P_{xx}, P_{zz}$, while $P_{xx}\not=P_{zz}, P_{xy}=P_{yy}-P_{zz}=0$, so that $T_1^1 \not= T_3^3$ and so on, where 
$$ T_{\mu\nu}=(\mathsf{c}^2\rho+P)U_{\mu}U_{\nu}-Pg_{\mu\nu}-\mathsf{P}_{\mu\nu} $$
with not all $\mathsf{P}_{\mu\nu}$ vanishing.
On the other hand, we are considering isotropic fluid state in this article. This research \cite{Hernandez2017} has preceding achievements \cite{Hernandez2016} and \cite{Herrera2013}
which has the root of way of research at `the pioneering work by Hernandez Jr \cite{Hernandez1967}, where a general method for obtaining solutions describing axially symmetric sources is presented; Such a method, or some of its modifications were used to find sources of different Weyl space-times.' (quoted from \cite[Page 024014-1 Right, Line 3-, ]{Herrera2013}, or \cite[Page 1 Left, Line 7 from the bottom -]{Hernandez2016}) Here `Weyle metric' is that of the form
$$ds^2=e^{2F}(\mathsf{c}dt+Ad\phi)^2 -e^{-2F}[e^{2K}(d\varpi^2+dz^2)
+\varpi^2d\phi^2], $$
which is \eqref{2} with $\Pi=\varpi$ provided that $P=0$. The Kerr metric turns out to be, outside the ergo space, of this form described by the canonical co-ordinate system. (See 
\cite[III. A.]{Hernandez2017}.) Actually the strategy of Hernandez Jr reads

\begin{quote}
In order to obtain more insight into the Kerr metric and also the more general problem of rotation in general relativity we have indicated here a simple method for constructing interior solutions which might serve as sources for the Kerr metric. These solutions are not obtained in the conventional manner of first choosing an equation of state and then solving the field equations. Instead the interior metrics are obtained by a method of guessing (g method) [ Footnote: \cite[Chapter VIII, \S 1]{Synge} ] which can be used for finding sources of any metric, which indicates that it possesses a positive mass sources. The essential idea of the g method is that one guesses an ``interior'' metric $g_{\mu\nu}$ and then calculates the resulting stress-energy tensor
$T_{\mu\nu}$ using the Einstein equations. ( \cite[Page 170 Left, Line 9 -]{Hernandez1967b})
\end{quote}

We see that the achievements \cite{Herrera2013}, \cite{Hernandez2016}, \cite{Hernandez2017} have been done along this strategy established by Hernandez Jr. Clearly the strategy of this article is of the opposite direction, namely, firstly given the equation of state of barotropic isotropic gas, we are trying to solve the Einstein equations to construct an asymptotically flat metric generated by the rotating gaseous mass obeying the given equation of state. However the author is not yet able to claim that the so constructed exterior metric would admit a suitable co-ordinate transformation by which the exterior metric turns out to
be the Kerr, or the canonical Weyl metric. In this sense a further achievement along the Hernandez Jr's strategy may be expected so that it will cover physically meaningful isotropic and non-homogeneous fluid sources as that of the interior metric matched to the exterior Kerr, or Weyl metric on the boundary surface. Regrettably this task is beyond the scope of the present ability of the author. Of course the anisotropic source matter provided by J. L. Hernandez-Pastora and L. Herrera in \cite{Hernandez2017} is physically meaningful, that is, all physical variables are regular within the matter distribution, the  density is positive, and so on, and, relationally, from the mathematical point of view, we should refer the theory on relativistic elasticity developed in \cite{Beig}, \cite{Andersson2008}, \cite{Andersson2016} by L. Andersson, R. Beig, T. Oliynyk, and B. G. Schmidt. 
Additionally let us note that \cite{Herrera2013} constructed isotropic but incompressible, or, homogeneous, spheroidal interior solutions, besides anisotropic interior solutions. As for this solution, the pressures are isotropic in the interior $r < r_{\Sigma}$ and vanishes on the boundary $r=r_{\Sigma}$, but the density $\rho$ is homogeneous, that is, $\rho=\rho_{\mathsf{O}}$, being a positive constant, throughout everywhere on $r <r_{\Sigma}$. It contrasts to our situation in which $\rho$ is variable and vanishes in $C^1$-way on the vacuum boundary. Remarkably \cite{Herrera2013} proved that this homogeneous isotropic interior solution cannot be matched to any Weyl exterior, even though it has the boundary surface of vanishing pressure. (\cite[Page 024014-6 Right, Line 16]{Herrera2013} )

(2) W. Israel \cite{Israel} clarified the unique source of the causally maximal extention of the Kerr's space-time. According to the result, it is a disk which has a negative surface density, rotates with supra-luminal speed and angular velocity, and has positively infinite mass and angular momentum. The motivation of the research reads

\begin{quote}
Several papers have been devoted to the problem of finding the Kerr exterior field to rotating material sources of various kinds. Naturally, this problem does not have a unique solution, since one is at liberty (for example) to choose arbitrarily the boundary between the exterior vacuum and the source. The arbitrariness disappears if one ask about the source of the {\it maximally extended} vacuum metric, but here we are thwarted by our
present inability to interpret singularities of Einstein's field equations. Newman and Janis conjectured that Kerr's solution is the field of spinning ring of mass, but this fails to take account of the peculiar geometry of the disk. (\cite[Page 642 Right, Line 13 - ]{Israel})
\end{quote}

Of course the matter in question and results are quite different from those of our present research, for we are not considering the maximal extension of the Kerr metric as a vacuum metric. However there is suggested a task of the further research in the future. Namely our results cover slowly rotating gaseous sources whose configurations are near to spheres. Sources with configurations near thin disks should be considered, too. Actually, if we consider rotating homogeneous self-gravitating liquid masses in the Newtonian non-relativistic framework, then, for any small angular velocity, there appear thin disk solutions besides the well known ellipsoidal solutions slightly distorted. Thus we have a task which should be done in the future.

(3) The survey \cite{Krasinski1978} by A. Krasi\'{n}ski, 1978, gave the definition of the `ellipsoidal space-time' and expressed the opinion that one should continue the effort to find an ellipsoidal interior solution of material source for the Kerr metric, rejecting many erroneous negative criticisms. According this paper, an ellipsoidal space-time is such that, by a suitable co-ordinate system $(r,\vartheta, \phi)$, the metric takes the form
\begin{align*}
ds^2=&\Big(\frac{1-kV}{U}dt+kd\phi\Big)^2-
f^2dr^2 \\
&-(r^2+a^2\cos^2\vartheta)d\vartheta^2
-(r^2+a^2)\sin^2\vartheta\Big(d\phi-\frac{U}{V}dt\Big)^2.
\quad\mbox{\cite[(3.14)]{Krasinski1978}} 
\end{align*}
Here $U, V, k, f$ are functions of $(r,\vartheta)$, and, although not explicitly written,
$a$ should be a constant and $U, V$ may be related to the 4-velocity of the matter as
$$U^{\mu}\frac{\partial}{\partial x^{\mu}}=U\frac{\partial}{\partial t}+
V\frac{\partial}{\partial \phi}, $$
(See \cite[(2.11)]{Krasinski1978}), the quoter is not sure. Anyway the author T. M. of the present article does not know whether the metric constructed in the present article is `spheroidal' in this sense of \cite{Krasinski1978} or not. The author feels that it is not the case. The reason is as following. 

A. Krasi\'{n}ski wrote:

\begin{quote}
If the observer resides inside the body and wants to determine the shape of equipressure surfaces, then a third force, the hydrostatic pressure of layers above him, appears, but the problem is still that of statics. Now, statics pursued in general relativistic language and in the rest frame of the static system is just identical to Newtonian statics. We can therefore reasonably expect that, in the local rest frame of an observer comoving with matter, what was an ellipsoid in Newtonian theory will stay an ellipsoid in general relativity. (\cite[Page 32, Line 7 from the bottom -]{Krasinski1978})
\end{quote}

\noindent Let us note that A. Krasi\'{n}ski keeps in mind only rotating liquid mass as Newtonian correspondences 
with Kerr interior, namely, he wrote:

\begin{quote}
In Newtonian hydrodynamics one is able to gain any exact information about a finite rotating portion of fluid only in case it has the shape of an ellipsoid [R. A. Littleton: The Stability of Rotating Liquid Masses, S. Chandrasekhar: Ellipsoidal Figures of Equilibrium ]. Therefore, it will be natural to start analogous investigations in general relativity considering the same simplest case. It will be most convenient to use the ellipsoids as defining a coordinate system in the space-time.
(\cite[Page 30, Line 2 from the bottom -]{Krasinski1978})
\end{quote}

\noindent However, in the Newtonian theory, the situation is different for gaseous masses than liquid masses, while the books by R. A. Littleton and by S. Chandrasekhar
deal with liquid (incompressible, homogeneous ) fluid masses under self-gravitation. The following ``no-go theorem''
has been known for more than 120 years:\\

{\it Theorem of Hamy-Pizzetti: An ellipsoidal stratification} ( that is, the situation that all level surfaces 
$\{\rho =\mbox{Const.}\}$ are ellipsoids ) {\it  is impossible for heterogeneous } (
that is, with non-constant $\rho$ ), {\it  rotating symmetric figures of equilibrium} ( that is, axisymmetric stationary solutions of the Euler-Poisson equations rotating with non-zero angular velocity). (\cite[Sec. 3.2]{Moritz}, \cite{Pizzetti}, \cite{Wavre})\\

\noindent  By this no-go theorem it is impossible that all level surfaces $\{\rho =\mbox{Const.}\}$ of the rotating polytrope to be considered are ellipsoids, provided that the angular velocity $\Omega\not=0$.
Although this no-go theorem does not claim that it is impossible that the vacuum boundary, an individual level surface, $\partial\{ \rho >0\}$ is an ellipsoid, it seems to be the case. This is the situation in the Newtonian theory. So, the author of the present article feels that it is hard to expect that the metric of the rotating gaseous mass constructed in general relativity here could be an `ellipsoidal' metric in the sense of \cite{Krasinski1978}. Of course this guess is not an objection to the strategy whose successfulness  A. Krasi\'{n}ski believed. 
Actually the author does not know a systematic method to judge whether the given metric to be considered is `spheroidal' or not without trying all possible co-ordinate systems he can lay his hands on.

\section{The system of equations for the potentials and density distribution}

\subsection{Einstein-Euler equations}

The non-zero components of the energy-momentum tensor $T^{\mu\nu}$ are:
\begin{align*}
&T^{00}=e^{-2G}(\epsilon+P)+
\frac{P}{\Pi^2}(e^{2F}A^2-e^{-2F}\Pi^2), \\
&T^{02}=T^{20}=
e^{-2G}(\epsilon+P)\frac{\Omega}{\mathsf{c}}-
\frac{P}{\Pi^2}e^{2F}A, \\
&T^{11}=T^{33}=Pe^{2F-2K}, \\
&T^{22}=
e^{-2G}(\epsilon+P)\frac{\Omega^2}{\mathsf{c}^2}+
\frac{P}{\Pi^2}e^{2F}.
\end{align*}\\

Recall that the Christoffel symbols are defined by
$$ \Gamma_{\nu\lambda}^{\mu}=\frac{1}{2}g^{\mu\alpha}
(\partial_{\lambda}g_{\alpha\nu}+
\partial_{\nu}g_{\alpha\lambda}-\partial_{\alpha}g_{\nu\lambda} ).$$
and hereafter $\partial_{\mu}$ stands for 
$\displaystyle \frac{\partial}{\partial x^{\mu}}$, while $\nabla_{\mu}$ will stand for
the covariant derivative
$$\nabla_{\lambda}A^{\mu\nu}=
\partial_{\lambda}A^{\mu\nu}
+\Gamma_{\alpha\lambda}^{\mu}A^{\alpha\nu}+
\Gamma_{\alpha\lambda}^{\nu}A^{\mu\alpha}.
$$\\

The Euler equations are $\nabla_{\mu}T^{\mu\nu}=0$. 
Let us consider them.

Keeping in mind \eqref{18}, we can show through tedious calculations that 
$$\nabla_{\mu}T^{\mu 1}=
e^{2F-2K}\Big(\partial_1P+(\epsilon+P)\partial_1G
-(\epsilon+P)U^0U_2\frac{\partial_1\Omega}{\mathsf{c}}\Big) $$
and
$$\nabla_{\mu}T^{\mu 3}=
e^{2F-2K}\Big(\partial_3P+(\epsilon+P)\partial_3G
-(\epsilon+P)U^0U_2\frac{\partial_3\Omega}{\mathsf{c}}\Big). $$
Here
\begin{align}
U^0U_2&=e^{2F-2G}\Big(A\Big(1+\frac{\Omega}{\mathsf{c}}\Big)-e^{-4F}\frac{\Omega}{\mathsf{c}}\Pi^2\Big)  \nonumber \\
&=\Big(
\Big(1+\frac{\Omega}{\mathsf{c}}A\Big)^2
-e^{-4F}\frac{\Omega^2}{\mathsf{c}^2}\Pi^2\Big)^{-1}
\Big(A\Big(1+\frac{\Omega}{\mathsf{c}}\Big)-e^{-4F}\frac{\Omega}{\mathsf{c}}\Pi^2\Big).
\end{align}
Hence the identity $\nabla_{\mu}T^{\mu 1}=\nabla_{\mu}T^{\mu 3}=0$ reduces to
\begin{align}
&\partial_1P+(\epsilon+P)\partial_1G
-(\epsilon+P)U^0U_2\frac{\partial_1\Omega}{\mathsf{c}} 
=0,  \nonumber \\
&\partial_3P+(\epsilon+P)\partial_3G
-(\epsilon+P)U^0U_2\frac{\partial_3\Omega}{\mathsf{c}}=0. \label{2.7}
\end{align}
On the other hand we see that the other Euler equations
$$\nabla_{\mu}T^{\mu 0}=0,\qquad \nabla_{\mu}T^{\mu 2}=0 $$
hold automatically, for we are assuming $\partial_0\Omega=\partial_2\Omega=0$.\\

Therefore, defining the `relativistic enthalpy density' $u$ by
\begin{equation}
u=f^{u}(\rho):=\mathsf{c}^2\int_0^{\rho}\frac{dP}{\epsilon+P}=
\int_0^{\rho}\frac{dP}{\rho+P/\mathsf{c}^2}, \label{Def.u}
\end{equation}
we have
\begin{equation}
\frac{u}{\mathsf{c}^2}+G=\mbox{Const.},
\end{equation}
while $\rho >0$, provided that $\Omega$ is constant while $\rho>0$. \\

Later we shall use the following

\begin{Definition} \label{Def.frho}
Defining $u=f^u(\rho)$ by \eqref{Def.u} for $\rho>0$, we define $f^{\rho}$ on $\mathbb{R}$ of the form
$$\rho=f^{\rho}(u)= 
\Big(\frac{\gamma-1}{\mathsf{A}\gamma}\Big)^{\frac{1}{\gamma-1}}
(u\vee 0)^{\frac{1}{\gamma-1}}(1+\Upsilon_{\rho}(u/\mathsf{c}^2))$$
with an analytic function $\Upsilon_{\rho}$ defined near $0$ such that
$\Upsilon_{\rho}(0)=0$. Here $u \vee 0$ stands for $\max\{ u, 0\}$.
We rewrite
\begin{align}
&f^{\rho}(u)=f_{\mathsf{N}}^{\rho}(u)(1+\Upsilon_{\rho}(u/\mathsf{c}^2)), \\
&f_{\mathsf{N}}^{\rho}(u)=\Big(\frac{\gamma-1}{\mathsf{A}\gamma}\Big)^{\frac{1}{\gamma-1}}
(u\vee 0)^{\frac{1}{\gamma-1}}.
\end{align}
Also we denote
\begin{align}
&f^P(u)=f_{\mathsf{N}}^P(u)(1+\Upsilon_P(u/\mathsf{c}^2)), \\
&f_{\mathsf{N}}^P(u)=\mathsf{A}\Big(\frac{\gamma-1}{\mathsf{A}\gamma}\Big)^{\frac{\gamma}{\gamma-1}}
(u\vee 0)^{\frac{\gamma}{\gamma-1}},
\end{align}
where $\Upsilon_P$ is an analytic function defined near $0$ such that $\Upsilon_P(0)=0$.
\end{Definition}

Hereafter we shall use

\begin{Notation}
We denote
\begin{equation}
Q \vee Q'=\max\{ Q, Q'\}, \qquad
Q\wedge Q'=\min\{ Q, Q'\}.
\end{equation}
\end{Notation}

In order to write down the Einstein equations
\begin{equation}
R_{\mu\nu}-\frac{1}{2}g_{\mu\nu}R=
\frac{8\pi\mathsf{G}}{\mathsf{c}^4}T_{\mu\nu}
\end{equation}
or
\begin{equation}
R_{\mu\nu}=\frac{8\pi\mathsf{G}}{\mathsf{c}^4}
(T_{\mu\nu}-\frac{1}{2}g_{\mu\nu}T),
\end{equation}
where $R$ stands for $g^{\alpha\beta}R_{\alpha\beta}$ and
$T$ stands for $g_{\alpha\beta}T^{\alpha\beta}$, let us compute $U_{\mu}$,
$T_{\mu\nu}$ and $T$. The result is as following:

\begin{align*}
&U_0=e^{2F-G}\Big(1+\frac{\Omega}{\mathsf{c}}A\Big), \\
&U_1=U_3=0, \\
&U_2=e^{-G+2F}\Big(
A\Big(1+\frac{\Omega}{\mathsf{c}}A\Big)-e^{-4F}\frac{\Omega}{\mathsf{c}}\Pi^2\Big) ;
\end{align*}
\begin{align*}
T_{00}&=(\epsilon+P)e^{4F-2G}\Big(1+\frac{\Omega}{\mathsf{c}}A\Big)^2
-Pe^{2F}, \\
T_{02}&=T_{20}= \\
&=(\epsilon+P)e^{4F-2G}\Big(1+\frac{\Omega}{\mathsf{c}}A\Big)
\Big(A\Big(1+\frac{\Omega}{\mathsf{c}}A\Big)-
e^{-4F}\frac{\Omega}{\mathsf{c}}\Pi^2\Big)
-Pe^{2F}A, \\
T_{11}&=T_{33}=Pe^{-2F+2K}, \\
T_{22}&=
(\epsilon+P)e^{-2G+4F}
\Big(A\Big(1+\frac{\Omega}{\mathsf{c}}A\Big)-e^{-4F}\frac{\Omega}{\mathsf{c}}\Pi^2\Big)^2
-P(e^{2F}A^2-e^{-2F}\Pi^2)
\end{align*}
and other $T_{\mu\nu}$'s are zero.

 Let us note that 
$$
T=(\epsilon+P)e^{-2G}\Big[e^{2F}\Big(1+\frac{\Omega}{\mathsf{c}}A\Big)^2
-e^{-2F}\frac{\Omega^2}{\mathsf{c}^2}\Pi^2\Big]
-4P
$$
reduces to
\begin{equation}
T=\epsilon-3P. \label{TeP}
\end{equation}\\

In order to calculate the Ricci tensor
$$ R_{\mu\nu}=\partial_{\alpha}\Gamma_{\mu\nu}^{\alpha}
-\partial_{\nu}\Gamma_{\mu\alpha}^{\alpha}
+\Gamma_{\mu\nu}^{\alpha}\Gamma_{\alpha\beta}^{\beta}
-\Gamma_{\mu\alpha}^{\beta}\Gamma_{\nu\beta}^{\alpha}
$$
 and to write down explicitly the Einstein equations, it is convenient to write the metric as
\begin{equation}
ds^2=f\mathsf{c}^2dt^2-2k\mathsf{c}dtd\phi-ld\phi^2
-e^m(d\varpi^2+dz^2), \label{25}
\end{equation}
that is, to put
\begin{equation}
f=e^{2F},\quad k=-e^{2F}A, \quad m=2(-F+K),\quad l=-e^{2F}A^2+
e^{-2F}\Pi^2.
\end{equation}

The expression \eqref{25} is called `Lewis metric' after \cite{Lewis}\\

First let us note the identity
\begin{equation}
\Pi^2=fl+k^2. \label{27}
\end{equation}
and that the factor $e^{-G}=U^0$ is given by
\begin{equation}
(U^0)^{-2}=e^{2G}=f-2\frac{\Omega}{\mathsf{c}}k
-\frac{\Omega^2}{\mathsf{c}^2}l.
\end{equation}\\

The components of the Ricci tensor are as following (other $R_{\mu\nu}$ are zero ):
\begin{subequations}
\begin{align}
\frac{2e^m}{\Pi}R_{00}&=\partial_1\frac{\partial_1f}{\Pi}+\partial_3\frac{\partial_3f}{\Pi}+
\frac{1}{\Pi^3}f\Sigma, \label{30a}\\
-\frac{2e^m}{\Pi}R_{02}&=-\frac{2e^m}{\Pi}R_{20}=
\partial_1\frac{\partial_1k}{\Pi}+\partial_3\frac{\partial_3k}{\Pi}+\frac{1}{\Pi^3}k\Sigma, \label{30b} \\
-\frac{2e^m}{\Pi}R_{22}&=
\partial_1\frac{\partial_1l}{\Pi}+\partial_3\frac{\partial_3l}{\Pi}+\frac{1}{\Pi^3}l\Sigma, \label{30c}\\
2R_{11}&=-\partial_1^2m-\partial_3^2m-2\frac{\partial_1^2\Pi}{\Pi}
+\frac{1}{\Pi}[(\partial_1m)(\partial_1\Pi)-(\partial_3m)(\partial_3\Pi)] + \nonumber \\
&+\frac{1}{\Pi^2}[(\partial_1f)(\partial_1l)+(\partial_1k)^2], \label{30d}\\
2R_{33}&=-\partial_1^2m-\partial_3^2m-2\frac{\partial_3^2\Pi}{\Pi}
-\frac{1}{\Pi}[(\partial_1m)(\partial_1\Pi)-(\partial_3m)(\partial_3\Pi)] + \nonumber \\
&+\frac{1}{\Pi^2}[(\partial_3f)(\partial_3l)+(\partial_3k)^2] \label{30e}\\
2R_{13}&=2R_{31}=-2\frac{\partial_1\partial_3\Pi}{\Pi}+
\frac{1}{\Pi}[(\partial_3m)(\partial_1\Pi)+(\partial_1m)(\partial_3\Pi)]
\nonumber \\
&+\frac{1}{2\Pi^2}[(\partial_1f)(\partial_3l)+(\partial_1l)(\partial_3f)
+2(\partial_1k)(\partial_3k)].\label{30f}
\end{align}
\end{subequations}

Here we have introduced the quantity $\Sigma$ defined by
\begin{equation}
\Sigma:=(\partial_1f)(\partial_1l)+(\partial_3f)(\partial_3l)+
(\partial_1k)^2+(\partial_3k)^2.
\end{equation}\\

Now the components of the 4-velocity vector are:
\begin{align*}
&U^0=e^{-G}=\Big(f-
2\frac{\Omega}{\mathsf{c}}k-\frac{\Omega^2}{\mathsf{c}^2}l\Big)^{-1/2},
\quad U^1=U^3=0,\quad U^2=\frac{\Omega}{\mathsf{c}}U^0; \\
&U_0=\Big(f-\frac{\Omega}{\mathsf{c}}k\Big)U^0,\quad
U_1=U_3=0,\quad
U_2=-\Big(k+\frac{\Omega}{\mathsf{c}}l\Big)U^0.
\end{align*}

The Einstein equations are
\begin{equation}
R_{\mu\nu}=\frac{8\pi\mathsf{G}}{\mathsf{c}^4}S_{\mu\nu}, \quad\mbox{where}\quad
S_{\mu\nu}:=T_{\mu\nu}-\frac{1}{2}g_{\mu\nu}T.
\end{equation}

The components $S_{\mu\nu}$ turn out to be as following (other
$S_{\mu\nu}$ are zero) :
\begin{subequations}
\begin{align}
&S_{00}=
\frac{1}{2}(\epsilon+P)e^{-2G}
\Big[\Big(f-\frac{\Omega}{\mathsf{c}}k\Big)^2+\frac{\Omega^2}{\mathsf{c}^2}\Pi^2\Big]
+Pf, \label{Ta}\\
&S_{02}=S_{20}=
\frac{1}{2}(\epsilon+P)e^{-2G}
\Big[-kf-2\frac{\Omega}{\mathsf{c}}fl+\frac{\Omega^2}{\mathsf{c}^2}kl\Big]-Pk, \label{Tb}\\
&S_{22}=
\frac{1}{2}(\epsilon+P)e^{-2G}\Big[\Pi^2+
\Big(k+\frac{\Omega}{\mathsf{c}}l\Big)^2\Big]-Pl, \label{Tc}\\
&S_{11}=S_{33}=\frac{e^m}{2}(\epsilon -P). \label{Td}
\end{align}
\end{subequations}

Recall \eqref{TeP}.\\

Thus the full set of Einstein equations is:
\begin{subequations}
\begin{align}
&R_{00}=\frac{8\pi\mathsf{G}}{\mathsf{c}^4}S_{00}, \label{36a}\\
&R_{02}=\frac{8\pi\mathsf{G}}{\mathsf{c}^4}S_{02}, \label{36b}\\
&R_{22}=\frac{8\pi\mathsf{G}}{\mathsf{c}^4}S_{22}, \label{36c}\\
&R_{11}=\frac{8\pi\mathsf{G}}{\mathsf{c}^4}S_{11}, \label{36d} \\
&R_{33}=\frac{8\pi\mathsf{G}}{\mathsf{c}^4}S_{33}, \label{36e}\\
&R_{13}=0 \label{36f}
\end{align}
\end{subequations}
. 

\subsection{Equations for the potentials $F, A, \Pi, K$}

{\bf 1) } From  \eqref{30a} and \eqref{36a} we can derive
\begin{align}
&\partial_1^2F+\partial_3^2F+\sum_{j=1,3}
\Big[\frac{1}{\Pi}(\partial_j\Pi)(\partial_jF)+\frac{e^{4F}}{2\Pi^2}(\partial_jA)^2\Big] = \nonumber \\
&=\frac{4\pi\mathsf{G}}{\mathsf{c}^4}
e^{2(-F+K)}
\Big[(\epsilon+P)
\frac{e^{2F}\Big(1+\frac{\Omega}{\mathsf{c}}A\Big)^2+e^{-2F}\frac{\Omega^2}{\mathsf{c}^2}\Pi^2}{e^{2F}\Big(1+\frac{\Omega}{\mathsf{c}}A\Big)^2-e^{-2F}\frac{\Omega^2}{\mathsf{c}^2}\Pi^2}
+2P\Big] \label{N.2}
\end{align}\\

{\bf 2)} From \eqref{30b}  and \eqref{36b}  we can derive 

\begin{align}
&\partial_1^2A+\partial_3^2A+\sum_{j=1,3}\Big[-\frac{1}{\Pi}(\partial_j\Pi)(\partial_jA)+
4(\partial_jF)(\partial_jA) + \nonumber \\
&+\Big(\frac{e^{4F}}{\Pi^2}(\partial_jA)^2
+2\partial_j^2F+2\frac{1}{\Pi}(\partial_j\Pi)(\partial_jF)\Big)\cdot A\Big] = \nonumber \\
&=\frac{8\pi\mathsf{G}}{\mathsf{c}^4}
e^{2(-F+K)}\Big[(\epsilon+P)
\Big(
A-2\frac{\Omega}{\mathsf{c}}\frac{e^{-2F}\Pi^2}{e^{2F}\Big(1+\frac{\Omega}{\mathsf{c}}A\Big)^2-e^{-2F}\frac{\Omega^2}{\mathsf{c}^2}\Pi^2}
\Big)+2AP\Big]. \label{N.4}
\end{align}

Taking \eqref{N.4}$- 2A\times$\eqref{N.2}, we get
\begin{align}
&\frac{2e^m}{f}(R_{02}-2A\cdot R_{00})= \nonumber \\
&=\partial_1^2A+\partial_3^2A+\sum_{j=1,3}\Big[-\frac{1}{\Pi}(\partial_j\Pi)(\partial_jA)+
4(\partial_jF)(\partial_jA) \Big]= \nonumber \\
&=-\frac{16\pi\mathsf{G}}{\mathsf{c}^4}
e^{2(-F+K)}(\epsilon+P)
\frac{e^{-2F}\frac{\Omega}{\mathsf{c}}\Pi^2\Big(1+\frac{\Omega}{\mathsf{c}}A\Big)}{e^{2F}\Big(1+\frac{\Omega}{\mathsf{c}}A\Big)^2-e^{-2F}\frac{\Omega^2}{\mathsf{c}^2}\Pi^2}. \label
{N.5}
\end{align}\\

{\bf 3)} We have the identity
\begin{equation}
lS_{00}-2kS_{02}
-fS_{22}=2P\Pi^2, \label{N.6}
\end{equation}
which can be verified from \eqref{Ta}\eqref{Tb}\eqref{Tc} thanks to \eqref{27}. 
On the other hand,
we have the identity
\begin{equation}
\frac{e^m}{\Pi}(lR_{00}-2kR_{02}-fR_{22})=
\partial_1^2\Pi+\partial_3^2\Pi,\label{N.7}
\end{equation}
which can be verified from \eqref{30a}\eqref{30b}\eqref{30c} thanks to \eqref{27}.
Hence \eqref{N.6} and \eqref{N.7} leads us to the equation
\begin{equation}
\partial_1^2\Pi+\partial_3^2\Pi=\frac{16\pi\mathsf{G}}{\mathsf{c}^4}
e^{2(-F+K)}P\Pi, \label{N.8}
\end{equation}
if \eqref{36a}\eqref{36b}\eqref{36c} hold. 

\begin{Remark} If we consider the vacuum or the dust for which $P=0$, then \eqref{N.8} says that $\Pi(\varpi, z)$ is a harmonic function and we can assume $\Pi=\varpi$ by conformal change of co-ordinates.
See \cite[p. 26]{Islam}. (This was first used by 
\cite{Weyl} and generalized to the present case by
\cite{Lewis}.)  But it is not the case when $P\not=0$. 
\end{Remark}

Summing up, we have the following

\begin{Proposition}
The set of equations \eqref{36a} \eqref{36b} \eqref{36c}
implies that
\begin{subequations}
\begin{align}
&\frac{\partial^2F}{\partial\varpi^2}+\frac{\partial^2F}{\partial z^2}
+\frac{1}{\Pi}\Big(\frac{\partial \Pi}{\partial\varpi}\frac{\partial F}{\partial\varpi}
+\frac{\partial \Pi}{\partial z}\frac{\partial F}{\partial z}\Big)+
\frac{e^{4F}}{2\Pi^2}\Big[\Big(\frac{\partial A}{\partial\varpi}\Big)^2+
\Big(\frac{\partial A}{\partial z}\Big)^2\Big] = \nonumber \\
&=\frac{4\pi\mathsf{G}}{\mathsf{c}^4}
e^{2(-F+K)}
\Big[(\epsilon+P)
\frac{e^{2F}\Big(1+\frac{\Omega}{\mathsf{c}}A\Big)^2+e^{-2F}\frac{\Omega^2}{\mathsf{c}^2}\Pi^2}{e^{2F}\Big(1+\frac{\Omega}{\mathsf{c}}A\Big)^2-e^{-2F}\frac{\Omega^2}{\mathsf{c}^2}\Pi^2}
+2P\Big], \label{N.9a} \\
&\frac{\partial^2A}{\partial\varpi^2}+\frac{\partial^2A}{\partial z^2}
-\frac{1}{\Pi}\Big(\frac{\partial\Pi}{\partial\varpi}\frac{\partial A}{\partial\varpi}+\frac{\partial\Pi}{\partial z}\frac{\partial A}{\partial z}\Big)+
4\Big(\frac{\partial F}{\partial\varpi}\frac{\partial A}{\partial\varpi}+
\frac{\partial F}{\partial z}\frac{\partial A}{\partial z}\Big) = \nonumber \\
&=-\frac{16\pi\mathsf{G}}{\mathsf{c}^4}
e^{2(-F+K)}(\epsilon+P)
\frac{e^{-2F}\frac{\Omega}{\mathsf{c}}\Pi^2\Big(1+\frac{\Omega}{\mathsf{c}}A\Big)}{e^{2F}\Big(1+\frac{\Omega}{\mathsf{c}}A\Big)^2-e^{-2F}\frac{\Omega^2}{\mathsf{c}^2}\Pi^2}, \label{N.9b} \\
&\frac{\partial^2\Pi}{\partial \varpi^2}
+\frac{\partial^2\Pi}{\partial z^2}=\frac{16\pi\mathsf{G}}{\mathsf{c}^4}
e^{2(-F+K)}P\Pi. \label{N.9c}
\end{align}
\end{subequations}
\end{Proposition}

Inversely we can show that \eqref{N.9a}\eqref{N.9b}\eqref{N.9c} imply 
\eqref{36a}\eqref{36b}\eqref{36c}.
In fact, if we put
$$Q_{\mu\nu}:=R_{\mu\nu}-\frac{8\pi\mathsf{G}}{\mathsf{c}^4}S_{\mu\nu},
$$
then \eqref{36a}\eqref{36b}\eqref{36c} claim that $Q_{00}=Q_{02}=Q_{22}=0$.
On the other hand \eqref{N.9a}\eqref{N.9b}\eqref{N.9c} are nothing but
\begin{align*}
&Q_{00}=0, \\
&-2AQ_{00}
+Q_{02} =0, \\
&lQ_{00}-2kQ_{02}-fQ_{22}=0.
\end{align*}
Here we see 
\begin{equation}
\Delta:=\det
\begin{bmatrix}
1 & 0 & 0 \\
-2A & 1 & 0 \\
l & -2k & -f 
\end{bmatrix} =-f=-e^{2F}\not= 0.
\end{equation}

Thus we can claim

\begin{Proposition}\label{PropositionNB}
The set of equations
\eqref{N.9a}, \eqref{N.9b}, \eqref{N.9c} implies 
\eqref{36a}, \eqref{36b}, \eqref{36c}.
\end{Proposition}

Proof. $\Delta \not=0$ guarantees
 $Q_{00}=Q_{02}=Q_{22}=0$.
This means \eqref{36a},\eqref{36b},\eqref{36c}. $\square$. \\

{\bf 4)} It follows from \eqref{30d}\eqref{30e} that
\begin{align*}
2(R_{11}-R_{33})&=-\frac{2}{\Pi}(\partial_1^2\Pi-\partial_3^2\Pi)+
\frac{2}{\Pi}[(\partial_1m)(\partial_1\Pi)-(\partial_3m)(\partial_3\Pi)]+ \\
&+\frac{1}{\Pi^2}
[(\partial_1f)(\partial_1l)+(\partial_1k)^2
-(\partial_3f)(\partial_3l)-(\partial_3k)^2)].
\end{align*}
Therefore \eqref{Td} implies
\begin{align}
-\frac{2}{\Pi}(\partial_1^2\Pi-\partial_3^2\Pi)+&
\frac{2}{\Pi}[(\partial_1m)(\partial_1\Pi)-(\partial_3m)(\partial_3\Pi)]+ \nonumber \\
&+\frac{1}{\Pi^2}
[(\partial_1f)(\partial_1l)+(\partial_1k)^2
-(\partial_3f)(\partial_3l)-(\partial_3k)^2]=0.\label{N.53}
\end{align}
So,
 \eqref{N.53} reads
\begin{align}
(\partial_1\Pi)(\partial_1K)-(\partial_3\Pi)(\partial_3K)&=
\frac{1}{2}(\partial_1^2\Pi-\partial_3^2\Pi)+
\Pi[(\partial_1F)^2-(\partial_3F)^2] + \nonumber \\
&-\frac{e^{4F}}{4\Pi}
[(\partial_1A)^2-(\partial_3A)^2].
\end{align}\\

{\bf 5)} \eqref{36f} and \eqref{30f} imply
that
\begin{align*}
-2\frac{\partial_1\partial_3\Pi}{\Pi}&+
\frac{1}{\Pi}[(\partial_3m)(\partial_1\Pi)+
(\partial_1m)(\partial_3\Pi)] + \\
&+\frac{1}{2\Pi^2}
[(\partial_1f)(\partial_3l)+(\partial_1l)(\partial_3f)+2(\partial_1k)(\partial_3k)]=0.
\end{align*}
But we see that
\begin{align}
&(\partial_1f)(\partial_3l)+(\partial_1l)(\partial_3f)+2(\partial_1k)(\partial_3k)= \nonumber \\
&=-8(\partial_1F)(\partial_3F)\Pi^2+
4\Pi[(\partial_1F)(\partial_3\Pi)+(\partial_3F)(\partial_1\Pi)]+
2e^{4F}(\partial_1A)(\partial_3A ).
\end{align}

Therefore, since $m=-2F+2K$, \eqref{36f} reads
\begin{equation}
(\partial_3\Pi)(\partial_1K)+(\partial_1\Pi)(\partial_3K)=
\partial_1\partial_3\Pi+2\Pi(\partial_1F)(\partial_3F)-
\frac{e^{4F}}{2\Pi}(\partial_1A)(\partial_3A ).
\end{equation}\\

Thus we can claim

\begin{Proposition}
The set of equations \eqref{36d} \eqref{36e} \eqref{36f}, provided \eqref{Td}, implies
\begin{subequations}
\begin{align}
(\partial_1\Pi)(\partial_1K)-(\partial_3\Pi)(\partial_3K)&=
\frac{1}{2}(\partial_1^2\Pi-\partial_3^2\Pi)+
\Pi[(\partial_1F)^2-(\partial_3F)^2] + \nonumber \\
&-\frac{e^{4F}}{4\Pi}
[(\partial_1A)^2-(\partial_3A)^2], \label{Ya} \\
(\partial_3\Pi)(\partial_1K)+(\partial_1\Pi)(\partial_3K)=&
\partial_1\partial_3\Pi+2\Pi(\partial_1F)(\partial_3F)-
\frac{e^{4F}}{2\Pi}(\partial_1A)(\partial_3A ). \label{Yb}
\end{align}
\end{subequations}
\end{Proposition}

However the inverse is doubtful, namely, we are not sure that the system of equations \eqref{Ya} and \eqref{Yb} can recover both \eqref{36d} and  \eqref{36e}, separately, when $\Omega$ is not a constant. \\

Be that as it may, we suppose the following assumption:\\

{\bf (C):} {\it It holds that }
\begin{equation}
\Big(\frac{\partial\Pi}{\partial\varpi}\Big)^2+\Big(\frac{\partial \Pi}{\partial z}\Big)^2
\not=0. 
\end{equation}\\

Under this assumption the set of equations \eqref{Ya}\eqref{Yb} is equivalent to 
\begin{subequations}
\begin{align}
&\partial_1K=((\partial_1\Pi)^2+(\partial_3\Pi)^2)^{-1}\Big[
(\partial_1\Pi)\mathrm{RH}\eqref{Ya}+(\partial_3\Pi)\mathrm{RH}\eqref{Yb}\Big], \label{Yc} \\
&\partial_3K=((\partial_1\Pi)^2+(\partial_3\Pi)^2)^{-1}\Big[
-(\partial_3\Pi)\mathrm{RH}\eqref{Ya}+(\partial_1\Pi)\mathrm{RH}\eqref{Yb}\Big], \label{Yd}
\end{align}
\end{subequations}
where RH\eqref{Ya}, RH\eqref{Yb}  stand for the right-hand sides of 
\eqref{Ya}, \eqref{Yb}, respectively.

We can claim
\begin{Proposition}\label{Prop4}
Suppose that the assumption {\bf (C)} holds on the considered domain $\mathfrak{O}$, and that the domain $\mathfrak{O}$ is the union of a domain $\mathfrak{O}_1$ on which $\Omega$ is a constant and a domain $\mathfrak{O}_2$ on which $\rho=P=0$, namely vacuum.
 Then in $\mathfrak{O}$ the set of equations
\eqref{Ya},\eqref{Yb} implies \eqref{36d},\eqref{36e},\eqref{36f}, provided that the equation
\eqref{Td} and the set of equations \eqref{N.9a},\eqref{N.9b}, \eqref{N.9c} hold. 
\end{Proposition}

Proof.  Since the set of equations \eqref{Ya}\eqref{Yb} is equivalent to the set of equations
\eqref{36d}$-$\eqref{36e}, \eqref{36f}, we have to prove \eqref{36d}$+$\eqref{36e}, namely, we have to prove
\begin{equation}
\frac{1}{2}(R_{11}+R_{33})=\frac{8\pi\mathsf{G}}{\mathsf{c}^4}S_{11}=
\frac{8\pi\mathsf{G}}{\mathsf{c}^4}S_{33}=\frac{4\pi\mathsf{G}}{\mathsf{c}^4}e^m(\epsilon-P).
\label{EqProp4}
\end{equation}

So, we consider 
$$-\frac{1}{2}(R_{11}+R_{33})\cdot (\Pi_1^2+\Pi_3^2)=\clubsuit, $$
where
$$\clubsuit:=(\Pi_1^2+\Pi_3^2)\Big[-\triangle F+\triangle K+
\frac{\triangle\Pi}{2\Pi}-\frac{\Sigma}{4\Pi^2}\Big].$$
Here $\Pi_j, (j=1,3), \triangle\Pi, \triangle K, \triangle F$ stand for
$\partial_j\Pi$,
$$\frac{\partial^2 \Pi}{\partial\varpi^2}+\frac{\partial^2\Pi}{\partial z^2}, \quad
\frac{\partial^2 K}{\partial\varpi^2}+\frac{\partial^2K}{\partial z^2}, \quad
\frac{\partial^2F}{\partial\varpi^2}+\frac{\partial^2F}{\partial z^2}. $$

Differentiating \eqref{Yc} by $\varpi$ and \eqref{Yd} by $z$, we have
\begin{align*}
&(\Pi_1^2+\Pi_3^2)\triangle K= \\
&=-(\Pi_1K_1+\Pi_3K_3)\triangle\Pi+
\frac{1}{2}(\Pi_1\partial_1\triangle\Pi+\Pi_3\partial_3\triangle\Pi)+ \\
&=(\Pi_1^2-\Pi_3^2)(F_1^2-F_3^2)+4\Pi_1\Pi_3F_1F_3+
2\Pi(\Pi_1F_1+\Pi_3F_3)\triangle F + \\
&+\frac{e^{4F}}{4\Pi^2}\Big[(\Pi_1^2-\Pi_3^2+4\Pi(-\Pi_1F_1+\Pi_3F_3))(A_1^2-A_3^2) + \\
&+(2\Pi_1\Pi_3-4\Pi(\Pi_1F_3+\Pi_3F_1))(2A_1A_3) + \\
&-2\Pi(\Pi_1A_1+\Pi_3A_3)\triangle A \Big]
\end{align*}
after tedious calculations. Here $ F_j, A_j, ( j=1,3,)$ stand for
$\partial_jF, \partial_jA$, and
$\triangle A$ means
$ \displaystyle
\frac{\partial^2A}{\partial\varpi^2}+\frac{\partial^2A}{\partial z^2} $.

Denoting by $[S\mathrm{a}], [S\mathrm{b}]$ the right-hand sides of the equations
\eqref{N.9a}, \eqref{N.9b}, we eliminate $\triangle F, \triangle A$. The result is
\begin{align*}
&(\Pi_1^2+\Pi_3^2)\triangle K= \\
&=-(\Pi_1K_1+\Pi_3K_3)\triangle\Pi+
\frac{1}{2}(\Pi_1\partial_1\triangle\Pi+\Pi_3\partial_3\triangle\Pi)+ \\
&-(\Pi_1^2+\Pi_3^2)(F_1^2+F_3^2)-
\frac{e^{4F}}{4\Pi^2}(\Pi_1^2+\Pi_3^2)(A_1^2+A_3^2) + \\
&+2\Pi (\Pi_1F_1+\Pi_3F_3)[S\mathrm{a}]
-\frac{e^{4F}}{2\Pi}(\Pi_1A_1+\Pi_3A_3)[S\mathrm{b}].
\end{align*}

Let us consider the case of constant $\Omega$. Put
\begin{align}
f':=e^{2G}&=e^{2F}\Big(1+\frac{\Omega}{\mathsf{c}}A\Big)^2
-e^{-2F}\frac{\Omega^2}{\mathsf{c}^2}\Pi^2 \nonumber \\
&=f-2\frac{\Omega}{\mathsf{c}}k-\frac{\Omega^2}{\mathsf{c}}l.
\end{align}
Then, using the identity
$$\partial_j\triangle\Pi=
\Big(2(-F_j+K_j)+\frac{\partial_jP}{P}+\frac{\Pi_j}{\Pi}\Big)\triangle\Pi $$
with
$$
\partial_jP=(\epsilon+P)\frac{1}{2f'}\Big(-\partial_jf+
2\frac{\Omega}{\mathsf{c}}\partial_jk+\frac{\Omega^2}{\mathsf{c}^2}\partial_jl\Big),
$$
 provided that $\Omega$ is a constant, and using
\begin{align*}
[S\mathrm{a}]&=\frac{4\pi\mathsf{G}}{\mathsf{c}^4}\frac{e^m}{f}\Big[(\epsilon+P)\frac{1}{f'}((f-\frac{\Omega}{\mathsf{c}}k)^2+\frac{\Omega^2}{\mathsf{c}^2}\Pi^2))+2Pf\Big], \\
[S\mathrm{b}]&=\frac{16\pi\mathsf{G}}{\mathsf{c}^4}\frac{e^m\Pi^2}{f^2f'}
(\epsilon+P)\frac{\Omega}{\mathsf{c}}\Big(\frac{\Omega}{\mathsf{c}}k-f\Big),
\end{align*}
we can deduce
$$\clubsuit =(\Pi_1^2+\Pi_3^2)(-\epsilon+P)\frac{4\pi\mathsf{G}}{\mathsf{c}^4}e^m.$$
Here we have used the identity
$$\Sigma=e^{4F}(A_1^2+A_3^2)-4\Pi^2(F_1^2+F_3^2)+4\Pi(\Pi_1F_1+\Pi_3F_3).$$
This gives the desired equation \eqref{EqProp4}.\\

When the vacuum is considered, we have $\triangle\Pi=[S\mathrm{a}]=[S\mathrm{b}]=0$, and we see $\clubsuit =0$. This completes the proof. $\square$.\\

{\bf 6)} Summing up, we can claim the following

\begin{Theorem}\label{Th.1}

Suppose that the assumption {\bf (C)} holds on the considered domain $\mathfrak{O}$, and that the domain $\mathfrak{O}$ is the union of a domain $\mathfrak{O}_1$ on which $\Omega$ is a constant and a domain $\mathfrak{O}_2$ on which $\rho=P=0$, namely vacuum.
Then on the domain $\mathfrak{O}$  the system of Einstein equations
\eqref{36a} $\sim$ \eqref{36f} is equivalent to the following system:
\begin{subequations}
\begin{align}
&\frac{\partial^2F}{\partial\varpi^2}+\frac{\partial^2F}{\partial z^2}
+\frac{1}{\Pi}\Big(\frac{\partial \Pi}{\partial\varpi}\frac{\partial F}{\partial\varpi}
+\frac{\partial \Pi}{\partial z}\frac{\partial F}{\partial z}\Big)+
\frac{e^{4F}}{2\Pi^2}\Big[\Big(\frac{\partial A}{\partial\varpi}\Big)^2+
\Big(\frac{\partial A}{\partial z}\Big)^2\Big] = \nonumber \\
&=\frac{4\pi\mathsf{G}}{\mathsf{c}^4}
e^{2(-F+K)}
\Big[(\epsilon+P)
\frac{e^{2F}\Big(1+\frac{\Omega}{\mathsf{c}}A\Big)^2+e^{-2F}\frac{\Omega^2}{\mathsf{c}^2}\Pi^2}{e^{2F}\Big(1+\frac{\Omega}{\mathsf{c}}A\Big)^2-e^{-2F}\frac{\Omega^2}{\mathsf{c}^2}\Pi^2}
+2P\Big], \label{N.EQa} \\
&\frac{\partial^2A}{\partial\varpi^2}+\frac{\partial^2A}{\partial z^2}
-\frac{1}{\Pi}\Big(\frac{\partial\Pi}{\partial\varpi}\frac{\partial A}{\partial\varpi}+\frac{\partial\Pi}{\partial z}\frac{\partial A}{\partial z}\Big)+
4\Big(\frac{\partial F}{\partial\varpi}\frac{\partial A}{\partial\varpi}+
\frac{\partial F}{\partial z}\frac{\partial A}{\partial z}\Big) = \nonumber \\
&=-\frac{16\pi\mathsf{G}}{\mathsf{c}^4}
e^{2(-F+K)}(\epsilon+P)
\frac{e^{-2F}\frac{\Omega}{\mathsf{c}}\Pi^2\Big(1+\frac{\Omega}{\mathsf{c}}A\Big)}{e^{2F}\Big(1+\frac{\Omega}{\mathsf{c}}A\Big)^2-e^{-2F}\frac{\Omega^2}{\mathsf{c}^2}\Pi^2}, \label{N.EQb} \\
&\frac{\partial^2\Pi}{\partial \varpi^2}
+\frac{\partial^2\Pi}{\partial z^2}=\frac{16\pi\mathsf{G}}{\mathsf{c}^4}
e^{2(-F+K)}P\Pi, \label{N.EQc} \\
&\frac{\partial\Pi}{\partial \varpi}
\frac{\partial K}{\partial\varpi}-
\frac{\partial\Pi}{\partial z}\frac{\partial K}{\partial z}=
\frac{1}{2}\Big(\frac{\partial^2\Pi}{\partial\varpi^2}
-\frac{\partial^2\Pi}{\partial z^2}\Big)+
\Pi\Big[\Big(\frac{\partial F}{\partial\varpi}\Big)^2
-\Big(\frac{\partial F}{\partial z}\Big)^2\Big] + \nonumber \\
&-\frac{e^{4F}}{4\Pi}
\Big[\Big(\frac{\partial A}{\partial\varpi}\Big)^2
-\Big(\frac{\partial A}{\partial z}\Big)^2\Big], \label{N.EQd} \\
&\frac{\partial\Pi}{\partial z}
\frac{\partial K}{\partial\varpi}
+\frac{\partial\Pi}{\partial \varpi}\frac{\partial K}{\partial z}=
\frac{\partial^2\Pi}{\partial \varpi\partial z}
+2\Pi\frac{\partial F}{\partial\varpi}
\frac{\partial F}{\partial z}-
\frac{e^{4F}}{2\Pi}
\frac{\partial A}{\partial\varpi}\frac{\partial A}{\partial z} , \label{N.EQe} \\
&G=
\Big( =\frac{1}{2}\log\Big[ e^{2F}\Big(1+\frac{\Omega}{\mathsf{c}}A\Big)^2-
e^{-2F}\frac{\Omega^2}{\mathsf{c}^2}\Pi^2\Big]\quad \Big) = \nonumber \\
&=-\frac{u}{\mathsf{c}^2}+\mbox{Const.}, \mbox{while}\quad \rho >0.\label{N.EQf}
\end{align}
\end{subequations}

Here $u, P, \epsilon=\mathsf{c}^2\rho$ are given functions of $\rho$, and $G$ is
determined by $F, A, \Pi, \Omega$ through \eqref{18}.

\end{Theorem}

{\bf 7)}  Now we have a question of the consistency of the first order system of equations for $K$. In order that there exists $K$ which satisfies \eqref{N.EQd}, \eqref{N.EQe}, or there exists $K$ which satisfies 
\begin{align}
\frac{\partial K}{\partial\varpi}&=\Big(
\Big(\frac{\partial\Pi}{\partial\varpi}\Big)^2+\Big(\frac{\partial \Pi}{\partial z}\Big)^2\Big)^{-1}
\Big[\frac{\partial\Pi}{\partial\varpi}\mbox{RH\eqref{N.EQd}}
+\frac{\partial\Pi}{\partial z}\mbox{RH\eqref{N.EQe}}\Big], \label{Z.Ka} \\
\frac{\partial K}{\partial z}&=\Big(
\Big(\frac{\partial\Pi}{\partial\varpi}\Big)^2+\Big(\frac{\partial \Pi}{\partial z}\Big)^2\Big)^{-1}
\Big[-\frac{\partial\Pi}{\partial z}\mbox{RH\eqref{N.EQd}}
+\frac{\partial\Pi}{\partial\varpi}\mbox{RH\eqref{N.EQe}}\Big], \label{Z.Kb}
\end{align}
where RH\eqref{N.EQd}, RH\eqref{N.EQe} stand for the right-hand sides of 
\eqref{N.EQd}, \eqref{N.EQe}, provided the assumption {\bf (C)}, it is necessary that the ` consistency condition' 
\begin{equation}
\frac{\partial \tilde{K}_1}{\partial z}=\frac{\partial \tilde{K}_3}{\partial\varpi} \label{Z.Cons}
\end{equation}
holds, where $\tilde{K}_1, \tilde{K}_3$ stand for the right-hand sides of 
\eqref{Z.Ka}, \eqref{Z.Kb}. 

It is claimed in \cite[Section 4.2, p.56]{Islam} that it is the case when $P=0$ and $\Omega$ is a constant in the considered domain, if we take $\Pi=\varpi$. Actually we can claim
the following: 

\begin{Proposition}\label{PropN5}
Suppose that the assumption {\bf (C)} holds on the considered domain $\mathfrak{O}$, and that the domain $\mathfrak{O}$ is the union of a domain $\mathfrak{O}_1$ on which $\Omega$ is a constant and a domain $\mathfrak{O}_2$ on which $\rho=P=0$, namely vacuum.
Let $K$ be arbitrarily given, and let $F, A, \Pi$ satisfy \eqref{N.EQa}, \eqref{N.EQb}, \eqref{N.EQc}, and \eqref{N.EQf} for $\rho>0, \Omega$ and this given $K$. Denote by 
$\tilde{K}_1, \tilde{K}_3$ the right-hand sides of \eqref{Z.Ka}, \eqref{Z.Kb}, respectively evaluated by these $F, A, \Pi$. Then it holds on $\mathfrak{O}$ that 
\begin{align}
\frac{\partial\tilde{K}_1}{\partial z}-
\frac{\partial\tilde{K}_3}{\partial\varpi}&=
\frac{16\pi \mathsf{G}}{\mathsf{c}^4}
e^{2(-F+K)}P\Pi\Big[\Big(\frac{\partial\Pi}{\partial\varpi}\Big)^2+
\Big(\frac{\partial\Pi}{\partial z}\Big)^2\Big]^{-1}\times \nonumber \\
&\times \Big[\Big(\frac{\partial K}{\partial\varpi}-\tilde{K}_1\Big)
\frac{\partial\Pi}{\partial z}-
\Big(\frac{\partial K}{\partial z}-\tilde{K}_3\Big)
\frac{\partial\Pi}{\partial\varpi}\Big].\label{N3.39}
\end{align}
\end{Proposition}

Proof. By a tedious calculation, we get
\begin{align*}
\frac{\partial\tilde{K}_1}{\partial z}-\frac{\partial\tilde{K}_3}{\partial\varpi}&=
-(\Pi_3\tilde{K}_1-\Pi_1\tilde{K}_3)
(\Pi_1^2+\Pi_3^2)^{-1}\cdot \triangle\Pi + \\
&+(\Pi_1^2+\Pi_3^2)^{-1}Z
\end{align*}
with
\begin{align*}
Z:=&\frac{1}{2}(-\Pi_1\partial_3\triangle\Pi+\Pi_3\partial_1\triangle\Pi)+ \\
&+2\Pi(-\Pi_1F_3+\Pi_3F_1)[S\mathrm{a}]
+\frac{e^{4F}}{2\Pi}(\Pi_1A_3-\Pi_3A_1)[S\mathrm{b}],
\end{align*}
where $\Pi_j, F_j, A_j, (j=1,3,) \triangle\Pi$ stand for $\partial_j\Pi, \partial_jF, \partial_jA, \partial_1^2\Pi+\partial_3^2\Pi$ respectively, and $[S\mathrm{a}], [S\mathrm{b}] $ stand for the right-hand sides of the equations \eqref{N.EQa}, \eqref{N.EQb}. 

Consider $Z$ on the domain $\mathfrak{D}_0$ on which $\Omega$ is a constant. We can write
\begin{align*}
Z=&\frac{1}{2}(-\Pi_1\partial_3\triangle\Pi+\Pi_3\partial_1\triangle\Pi) + \\
&+(-\Pi_1f_3+\Pi_3f_1)\Big(\frac{\Pi}{f}[S\mathrm{a}]-\frac{k}{2\Pi}[S\mathrm{b}]\Big) + \\
&+(-\Pi_1k_3+\Pi_3k_1)\frac{f}{2\Pi}[S\mathrm{b}],
\end{align*}
where $f_j, k_j$ stand for $\partial_jf, \partial_jk$. We have
\begin{align*} 
&\frac{\Pi}{f}[S\mathrm{a}]-\frac{k}{2\Pi}[S\mathrm{b}]=
\frac{8\pi\mathsf{G}}{\mathsf{c}^4}e^m\Big[
(\epsilon+P)\frac{\Pi}{2f^2f'}\Big(f^2+\frac{\Omega^2}{\mathsf{c}^2}(\Pi^2-k^2)\Big)
+\frac{\Pi P}{f}\Big], \\
&[S\mathrm{b}]=\frac{16\pi\mathsf{G}}{\mathsf{c}^4}\frac{e^m\Pi^2}{f^2f'}
(\epsilon+P)\frac{\Omega}{\mathsf{c}}\Big(\frac{\Omega}{\mathsf{c}}k-f\Big),
\end{align*}
where $f'=e^{2G}=f-2\frac{\Omega}{\mathsf{c}}k-\frac{\Omega^2}{\mathsf{c}^2}l$.  Using
the identity
$$\partial_j\triangle\Pi=
\Big(2(-F_j+K_j)+\frac{\partial_jP}{P}+\frac{\Pi_j}{\Pi}\Big)\triangle\Pi
$$
with
$$
\partial_jP=(\epsilon+P)\frac{1}{2f'}\Big(-f_j+
2\frac{\Omega}{\mathsf{c}}k_j+\frac{\Omega^2}{\mathsf{c}^2}\partial_jl\Big) ,
$$
provided that $\Omega$ is a constant, we can deduce that
$$Z=(-\Pi_1\partial_3K+\Pi_3\partial_1K)\triangle\Pi.$$
Here we have used the identity
$$\partial_jl=\frac{1}{f}\Big(2\Pi\Pi_j-\frac{\Pi^2-k^2}{f}f_j-2kk_j\Big),$$
which can be derived from the identity
$\Pi^2=fl+k^2$.
This implies the desired identity. 

On the domain on which $\rho=P=0$, we have $\triangle\Pi=[S\mathrm{a}]=
[S\mathrm{b}]=0$ so that $Z=0$ and
$$\frac{\partial\tilde{K}_1}{\partial z}-\frac{\partial \tilde{K}_3}{\partial \varpi}=0.$$
This completes the proof. $\square$\\

Therefore, as a conclusion of Proposition \ref{PropN5}, if $K$ satisfies \eqref{Z.Ka} \eqref{Z.Kb},
then the consistency condition 
\begin{equation}
\frac{\partial}{\partial z}\mbox{RH\eqref{Z.Ka}}=\frac{\partial}{\partial\varpi}\mbox{RH\eqref{Z.Kb}}
\end{equation}
holds, since 
$$\frac{\partial K}{\partial\varpi}=\tilde{K}_1,\qquad
\frac{\partial K}{\partial z}=\tilde{K}_3.
$$ Of course this conclusion in itself is a vicious circular argument of no use. However the following argument is useful: 

\begin{Theorem} \label{ThN2}
Suppose that the assumption {\bf (C)} holds on the closure $\bar{\mathfrak{D}}$ of the domain $\mathfrak{D}=\{ (\varpi, z)| r=\sqrt{\varpi^2+z^2}<R\}$. Let us suppose that
$\Omega$ is constant on $\{ r \leq R_0\}$, and $\rho=0$ on $\{ r \geq R_1  \}$, where $0<R_1<R_0<R$. 
Suppose that $K \in C^0(\bar{\mathfrak{D}})$ is given and that $F, A, \Pi, u \in C^2(\bar{\mathfrak{D}})$ satisfy \eqref{N.EQa},\eqref{N.EQb},\eqref{N.EQc} and \eqref{N.EQf} with $\rho=f^{\rho}(u), \Omega$ and this $K$. 
 Let us denote by
$\tilde{K}_1,\tilde{K}_3$ the right-hand sides of
\eqref{Z.Ka},\eqref{Z.Kb}, respectively, evaluated by these
$F, A, \Pi$. 
Put
\begin{equation}
\tilde{K}(\varpi, z):=K_O+
\int_0^z\tilde{K}_3(0,z')dz'+
\int_0^{\varpi}\tilde{K}_1(\varpi',z)d\varpi'\label{N3.40}
\end{equation}
for $(\varpi, z) \in \mathfrak{D}$. Here $K_O$ is an arbitrary constant. If $\tilde{K}=K$, then $K$
satisfies 
\begin{equation}
\frac{\partial K}{\partial\varpi}=\tilde{K}_1,\qquad
\frac{\partial K}{\partial z}=\tilde{K}_3,\label{N3.41}
\end{equation}
that is, the equations \eqref{N.EQd}\eqref{N.EQe} are satisfied.
\end{Theorem}

 Proof.  Suppose that $\tilde{K}=K$. It follows from \eqref{N3.40} with $\tilde{K}=K$ that
\begin{align}
\frac{\partial K}{\partial \varpi}(\varpi, z)&=\tilde{K}_1(\varpi, z), \label{N3.42} \\
\frac{\partial K}{\partial z}(\varpi, z)&=
\tilde{K}_3(0,z)+\int_0^{\varpi}
\frac{\partial\tilde{K}_1}{\partial z}(\varpi', z)d\varpi'.\label{N3.43}
\end{align}
Put
\begin{equation}
L(\varpi, z):=\frac{\partial \tilde{K}_1}{\partial z}-
\frac{\partial\tilde{K}_3}{\partial\varpi},
\end{equation}
which is a continuous function on $\bar{\mathfrak{D}}$. Then \eqref{N3.43} reads
\begin{equation}
\frac{\partial K}{\partial z}(\varpi, z)=
\tilde{K}_3(\varpi, z)+
\int_0^{\varpi}
L(\varpi', z)d\varpi'.\label{N3.45}
\end{equation}
Now therefore \eqref{N3.39} of Proposition \ref{PropN5} reads
\begin{equation}
L(\varpi, z)=-
\frac{16\pi\mathsf{G}}{\mathsf{c}^4}
e^{2(-F+K)}P\Pi\Big[\Big(\frac{\partial\Pi}{\partial\varpi}\Big)^2+
\Big(\frac{\partial\Pi}{\partial z}\Big)^2\Big]^{-1}\frac{\partial\Pi}{\partial\varpi}
\int_0^{\varpi}L(\varpi',z)d\varpi'.
\end{equation}
Since the function
$$
\frac{16\pi\mathsf{G}}{\mathsf{c}^4}
e^{2(-F+K)}P\Pi\Big[\Big(\frac{\partial\Pi}{\partial\varpi}\Big)^2+
\Big(\frac{\partial\Pi}{\partial z}\Big)^2\Big]^{-1}\frac{\partial\Pi}{\partial\varpi}
$$
is bounded on the compact $\bar{\mathfrak{D}}$, the Gronwall's argument implies that
$L(\varpi, z)=0$ on $\mathfrak{D}$ so that \eqref{N3.45} reads
\begin{equation}
\frac{\partial K}{\partial z}(\varpi, z)=\tilde{K}_3(\varpi, z).\label{N3.47}
\end{equation}
Thus \eqref{N3.42} and \eqref{N3.47} complete the proof. $\square$

\section{Existence of asymptotically flat metrics }

In this section we suppose that $\Omega$ is a constant $\Omega_{\mathsf{O}}$ in a bounded domain.

\subsection{Post-Newtonian approximation}

We use the following notations:\\

\begin{Notation}

1) Let $n=3,4,5, \Xi >0$. We denote

\begin{align}
&{B}^{(n)}(\Xi)=\{ \mbox{\boldmath$\xi$}=(\xi_1,\cdots,\xi_n)\in \mathbb{R}^n \  |\  |\mbox{\boldmath$\xi$}|:=\sqrt{\sum_k(\xi_k)^2} < \Xi \}, \nonumber \\
&
\bar{B}^{(n)}(\Xi)=\{ \mbox{\boldmath$\xi$}=(\xi_1,\cdots,\xi_n)\in \mathbb{R}^n \  |\  |\mbox{\boldmath$\xi$}| \leq \Xi \}, \nonumber \\
&\partial {B}^{(n)}(\Xi)=\{ \mbox{\boldmath$\xi$}=(\xi_1,\cdots,\xi_n)\in \mathbb{R}^n \  |\  |\mbox{\boldmath$\xi$}|= \Xi \}.
\end{align}

For a continuous function $f$ on $\bar{B}^{(n)}(\Xi)$ and $l=0,1,2$, we put
$$
\|f ; C^{l}(\bar{B}^{(n)}(\Xi))\|:=
\sup_{|L|\leq l,|\mbox{\boldmath$\xi$}|\leq\Xi}|\partial_{\mbox{\boldmath$\xi$}}^Lf(\mbox{\boldmath$\xi$})|,
$$
where
$$ \partial_{\mbox{\boldmath$\xi$}}^L =\Big(\frac{\partial}{\partial\xi_1}\Big)^{L_1}\cdots
\Big(\frac{\partial}{\partial \xi_n}\Big)^{L_n} $$
for $L=(L_1,\cdots, L_n)$ and $|L|=L_1+\cdots +L_n$. 
This is the norm of the Banach space
$ C^l(\bar{B}^{(n)}(\Xi)) $.

Let us fix a number $\alpha$ such that
\begin{equation}
0<\alpha <\Big(\frac{1}{\gamma-1}-1\Big)\wedge  1. \label{alpha}
\end{equation}

For a continuous function $f$ on $\bar{B}^{(n)}(\Xi)$ and $l=0,1,2$, we put
\begin{align*}
\|f ; C^{l,\alpha}(\bar{B}^{(n)}(\Xi))\|:=&
\sup_{|L|\leq l,|\mbox{\boldmath$\xi$}|\leq\Xi}|\partial_{\mbox{\boldmath$\xi$}}^Lf(\mbox{\boldmath$\xi$})| + \nonumber \\
&+\sup_{|\mbox{\boldmath$\xi$}'|,|\mbox{\boldmath$\xi$}|\leq \Xi, 
0<|\mbox{\boldmath$\xi$}'-\mbox{\boldmath$\xi$}|\leq 1,
|L|=l}
\frac{|\partial_{\mbox{\boldmath$\xi$}}^Lf(\mbox{\boldmath$\xi$}')-\partial_{\mbox{\boldmath$\xi$}}^Lf(\mbox{\boldmath$\xi$})|}{|\mbox{\boldmath$\xi$}'-\mbox{\boldmath$\xi$}|^{\alpha}}.
\end{align*}
 This is the norm of the Banach space
$$C^{l,\alpha}(\bar{B}^{(n)}(\Xi))=\{ f \in C^l(\bar{B}^{(n)}(\Xi))\  |\  
\|f; C^{l,\alpha}(\bar{B}^{(n)}(\Xi))\|<\infty \}.
$$

2) Let $R >0$. We denote 
\begin{align}
&\mathfrak{D}(R)=\{ (\varpi,z) \ |\    r:=\sqrt{\varpi^2+z^2}<R\}, \nonumber \\
&\bar{\mathfrak{D}}(R) =\{ (\varpi, z) \  |\   r \leq R\}, \nonumber  \\
&\partial\mathfrak{D}(R)=\{ (\varpi, z) \  |\   r=R\}.
\end{align}\\

 For a function $Q$ of $(\varpi, z)$ and $n=3,4,5$, we denote by $Q^{\flat(n)}$ the
function of $\mbox{\boldmath$\xi$}=(\xi_1,\cdots, \xi_n)$ defined by
\begin{equation}
Q^{\flat(n)}(\mbox{\boldmath$\xi$})=Q(\varpi, z) 
\end{equation}
with
$$\varpi=\mathsf{a}\sqrt{(\xi_1)^2+\cdots+(\xi_{n-1})^2},\qquad z=\mathsf{a}\xi_n,$$
where $\mathsf{a}$ is a positive parameter specified later.\\

 We put
\begin{align*}
\mathfrak{C}^{l}(\bar{\mathfrak{D}}(R))&:=\{ Q \in C(\bar{\mathfrak{D}}(R)) \  |  \\
&   
Q(\varpi, -z)=Q(\varpi, z)\quad\forall \varpi, z \quad  
\mbox{and}\quad Q^{\flat(n)}\in C^{l}(\bar{B}^{(n)}(R/\mathsf{a})) \}, \\
\|Q; \mathfrak{C}^{l}(\bar{\mathfrak{D}}(R))\|&:=\|Q^{\flat(n)}; C^{l}(\bar{B}^{(n)}(R/\mathsf{a}))\|;
\end{align*}
\begin{align*}
\mathfrak{C}^{l,\alpha}(\bar{\mathfrak{D}}(R))&:=\{ Q \in C(\bar{\mathfrak{D}}(R)) \  |  \\
&   
Q(\varpi, -z)=Q(\varpi, z)\quad\forall \varpi, z \quad  
\mbox{and}\quad Q^{\flat(n)}\in C^{l,\alpha}(\bar{B}^{(n)}(R/\mathsf{a})) \}, \\
\|Q; \mathfrak{C}^{l,\alpha}(\bar{\mathfrak{D}}(R))\|&:=\|Q^{\flat(n)}; C^{l,\alpha}(\bar{B}^{(n)}(R/\mathsf{a}))\|.
\end{align*}
\end{Notation}

It is easy to see that the spaces 
$\mathfrak{C}^{l}(\mathfrak{D}(R))$, 
$\mathfrak{C}^{l,\alpha}(\mathfrak{D}(R))$ 
and the norms
 do not depend on the choice of $n$, but depend only on $\mathsf{a}$.\\

Now let us recall the result on the stationary rotating solution of the Euler-Poisson equations, which is the Newtonian limit, obtained by \cite{JJTM1} and \cite{JJTM2}. Namely, it is a solution $(\rho, \vec{v})$ of the form 
$\displaystyle \rho=\rho(\varpi, z), \vec{v}=\Omega_{\mathsf{O}}\frac{\partial}{\partial\phi}$ to the  system
\begin{align*}
& -\rho(1-\zeta^2)r\Omega_{\mathsf{O}}^2+\frac{\partial P}{\partial r}
+\rho\frac{\partial\Phi}{\partial r}=0, \\
&\rho\zeta r^2\Omega_{\mathsf{O}}^2+\frac{\partial P}{\partial\zeta}+
\rho\frac{\partial\Phi}{\partial\zeta}=0, \\
&\triangle\Phi=4\pi\mathsf{G}\rho.
\end{align*}
Here $r=\sqrt{\varpi^2+z^2}, \zeta=z/r$. The density distribution $\rho$ is given as
\begin{equation}
\rho=\rho_{\mathsf{N}}=\rho_{\mathsf{N}\mathsf{O}}\Big(\Theta\Big(\frac{r}{\mathsf{a}},\zeta;
 \frac{1}{\gamma-1}, \mathsf{b}\Big)\vee 0\Big)^{\frac{1}{\gamma-1}}.
\end{equation}

Here   $\Theta(|\mbox{\boldmath$\xi$}|,\zeta;\frac{1}{\gamma-1},\mathsf{b}),
\mbox{\boldmath$\xi$}=(\xi_1,\xi_2,\xi_3), \zeta=\xi_3/|\mbox{\boldmath$\xi$}|$ is the distorted Lane-Emden function with the following properties:

1) The function $\mbox{\boldmath$\xi$} \mapsto \Theta(|\mbox{\boldmath$\xi$}|,\zeta)$ belongs to $ C^{2,\alpha}(\bar{B}^{(3)}(\Xi_0))$, where 
$\Xi_0 := 4\xi_1(\frac{1}{\gamma-1})$, $\xi_1(\nu)$ being the zero of the Lane-Emden function $\theta(\xi; \nu)$ of index $\displaystyle \nu=\frac{1}{\gamma-1}$:\, that is, the solution of
$$-\frac{1}{\xi^2}\frac{d}{d\xi}\xi^2\frac{d\theta}{d\xi}=(\theta\vee 0)^{\nu},\quad \theta|_{\xi=0}=1.$$

2) $\Theta\Big(0,\zeta;\frac{1}{\gamma-1},\mathsf{b}\Big) =1$ and there is a curve
$\zeta\in [-1,1]\mapsto |\mbox{\boldmath$\xi$}|=\Xi_1(\zeta;\frac{1}{\gamma-1},\mathsf{b})$ such that $\Xi_1(\zeta;\frac{1}{\gamma-1},\mathsf{b}) < 2\xi_1(\frac{1}{\gamma-1})$ and, for $ 0\leq|\mbox{\boldmath$\xi$}|\leq \Xi_0$, it holds that
$$ 0<\Theta(|\mbox{\boldmath$\xi$}|,\zeta;\frac{1}{\gamma-1},\mathsf{b}) \quad\Leftrightarrow \quad 0\leq|\mbox{\boldmath$\xi$}|<\Xi_1(\zeta;\frac{1}{\gamma-1},\mathsf{b}).
$$
Moreover it holds that
$$\frac{\partial\Theta}{\partial|\mbox{\boldmath$\xi$}|} \leq -\frac{1}{C}|\mbox{\boldmath$\xi$}| $$
everywhere with a positive number $C$. Therefore $\Theta \leq 1$ everywhere. 

The positive number $\rho_{\mathsf{N}\mathsf{O}}$ is the central density, a positive number, and the parameters $\mathsf{a}, \mathsf{b}$ are specified by

\begin{equation}
\mathsf{a}=\sqrt{\frac{\mathsf{A}\gamma}{4\pi\mathsf{G}(\gamma-1)}}
\rho_{\mathsf{N}\mathsf{O}}^{-\frac{2-\gamma}{2}},
\quad
\mathsf{b}=\frac{\Omega_{\mathsf{O}}^2}{4\pi\mathsf{G}\rho_{\mathsf{N}\mathsf{O}}}. \label{Def.ab}
\end{equation}

Here we are supposing the following assumptions: \\

{\bf (D0)}: {\it It holds that }
\begin{equation}
\frac{6}{5} <\gamma <2,
\end{equation}\\

{\bf (D1)} : $\mathsf{b}$ {\it is sufficiently small, say,
$\mathsf{b}\leq \beta^0, \beta^0$ being a positive number depending on $ \gamma$.} \\

\begin{Remark}
We can take $\Xi_0$ arbitrarily large, but we should take $\beta^0$ small for large $\Xi_0$. In fact $\Theta$ is the solution of 
\begin{align*}
&\Theta=\frac{\mathsf{b}}{4}((\xi_1)^2+(\xi_2)^2)+\mathcal{G}(\Theta), \\
&\mathcal{G}(\Theta)=\mathcal{K}^{(3)}(\Theta\vee 0)^{\frac{1}{\gamma-1}}-\mathcal{K}^{(3)}(\Theta\vee 0)^{\frac{1}{\gamma-1}}(O)+1, \\
&\mathcal{K}^{(3)}g(\mbox{\boldmath$\xi$})=\frac{1}{4\pi}\int\frac{g(\mbox{\boldmath$\xi$}')}{|\mbox{\boldmath$\xi$}-\mbox{\boldmath$\xi$}'|}d\mbox{\boldmath$\xi$}',
\end{align*}
and
$\displaystyle \frac{\mathsf{b}}{4}((\xi_1)^2+(\xi_2)^2)$ grows infinitely large as $|\mbox{\boldmath$\xi$}|\rightarrow \infty$, while $\mathcal{G}(\Theta)$ is bounded.
\end{Remark}

We put 
\begin{equation}
u_{\mathsf{N}}=u_{\mathsf{O}}\Theta\Big(\frac{r}{\mathsf{a}},\zeta;\frac{1}{\gamma-1},\mathsf{b}\Big),\quad
u_{\mathsf{O}}=\frac{\mathsf{A}\gamma}{\gamma-1}\rho_{\mathsf{N}\mathsf{O}}^{\gamma-1}
\end{equation}
and we denote by $\Phi_{\mathsf{N}}$ the (Newtonian) gravitational potential generated by $\rho_{\mathsf{N}}$, namely,
$$\Phi_{\mathsf{N}}^{\flat(3)}(\mbox{\boldmath$\xi$})=-\mathsf{G}\int
\frac{\rho_{\mathsf{N}}^{\flat(3)}(\mbox{\boldmath$\xi$}')}{|\mbox{\boldmath$\xi$}-\mbox{\boldmath$\xi$}'|}d\mbox{\boldmath$\xi$}'. $$
Note that the integral is performed on $|\mbox{\boldmath$\xi$}'|\leq \Xi_0$ and $\Phi_{\mathsf{N}}$ is defined for $\forall r \in [0,+\infty[$.\\

So we are using the following
\begin{Definition}\label{Def.R0}
 We shall denote $r_1=\mathsf{a}\xi_1(\frac{1}{\gamma-1})$, and put $$R_0=4r_1=\mathsf{a}\Xi_0=
4\mathsf{a}\xi_1(\frac{1}{\gamma-1}),$$
where
 $\xi_1(\nu)$ stands for the zero of the Lane-Emden function $\theta(\xi; \nu)$ of index $\displaystyle \nu=\frac{1}{\gamma-1}$. 
We put
\begin{equation}
\mathsf{a}=\frac{1}{\sqrt{4\pi\mathsf{G}}}\Big(\frac{\mathsf{A}\gamma}{\gamma-1}\Big)^{\frac{1}{2(\gamma-1)}}
u_{\mathsf{O}}^{-\frac{2-\gamma}{2(\gamma-1)}}, \quad
\mathsf{b}=\frac{1}{4\pi\mathsf{G}}\Big(\frac{\mathsf{A}\gamma}{\gamma-1}\Big)^{\frac{1}{\gamma-1}}
\Omega_{\mathsf{O}}^2u_{\mathsf{O}}^{-\frac{1}{\gamma-1}}.
\end{equation}
\end{Definition}

When we consider the potentials and the associated metric on the whole space, we assume that $\Omega$ vanishes for large $r$, namely, we assume
$
\Omega=\Omega_{\mathsf{O}}\chi(r/R_0)$, where  $\chi \in C^{\infty}(\mathbb{R})$ is a cut off function such that
$\chi(t)=1$ for $t\leq 1$, $0 \leq \chi(t)\leq 1$ for $1<t<2$ and $\chi(t)=0$ for $t \geq 2$,
without loss of generality. Actually, since we are going to establish a density distribution $\rho$ such that $\{ \rho >0\} \subset \mathfrak{D}(3r_1)$,   the values of $\Omega$ in 
$\bar{\mathfrak{D}}(R_0)^{\mathsf{c}}=\bar{\mathfrak{D}}(4r_1)^{\mathsf{c}}$ do not affect the angular momentum
$T^{\mu\nu}=(\epsilon +P)U^{\mu}U^{\nu}-P g^{\mu\nu} $
generated by the 4-velocity
$$U^{\mu}\frac{\partial}{\partial x^{\mu}}=e^{-G}\Big(
\frac{1}{\mathsf{c}}\frac{\partial}{\partial t}+\frac{\Omega}{\mathsf{c}}\frac{\partial}{\partial\phi}\Big), $$
therefore we can assign arbitrary values for $\Omega$ on $\bar{\mathfrak{D}}(R_0)^{\mathsf{c}}$. 
Otherwise, if we do not perform cut off, e.g., if $\Omega=O(1)$ but $\Omega \rightarrow \Omega_{\infty}\not=0$ as $r \rightarrow +\infty$, and if the metric is asymptotically flat, then we have
$$e^{2F}\Big(1+\frac{\Omega}{\mathsf{c}}A\Big)^2-
e^{-2F}\frac{\Omega^2}{\mathsf{c}^2}\Pi^2=-\frac{\Omega_{\infty}^2}{\mathsf{c}^2}\varpi^2(1+o(1))+O(1)$$
as $r\rightarrow +\infty$. This turns out to be negative for large $\varpi$ and the condition {\bf (B)} would break down. 

In view of this extension of $\Omega$ on the whole space, we consider $u_{\mathsf{N}}$ is defined for $\forall r \in [0,+\infty[$ as the solution of the equation
\begin{equation}
u_{\mathsf{N}}=\frac{\Omega(r)^2}{2}\varpi^2+
\Phi_{\mathsf{N}}-\Phi_{\mathsf{N}}(O)+u_{\mathsf{O}}.
\end{equation}
on the whole space. 
Here recall that $\Phi_{\mathsf{N}}$ is already defined on the whole space as the Newtonian potential of the density distribution $\rho_{\mathsf{N}}$ whose support is included in $\mathfrak{D}(2r_1)$.
So, we can assume that $\sup_{r \geq R_0}u_{\mathsf{N}}<0$ 
by replacing  $\beta^0$ in the condition {\bf (D1)} by a smaller one if necessary,
since $u_{\mathsf{N}}$ tends to $u_{\mathsf{O}}\theta$ uniformly on $\mathfrak{D}(R_0)^{\mathsf{c}}$ as $\mathsf{b} \rightarrow 0$, where $\theta$ is the Lane-Emden function
of index $\displaystyle \nu=\frac{1}{\gamma-1}$, which satisfies $$\theta= -\mu_1
(\nu)\Big(\frac{1}{\xi}-
\frac{1}{\xi_1(\nu)}\Big)
 \quad \mbox{for}\quad \xi \geq \xi_1(\nu),$$
where $\mu_1 (\nu)>0$.
We have
$$ u_{\mathsf{N}}=-\Phi_{\mathsf{N}}(O)+u_{\mathsf{O}}+O\Big(\frac{1}{r}\Big)$$
as $r \rightarrow +\infty$.
In other words, we have extended the distorted Lane-Emden function $\Theta$ onto the whole space so that $\sup\{\Theta \  |\   |\mbox{\boldmath$\xi$}|\geq \Xi_0\} < 0$ holds and the  equation
$$\Theta=\frac{\mathsf{b}}{4}\chi\Big(\frac{|\mbox{\boldmath$\xi$}|}{\Xi_0}\Big)^2((\xi_1)^2+(\xi_2)^2)+\mathcal{G}(\Theta)$$
holds everywhere. 
We have 
$$\Theta=-\mathcal{K}^{(3)}(\Theta\vee 0)^{\frac{1}{\gamma-1}}(O)+1+O\Big(\frac{1}{|\mbox{\boldmath$\xi$}|}\Big) $$
as $|\mbox{\boldmath$\xi$}|\rightarrow +\infty$.\\

We are going to find a set of potentials of the form
\begin{subequations}
\begin{align}
&F=\frac{1}{\mathsf{c}^2}\Phi_{\mathsf{N}}-\frac{1}{\mathsf{c}^4}W, \quad W(O)=0,  \label{PNDef.a}\\
&A=\frac{1}{\mathsf{c}^3}\varpi^2Y, \label{PNDef.b}\\
&\Pi=\varpi\Big(1+\frac{1}{\mathsf{c}^4}X\Big), \label{PNDef.c}\\
&K=\frac{1}{\mathsf{c}^4}V, \label{PNDef.d}\\
&u=u_{\mathsf{N}}+\frac{1}{\mathsf{c}^2}w,\quad w(O)=0. \label{PNDef.e}
\end{align}
\end{subequations}
From \eqref{PNDef.e} $G$ is specified as
\begin{equation}
G=\frac{1}{\mathsf{c}^2}\Phi_{\mathsf{N}}-\frac{\Omega^2}{2\mathsf{c}^2}\varpi^2-\frac{w}{\mathsf{c}^4}. \label{Gw}
\end{equation}
Therefore the equation \eqref{18} leads us to the relation between $w$ and $W$ as
$$-\frac{\Omega^2}{2\mathsf{c}^2}\varpi^2+\frac{1}{\mathsf{c}^4}(W-w)=
\frac{1}{2}\log\Big(\Big(1+\frac{\Omega}{\mathsf{c}}A\Big)^2-e^{-4F}\frac{\Omega^2}{\mathsf{c}^2}\Pi^2\Big).$$
This relation can be written as
\begin{align}
w=&W-
\frac{\Omega^2\varpi^2}{2}\mathsf{c}^2\Big(1-e^{-\frac{4}{\mathsf{c}^2}(\Phi_{\mathsf{N}}-\frac{W}{\mathsf{c}^2})}
\Big(1+\frac{X}{\mathsf{c}^4}\Big)^2\Big) + \nonumber \\
&-\Omega\varpi^2Y-\frac{\Omega^2\varpi^4}{2\mathsf{c}^4}Y^2
-\frac{\mathsf{c}^4}{2}\sum_{k=2}^{\infty}\frac{(-1)^{k+1}}{k}\Big(\frac{Z}{\mathsf{c}^2}\Big)^k \label{0313}
\end{align}
for
\begin{equation}
Z:=\frac{2\Omega\varpi^2}{\mathsf{c}^2}Y+
\frac{\Omega^2\varpi^4}{\mathsf{c}^6}Y^2
-\Omega^2\varpi^2e^{-\frac{4}{\mathsf{c}^2}(\Phi_{\mathsf{N}}-\frac{W}{\mathsf{c}^2})}\Big(1+
\frac{X}{\mathsf{c}^4}\Big)^2, \label{0314}
\end{equation}
provided that $|Z|/\mathsf{c}^2 <1$. Actually we have put
$$
\Big(1+\frac{\Omega}{\mathsf{c}}A\Big)^2-e^{-4F}\frac{\Omega^2}{\mathsf{c}^2}\Pi^2
=1+\frac{Z}{\mathsf{c}^2}.
$$
We introduce the auxiliary quantity $\mathfrak{Q}_0$ defined by
\begin{equation}
w=W-2\Omega^2\varpi^2\Phi_{\mathsf{N}}-\Omega\varpi^2Y+
\frac{1}{4}\Omega^4\varpi^4+\frac{1}{\mathsf{c}^2}\mathfrak{Q}_0.\label{Ww}
\end{equation}
Note that
\begin{align}
\frac{1}{\mathsf{c}^2}\mathfrak{Q}_0&=\frac{\Omega^2\varpi^2}{2}\Big(e^{-\frac{4}{\mathsf{c}^2}(\Phi_{\mathsf{N}}-\frac{W}{\mathsf{c}^2})}\Big(1+\frac{X}{\mathsf{c}^4}\Big)+\frac{4}{\mathsf{c}^2}\Phi_{\mathsf{N}}\Big) + \nonumber \\
&+\frac{Z^2}{4}-\frac{1}{4}\Omega^2\varpi^2
-\frac{\mathsf{c}^4}{2}\sum_{k=3}^{\infty}
\frac{(-1)^{k+1}}{k}\Big(\frac{Z}{\mathsf{c}^2}\Big)^k, \label{WwZ}
\end{align}
while
$$Z=-\Omega^2\varpi^2+\frac{1}{\mathsf{c}^2}Z_1 $$
with
$$Z_1=2\Omega\varpi^2Y+
\frac{\Omega^2\varpi^4}{\mathsf{c}^4}Y+
\mathsf{c}^2\Omega^2\varpi^2\Big(1-
e^{-\frac{4}{\mathsf{c}^2}(\Phi_{\mathsf{N}}-\frac{W}{\mathsf{c}^2})}\Big(1+\frac{X}{\mathsf{c}^4}\Big)^2\Big).
$$\\

By tedious calculations the equations \eqref{N.EQa}, \eqref{N.EQb}, \eqref{N.EQc}
read as 
\begin{subequations}
\begin{align}
&\Big[\frac{\partial^2}{\partial\varpi^2}+\frac{1}{\varpi}\frac{\partial}{\partial\varpi}
+\frac{\partial^2}{\partial z^2}
+4\pi\mathsf{G}\frac{1}{\gamma-1}\frac{\rho_{\mathsf{N}}}{u_{\mathsf{N}}}
\Big]W+ \nonumber \\
&-4\pi\mathsf{G}\frac{1}{\gamma-1}\frac{\rho_{\mathsf{N}}}{u_{\mathsf{N}}}
\Omega\varpi^2\Big(Y+2\Omega\Phi_{\mathsf{N}}-\frac{1}{4}\Omega^3\varpi^2\Big) + \nonumber \\
&+
4\pi\mathsf{G}\Upsilon_1\rho_{\mathsf{N}}u_{\mathsf{N}} -8\pi\mathsf{G}\rho_{\mathsf{N}}(\Phi_{\mathsf{N}}+2\Omega^2\varpi^2)+
12\pi\mathsf{G}P_{\mathsf{N}}+\mathfrak{R}_a=0, \label{PNEq.a} \\
&\Big[\frac{\partial^2}{\partial\varpi^2}+\frac{3}{\varpi}\frac{\partial}{\partial\varpi}
+\frac{\partial^2}{\partial z^2}\Big]Y+16\pi\mathsf{G}\Omega\rho_{\mathsf{N}}+\mathfrak{R}_b=0, \label{PNEq.b} \\
&\Big[\frac{\partial^2}{\partial\varpi^2}+\frac{2}{\varpi}\frac{\partial}{\partial\varpi}
+\frac{\partial^2}{\partial z^2}\Big]X
-16\pi\mathsf{G}P_{\mathsf{N}}+\mathfrak{R}_c=0. \label{PNEq.c}
\end{align}
\end{subequations}

Here $\Upsilon_1$ is the coefficient of the expansion 
$\Upsilon_{\rho}(\eta)=\sum_{k=1}^{\infty}\Upsilon_k\eta^k$.
 We define $\mathfrak{R}_a, \mathfrak{R}_b, \mathfrak{R}_c$ by using auxiliary quantities $\mathfrak{Q}_1 \sim \mathfrak{Q}_6 $ as follows:

\begin{equation}
\mathfrak{Q}_1:=e^{-4F}{\Omega^2}\varpi^2
\Big(1+\frac{X}{\mathsf{c}^4}\Big)^2
\Big(1+\frac{\Omega}{\mathsf{c}^4}\varpi^2Y\Big)^{-2}, 
\end{equation}
\begin{subequations}
\begin{align}
[Q1]:=&\rho-\rho_{\mathsf{N}}=\frac{1}{\mathsf{c}^2}(Df_{\mathsf{N}}^{\rho}(u_{\mathsf{N}})w+\Upsilon_1\rho_{\mathsf{N}}u_{\mathsf{N}})+H_{\rho}(w)+\frac{1}{\mathsf{c}^4}\mathfrak{Q}_2, \label{Def.Q2} \\
[Q2]:=&-\mathsf{c}^2\rho 
\Big(e^{2(-F+K)}-
\frac{ 1+\frac{1}{\mathsf{c}^2}\mathfrak{Q}_1 }{ 1-\frac{1}{\mathsf{c}^2}\mathfrak{Q}_1}\Big)= 2\rho_{\mathsf{N}}(\Phi_{\mathsf{N}}+2\mathfrak{Q}_1)+\frac{1}{\mathsf{c}^2}\mathfrak{Q}_3 = \nonumber \\
&=2\rho_{\mathsf{N}}(\Phi_{\mathsf{N}}+2\Omega^2\varpi^2)+\frac{1}{\mathsf{c}^2}\mathfrak{Q}_4, \\
-e^{2(-F+K)}&\Big(
\mathsf{c}^2\rho\frac{1+\frac{1}{\mathsf{c}^2}\mathfrak{Q}_1}{1-\frac{1}{\mathsf{c}^2}\mathfrak{Q}_1}+
P\frac{3-\frac{2}{\mathsf{c}^2}\mathfrak{Q}_1}{1-\frac{1}{\mathsf{c}^2}\mathfrak{Q}_1}\Big)
+\mathsf{c}^2\rho_{\mathsf{N}}= \nonumber \\
&=-\mathsf{c}^2[Q1] +[Q2] - 3P
-P\Big(
e^{2(-F+K)}-\frac{\frac{1}{\mathsf{c}^2}\mathfrak{Q}_1}{1-\frac{1}{\mathsf{c}^2}\mathfrak{Q}_1}\Big) = \nonumber \\
&=-(Df_{\mathsf{N}}^{\rho}(u_{\mathsf{N}})w+\Upsilon_1\rho_{\mathsf{N}}u_{\mathsf{N}})+
2\rho_{\mathsf{N}}(\Phi_N+2\Omega^2\varpi^2)+ \nonumber \\
&-3P_{\mathsf{N}}-
\mathsf{c}^2H_{\rho}(w)+\frac{1}{\mathsf{c}^2}\mathfrak{Q}_5;\\
\mathfrak{R}_a:=&-\frac{1}{\mathsf{c}^2}\Big(1+\frac{X}{\mathsf{c}^4}\Big)^{-1}
[(\partial_1X)\partial_1+(\partial_3X)\partial_3](\Phi_{\mathsf{N}}-\frac{W}{\mathsf{c}^2})+ \nonumber \\
&+\frac{1}{\mathsf{c}^2}\frac{e^{2F}}{2}\Big(1+\frac{X}{\mathsf{c}^4}\Big)^{-2}
((2Y+\varpi\partial_1Y)^2+
(\varpi\partial_3Y)^2) + \nonumber \\
&
+4\pi\mathsf{G}\frac{1}{\gamma-1}\frac{\rho_{\mathsf{N}}}{u_{\mathsf{N}}}\frac{1}{\mathsf{c}^2}\mathfrak{Q}_0
+4\pi\mathsf{G}\mathsf{c}^2H_{\rho}(w)-4\pi\mathsf{G}\frac{1}{\mathsf{c}^2}\mathfrak{Q}_5; 
\end{align}
\end{subequations}
\begin{subequations}
\begin{align}
\frac{1}{\mathsf{c}^2}e^{-6F+2K}&(\mathsf{c}^2\rho+P)
\Big(1-\frac{1}{\mathsf{c}^2}\mathfrak{Q}_1\Big)^{-1}
\Big(1+\frac{X}{\mathsf{c}^4}\Big)^2
\Big(1+\frac{\Omega}{\mathsf{c}^4}\varpi^2Y\Big)^{-1}= \nonumber \\
&=\rho_{\mathsf{N}}+\frac{1}{\mathsf{c}^2}\mathfrak{Q}_6;  \\
\mathfrak{R}_b:=&-\frac{1}{\mathsf{c}^4}\Big(1+\frac{X}{\mathsf{c}^4}\Big)^{-1}
[(\partial_1X)\partial_1+(\partial_3X)\partial_3-\frac{2}{\varpi}(\partial_1X)\cdot]Y +\nonumber \\
&+\frac{4}{\mathsf{c}^4}[\partial_1(\Phi_{\mathsf{N}}-\frac{W}{\mathsf{c}^2})\partial_1+
\partial_3(\Phi_{\mathsf{N}}-\frac{W}{\mathsf{c}^2})\partial_3+
\frac{2}{\varpi}\partial_1(\Phi_{\mathsf{N}}-\frac{W}{\mathsf{c}^2})\cdot]Y + \nonumber \\
&+16\pi\mathsf{G}\frac{\Omega}{\mathsf{c}^2}\mathfrak{Q}_6; 
\end{align}
\end{subequations}
\begin{equation}
\mathfrak{R}_c:=-16\pi\mathsf{G}\Big(
e^{2(-F+K)}P\Big(1+\frac{X}{\mathsf{c}^4}\Big)-P_{\mathsf{N}}\Big).
\end{equation}

Here we have used the following

\begin{Definition}
We put
\begin{align}
H_{\rho}(w)&:=f_{\mathsf{N}}^{\rho}(u_{\mathsf{N}}+\frac{w}{\mathsf{c}^2})
-f_{\mathsf{N}}^{\rho}(u_{\mathsf{N}})
-Df_{\mathsf{N}}^{\rho}(u_{\mathsf{N}})\frac{w}{\mathsf{c}^2} \nonumber \\
&=f_{\mathsf{N}}^{\rho}(u_{\mathsf{N}}+\frac{w}{\mathsf{c}^2})
-f_{\mathsf{N}}^{\rho}(u_{\mathsf{N}})
-
\frac{1}{\gamma-1}\frac{\rho_{\mathsf{N}}}{u_{\mathsf{N}}}
\frac{w}{\mathsf{c}^2}. \label{Def.H}
\end{align}
\end{Definition}
As for the definition of $f_{\mathsf{N}}^{\rho}$, we recall Definition \ref{Def.frho}.
Note that, therefore, the quantity $\mathfrak{Q}_2$ defined by \eqref{Def.Q2} can be explicitly 
written as
\begin{align}
\frac{1}{\mathsf{c}^4}\mathfrak{Q}_2&=
\frac{1}{\mathsf{c}^4}\rho_{\mathsf{N}}w+
\rho_{\mathsf{N}}\sum_{k=2}^{\infty}\Upsilon_k\Big(\frac{u_{\mathsf{N}}}{\mathsf{c}^2}+\frac{w}{\mathsf{c}^4}\Big)^k + \nonumber \\
&+\Big(f_{\mathsf{N}}^{\rho}\Big(u_{\mathsf{N}}+\frac{w}{\mathsf{c}^2}\Big)-f_{\mathsf{N}}^{\rho}(u_{\mathsf{N}})\Big)
\sum_{k=1}^{\infty}
\Upsilon_k\Big(\frac{u_{\mathsf{N}}}{\mathsf{c}^2}+\frac{w}{\mathsf{c}^4}\Big)^k,
\end{align}
provided that $\displaystyle \frac{u}{\mathsf{c}^2}=\frac{u_{\mathsf{N}}}{\mathsf{c}^2}+\frac{w}{\mathsf{c}^4} \ll 1$.
Of course, we read 
\begin{align*}
w=&W-
\frac{\Omega^2\varpi^2}{2}\mathsf{c}^2\Big(1-e^{-\frac{4}{\mathsf{c}^2}(\Phi_{\mathsf{N}}-\frac{W}{\mathsf{c}^2})}
\Big(1+\frac{X}{\mathsf{c}^4}\Big)^2\Big) +  \\
&-\Omega\varpi^2Y-\frac{\Omega^2\varpi^4}{2\mathsf{c}^4}Y^2
-\frac{\mathsf{c}^4}{2}\sum_{k=2}^{\infty}\frac{(-1)^{k+1}}{k}\Big(\frac{Z}{\mathsf{c}^2}\Big)^k
\end{align*}
with $$
Z=\frac{2\Omega\varpi^2}{\mathsf{c}^2}Y+
\frac{\Omega^2\varpi^4}{\mathsf{c}^6}Y^2
-\Omega^2\varpi^2e^{-\frac{4}{\mathsf{c}^2}(\Phi_{\mathsf{N}}-\frac{W}{\mathsf{c}^2})}\Big(1+
\frac{X}{\mathsf{c}^4}\Big)^2,$$
or
$$
w=W-2\Omega^2\varpi^2\Phi_{\mathsf{N}}-\Omega\varpi^2Y+
\frac{1}{4}\Omega^4\varpi^4+\frac{1}{\mathsf{c}^2}\mathfrak{Q}_0,
$$
and
$F=\frac{1}{\mathsf{c}^2}\Big(\Phi_{\mathsf{N}}-\frac{W}{\mathsf{c}^2}\Big)$ and
$K=\frac{V}{\mathsf{c}^4}$.

 Note that $\mathfrak{Q}_5=\mathfrak{Q}_6=\mathfrak{R}_c=0$ when both $\rho_{\mathsf{N}}$ and $\rho$ vanish, that is, when $u_{\mathsf{N}}\leq 0$ and
$u_{\mathsf{N}}+\frac{1}{\mathsf{c}^2}w\leq 0$.\\

The equations for $V$ are
\begin{subequations}
\begin{align}
&\frac{\partial V}{\partial\varpi}=\frac{\varpi}{2}\Big[\frac{\partial^2}{\partial\varpi^2}+\frac{2}{\varpi}\frac{\partial}{\partial \varpi}-\frac{\partial^2}{\partial z^2}\Big]X+
\varpi\Big(\Big(\frac{\partial\Phi_{\mathsf{N}}}{\partial\varpi}\Big)^2
-\Big(\frac{\partial\Phi_{\mathsf{N}}}{\partial z}\Big)^2\Big)+\mathfrak{R}_d, \label{Veq.a}\\
&\frac{\partial V}{\partial z}=\Big[1+\varpi\frac{\partial}{\partial\varpi}\Big]
\frac{\partial X}{\partial z}+
2\varpi\frac{\partial\Phi_{\mathsf{N}}}{\partial\varpi}\frac{\partial\Phi_{\mathsf{N}}}{\partial z}+\mathfrak{R}_e. \label{Veq.b}
\end{align}
\end{subequations}

We define  $\mathfrak{R}_d, \mathfrak{R}_e$ by using the auxiliary quantities $\mathfrak{Q}_7, \mathfrak{Q}_8$ as follows:
\begin{subequations}
\begin{align}
\Big(1+\frac{X}{\mathsf{c}^4}&+\frac{\varpi}{\mathsf{c}^4}\frac{\partial X}{\partial \varpi}\Big)^{-1}\times\Bigg[\quad
\frac{1}{2}(2\partial_1X+\varpi\partial_1^2X-\varpi\partial_3^2X) + \nonumber \\
&+\varpi\Big(1+\frac{X}{\mathsf{c}^4}\Big)\Big(\Big(\partial_1(\Phi_{\mathsf{N}}-\frac{W}{\mathsf{c}^4})\Big)^2
-\Big(\partial_3(\Phi_{\mathsf{N}}-\frac{W}{\mathsf{c}^4})\Big)^2\Big) + \nonumber \\
&-\frac{e^{4F}}{4}\Big(1+\frac{X}{\mathsf{c}^4}\Big)^{-1}\frac{\varpi}{\mathsf{c}^2}
\Big((2Y+\varpi\partial_1Y)^2-(\varpi\partial_3Y)^2\Big) \quad \Big) + \nonumber \\
&+\varpi\frac{\partial X}{\partial z}\times \frac{1}{\mathsf{c}^4}\Big(\quad \partial_3X+
\varpi\partial_1\partial_3X + \nonumber  \\
&2\varpi\Big(1+\frac{X}{\mathsf{c}^4}\Big)\partial_1(\Phi_{\mathsf{N}}-\frac{W}{\mathsf{c}^2})
\partial_3(\Phi_{\mathsf{N}}-\frac{W}{\mathsf{c}^2}) + \nonumber \\
&-\frac{e^{4F}}{2}\Big(1+\frac{X}{\mathsf{c}^4}\Big)\frac{\varpi^2}{\mathsf{c}^2}(2Y+\partial_1Y)(\partial_3Y) \quad \Bigg] \nonumber \\
&=
\frac{1}{2}(2\partial_1X+\varpi\partial_1^2X-\varpi\partial_3^2X)+
\varpi((\partial_1\Phi_{\mathsf{N}})^2-(\partial_3\Phi_{\mathsf{N}})^2)
+\frac{1}{\mathsf{c}^2}\mathfrak{Q}_7, \\
\mathfrak{R}_d&:=\frac{\hat{X}}{\mathsf{c}^4}\Big(
\frac{1}{2}(2\partial_1X+\varpi\partial_1^2X-\varpi\partial_3^2X)+
\varpi((\partial_1\Phi_{\mathsf{N}})^2-(\partial_3\Phi_{\mathsf{N}})^2) \Big)+ \nonumber \\
&+\Big(1+\frac{\hat{X}}{\mathsf{c}^4}\Big)\frac{1}{\mathsf{c}^2}\mathfrak{Q}_7;
\end{align}
\end{subequations}
\begin{subequations}
\begin{align}
-\frac{\varpi}{\mathsf{c}^4}\frac{\partial X}{\partial z}&\times \Bigg[\quad \frac{1}{2}(2\partial_1X+\varpi\partial_1^2X-\partial_3^2X)+ \nonumber \\
&+\varpi\Big(1+\frac{X}{\mathsf{c}^4}\Big)
\Big((\partial_1(\Phi_{\mathsf{N}}-\frac{W}{\mathsf{c}^2}))^2-(\partial_3(\Phi_{\mathsf{N}}-\frac{W}{\mathsf{c}^2}))^2\Big) + \nonumber \\
&
-\frac{e^{4F}}{\mathsf{c}^2}\Big(1+\frac{X}{\mathsf{c}^4}\Big)^{-1}\varpi((2Y+\varpi\partial_1Y)^2-(\varpi\partial_3Y)^2) + \nonumber \\
&+\Big(1+\frac{X}{\mathsf{c}^4}+\frac{\varpi}{\mathsf{c}^4}\frac{\partial X}{\partial\varpi}\Big)\times\Big(\quad \partial_3X+\varpi\partial_1\partial_3X + \nonumber \\
&+2\varpi\Big(1+\frac{X}{\mathsf{c}^4}\Big(\partial_1(\Phi_{\mathsf{N}}-\frac{W}{\mathsf{c}^2})\partial_3(\Phi_{\mathsf{N}}-\frac{W}{\mathsf{c}^2}) + \nonumber \\
&-\frac{e^{4F}}{2\mathsf{c}^2}\Big(1+\frac{X}{\mathsf{c}^4}\Big)^{-1}
\varpi(2Y+\varpi\partial_1Y)(\varpi\partial_3Y)\quad \Bigg] \nonumber \\
&=
\partial_3X+\varpi\partial_1\partial_3X+2\varpi(\partial_1\Phi_{\mathsf{N}})(\partial_3\Phi_{\mathsf{N}})+
\frac{1}{\mathsf{c}^2}\mathfrak{Q}_8, \\
\mathfrak{R}_e&:=\frac{\hat{X}}{\mathsf{c}^4}\Big(
\partial_3X+\varpi\partial_1\partial_3X+2\varpi(\partial_1\Phi_{\mathsf{N}})(\partial_3\Phi_{\mathsf{N}})
\Big) + \nonumber \\
&+\Big(1+\frac{\hat{X}}{\mathsf{c}^4}\Big)\frac{1}{\mathsf{c}^2}\mathfrak{Q}_8.
\end{align}
\end{subequations}
Here we put
\begin{equation}
\frac{\hat{X}}{\mathsf{c}^4}:=
\Big(\Big(1+\frac{X}{\mathsf{c}^4}+\frac{\varpi}{\mathsf{c}^4}\frac{\partial X}{\partial\varpi}\Big)^2+
\Big(\frac{\varpi}{\mathsf{c}^4}\frac{\partial X}{\partial z}\Big)^2\Big)^{-1}-1.
\end{equation}\\

\subsection{Functional spaces $\mathfrak{C}_{(n)}^l, \mathfrak{C}_{(n)}^{l,\alpha}$}

Let us introduce the following 

\begin{Definition}
Let $n=3,4,5$.
The Kelvin transformation of the co-ordinates $\mbox{\boldmath$\xi$} \mapsto \mbox{\boldmath$\xi$}^{\star}: \mathbb{R}^n
\rightarrow \mathbb{R}^n$ is defined by
\begin{equation}
\mbox{\boldmath$\xi$}^{\star}=
\Big(\frac{\Xi_0}{|\mbox{\boldmath$\xi$}|}\Big)^2\mbox{\boldmath$\xi$}.
\end{equation}
The Kelvin transformation $f_{\star(n)}$ of a function of $\mbox{\boldmath$\xi$}$ is defined by
\begin{equation}
f_{\star(n)}(\mbox{\boldmath$\xi$}^{\star})=\Big(\frac{|\mbox{\boldmath$\xi$}|}{\Xi_0}\Big)^{n-2}f(\mbox{\boldmath$\xi$})
=\Big(\frac{\Xi_0}{|\mbox{\boldmath$\xi$}^{\star}|}\Big)^{n-2}
f\Big(
\Big(\frac{\Xi_0}{|\mbox{\boldmath$\xi$}^{\star}|}\Big)^2\mbox{\boldmath$\xi$}^{\star}\Big).
\end{equation}
\end{Definition}

Then we have $|\mbox{\boldmath$\xi$}||\mbox{\boldmath$\xi$}^{\star}|=\Xi_0^2$ and the Kelvin transformation maps $\mathbb{R}^n\setminus B^{(n)}(\Xi_0)$ onto $\bar{B}^{(n)}(\Xi_0)\setminus\{\mathbf{0}\}$.
Note that $ (\mbox{\boldmath$\xi$}^{\star})^{\star}=\mbox{\boldmath$\xi$}$ and $(f_{{\star}(n)})_{{\star}(n)}=f$.

 It is well known that
\begin{equation}
\triangle_{\mbox{\boldmath$\xi$}^{\star}}^{(n)}f_{\star(n)}(\mbox{\boldmath$\xi$}^{\star})=
\Big(\frac{|\mbox{\boldmath$\xi$}|}{\Xi_0}\Big)^{n+2}
\triangle_{\mbox{\boldmath$\xi$}}^{(n)}f(\mbox{\boldmath$\xi$}).
\end{equation}
Actually we have
\begin{equation}
\frac{\partial}{\partial\xi^{\star}_j}f_{\star(n)}(\mbox{\boldmath$\xi$}^{\star})=
\frac{|\mbox{\boldmath$\xi$}|^{n-2}}{(\Xi_0)^n}\Big[
-(n-2)\xi_jf(\mbox{\boldmath$\xi$})+
\sum_k(|\mbox{\boldmath$\xi$}|^2\delta_{jk}-2\xi_j\xi_k)\frac{\partial}{\partial \xi_k}
f(\mbox{\boldmath$\xi$})\Big], \label{D*}
\end{equation}
and
\begin{align}
&\frac{\partial^2}{\partial \xi^{\star}_i\partial \xi^{\star}_j}f_{\star(n)}(\mbox{\boldmath$\xi$}^{\star})=
\frac{|\mbox{\boldmath$\xi$}|^{n-2}}{(\Xi_0)^{n+2}}\Big[\quad -(n-2)(|\mbox{\boldmath$\xi$}|^2\delta_{ij}-
n\xi_i\xi_j)f(\mbox{\boldmath$\xi$}) + \nonumber \\
&-n|\mbox{\boldmath$\xi$}|^2(\xi_j\frac{\partial}{\partial \xi_i}f(\mbox{\boldmath$\xi$})+
\xi_i\frac{\partial}{\partial \xi_j}f(\mbox{\boldmath$\xi$}))
-2(|\mbox{\boldmath$\xi$}|^2\delta_{ij}-2n\xi_i\xi_j)
\sum_k
\xi_k\frac{\partial}{\partial \xi_k}f (\mbox{\boldmath$\xi$})+ \nonumber \\
&\sum_{k,l}
(|\mbox{\boldmath$\xi$}|^2\delta_{jk}-2\xi_j\xi_k)
(|\mbox{\boldmath$\xi$}|^2\delta_{il}-2\xi_i\xi_l)
\frac{\partial^2}{\partial \xi_l\partial \xi_k}f (\mbox{\boldmath$\xi$})\qquad \Big]. \label{DD*}
\end{align}\\

We use the following
\begin{Definition}
We denote by $\mathfrak{C}_{(n)}^l(\overline{\mathfrak{D}^{\mathsf{c}}}(R))$
the set of all functions $Q$ of $(\varpi, z)$ such that
$Q(\varpi, -z)=Q(\varpi, z)$ for $\forall \varpi, z$ and
$$(Q^{\flat(n)})_{\star(n)}\in C^l(\bar{B}^{(n)}(\Xi^{\star})) \upharpoonright \bar{B}^{(n)}(\Xi^{\star})\setminus\{\mathbf{0}\}.$$ Here $l=0,1,2, R=\mathsf{a}\Xi$ and $\Xi^{\star}=(\Xi_0)^2/\Xi$.. We put
$$\|Q;\mathfrak{C}_{(n)}^l(\overline{\mathfrak{D}^{\mathsf{c}}}(R))\|:=\|F;C^l(\bar{B}(\Xi^{\star}))\| $$
by $F$ such that $(Q^{\flat(n)})_{\star(n)}=F \upharpoonright \bar{B}^{(n)}(\Xi^{\star})\setminus\{\mathbf{0}\}$.

We denote by $\mathfrak{C}_{(n)}^{l,\alpha}(\overline{\mathfrak{D}^{\mathsf{c}}}(R))$
the set of all functions $Q$ of $(\varpi, z)$ such that
$Q(\varpi, -z)=Q(\varpi, z)$ for $\forall \varpi, z$ and
$$(Q^{\flat(n)})_{\star(n)}\in C^{l,\alpha}(\bar{B}^{(n)}(\Xi^{\star})) \upharpoonright \bar{B}^{(n)}(\Xi^{\star})\setminus\{\mathbf{0}\}.$$ Here $l=0,1,2, R=\mathsf{a}\Xi$ and $\alpha$ is the fixed number satisfying  \eqref{alpha}. We put
$$\|Q;\mathfrak{C}_{(n)}^{l,\alpha}(\overline{\mathfrak{D}^{\mathsf{c}}}(R))\|:=\|F;C^{l,\alpha}(\bar{B}(\Xi^{\star}))\| $$
by $F$ such that $(Q^{\flat(n)})_{\star(n)}=F \upharpoonright \bar{B}^{(n)}(\Xi^{\star})\setminus\{\mathbf{0}\}$.
\end{Definition}

Note that, when
$Q \in \mathfrak{C}_{(n)}^{0}(\overline{\mathfrak{D}^{\mathsf{c}}}(R))$,
 the limit
$$C_{\infty}=\lim_{\mbox{\boldmath$\xi$}^{\star}\rightarrow \mathbf{0}}(Q^{\flat(n)})_{\star(n)}(\mbox{\boldmath$\xi$}^{\star})$$
exists by the definition, and then
$$Q=\Big(\frac{R_0}{r}\Big)^{n-2}\Big[C_{\infty}+O\Big(\frac{1}{r}\Big)\Big]\quad\mbox{as}\quad r
\rightarrow +\infty,$$
provided that $Q \in \mathfrak{C}_{(n)}^{1}(\overline{\mathfrak{D}^{\mathsf{c}}}(R))$.\\

The removable singularity theorem of harmonic functions (see \cite[p.269, Corollary to Chapter X, Theorem XI]{Kellog}, e.g., ) tells us 

\begin{Lemma} \label{removable}
If $f \in C^2(\bar{B}^{(n)}(\Xi)\setminus\{\mathbf{0}\}), \triangle f=g \in C^{0,\alpha}(\bar{B}^{(n)}(\Xi))
\upharpoonright \bar{B}^{(n)}(\Xi)\setminus\{\mathbf{0}\}$, and if
$f(\mbox{\boldmath$\xi$})=O(1)$ as $\mbox{\boldmath$\xi$} \rightarrow \mathbf{0}$, then 
$f \in C^{2,\alpha}(\bar{B}^{(n)}(\Xi))\upharpoonright \bar{B}^{(n)}(\Xi)\setminus\{\mathbf{0}\}$.
Moreover suppose  $g=0$. Then $f$ is given by the Poisson integral
$$f(\mbox{\boldmath$\xi$})=\frac{\Xi^2-|\mbox{\boldmath$\xi$}|^2}{\Sigma_n\Xi}\int_{\partial B^{(n)}(\Xi)}
\frac{f(\mbox{\boldmath$\xi$}')}{|\mbox{\boldmath$\xi$}-\mbox{\boldmath$\xi$}'|^n}d\sigma(\mbox{\boldmath$\xi$}'),$$
where $d\sigma$ is the volume element of $\partial B^{(n)}$, $\Sigma_n=2(n-2)\pi^{n/2}/\Gamma(n/2)$,
and it holds that 
\begin{align}
&\sup_{\bar{B}^{(n)}(\Xi))}|f| \leq \sup_{\partial B^{(n)}(\Xi))}|f|,\quad
\sup_{\bar{B}^{(n)}(\Xi))}|\partial_jf| \leq \sup_{\partial B^{(n)}(\Xi))}|\partial_jf|,(j=1,\cdots,n)
, \label{Pineq.1}\\
&\|f;C^{2,\alpha}(\bar{B}^{(n)}(\Xi))\| \leq C\|f\restriction \partial B^{(n)}(\Xi);C^{2,\alpha}(\partial B^{(n)}(\Xi))\|. \label{Pineq.2}
\end{align}
\end{Lemma}

Here, for a function $h$ defined on $\partial B^{(n)}(\Xi)$, $\|h; C^{2,\alpha}(\partial B^{(n)}(\Xi))\|$ means
$$\inf \Big\{ \|H;C^{2,\alpha}(\bar{B}^{(n)}(\Xi))\|\quad\Big|\quad H\restriction
\partial B^{(n)}(\Xi) =h \Big\}. $$
Clearly, for any $0<\delta<1$ we have
\begin{align*}
\|f\restriction \partial B^{(n)}(\Xi);C^{2,\alpha}(\partial B^{(n)}(\Xi))\|&\leq
\|\chi\Big(\frac{2}{\delta}\frac{\Xi-|\cdot|}{\Xi}\Big)f; C^{2,\alpha}(\bar{B}^{(n)}(\Xi))\| \\
&\leq C_{\delta}\Big\|f;C^{2,\alpha}\Big(\bar{B}^{(n)}(\Xi) \setminus B^{(n)}((1-\delta)\Xi)\Big)\Big\|
,
\end{align*}
where $C_{\delta}$ is a constant depending on $\delta$ only for fixed $\Xi, \alpha, n$ and the cut off function $\chi$. Therefore \eqref{Pineq.2} says
$$\|f;C^{2,\alpha}(\bar{B}^{(n)}(\Xi))\| \leq 
 C_{\delta}'\Big\|f;C^{2,\alpha}\Big(\bar{B}^{(n)}(\Xi) \setminus B^{(n)}((1-\delta)\Xi)\Big)\Big\|
$$
for any $0<\delta<1$.\\

Applying  Lemma \ref{removable} to the Kelvin transformed function, we can claim

\begin{Proposition}
If $ Q \in \mathfrak{C}_{(n)}^2(\mathfrak{D}(R)^{\mathsf{c}})$, 
$$\Big[\frac{\partial^2}{\partial\varpi^2}+\frac{n-1}{\varpi}\frac{\partial}{\partial \varpi}
+\frac{\partial^2}{\partial z^2}\Big]Q \in
\mathfrak{C}_{(n)}^{0,\alpha}(\overline{\mathfrak{D}^{\mathsf{c}}}(R)), $$
and if $Q=O(r^{-(n-2)})$ as $ r \rightarrow +\infty$, then
$Q \in
\mathfrak{C}_{(n)}^{2,\alpha}(\overline{\mathfrak{D}^{\mathsf{c}}}(R))$
and 
$$Q=\Big(\frac{R_0}{r}\Big)^{n-2}\Big[C_{\infty}+O\Big(\frac{1}{r}\Big)\Big]
\quad\mbox{as}\quad
r \rightarrow +\infty$$
with a constant $C_{\infty}$.
\end{Proposition}

Let us use again the cut off function $\chi \in C^{\infty}(\mathbb{R})$ such that
$\chi(t)=1$ for $t \leq 1$, $0<\chi(t)<1$ for $1<t<2$ and
$\chi(t)=0$ for $2\geq t$. We use the following

\begin{Definition}
1) For any continuous function $f$ on $\mathbf{R}^n$, we denote
\begin{subequations}
\begin{align}
&f^{[0]}(\mbox{\boldmath$\xi$}):=\chi(|\mbox{\boldmath$\xi$}|/\Xi_0)f(\mbox{\boldmath$\xi$}), \\
&f^{[\infty]}(\mbox{\boldmath$\xi$}):=(1-\chi(|\mbox{\boldmath$\xi$}|/\Xi_0))f(\mbox{\boldmath$\xi$})).
\end{align}
\end{subequations}
 For any continuous function $Q$ of $(\varpi, z)$ we put
\begin{equation}
Q^{[0]}=\chi(r/R_0)Q,\qquad
Q^{[\infty]}=(1-\chi(r/R_0))Q,
\end{equation}
while $R_0=4r_1=\mathsf{a}\Xi_0=4\mathsf{a}\xi_1(\frac{1}{\gamma-1})$. 

2) We denote by $\mathfrak{C}_{(n)}^{l }$ the set of all functions $Q$ of $(\varpi, z)$ such that $Q(\varpi, -z)=Q(\varpi, z)$ and
$$ Q^{[0]}\in \mathfrak{C}^l(\bar{\mathfrak{D}}(2R_0))\quad\mbox{and}\quad
Q^{[\infty]}\in \mathfrak{C}_{(n)}^{l }(\overline{\mathfrak{D}^{\mathsf{c}}}(R_0)).
$$ 
We put
$$\|Q;\mathfrak{C}_{(n)}^{l  }\|=
 \|Q^{[0]}; \mathfrak{C}^{l  }(\bar{\mathfrak{D}}(2R_0))
\| \vee \|Q^{[\infty]};\mathfrak{C}_{(n)}^{l  }(\overline{\mathfrak{D}^{\mathsf{c}}}(R_0))
\| .
$$

We denote by $\mathfrak{C}_{(n)}^{l,\alpha }$ the set of all functions $Q$ of $(\varpi, z)$ such that $Q(\varpi, -z)=Q(\varpi, z)$ and
$$ Q^{[0]}\in \mathfrak{C}^{l,\alpha}(\bar{\mathfrak{D}}(2R_0))\quad\mbox{and}\quad
Q^{[\infty]}\in \mathfrak{C}_{(n)}^{l,\alpha }(\overline{\mathfrak{D}^{\mathsf{c}}}(R_0)).
$$
We put
$$\|Q;\mathfrak{C}_{(n)}^{l,\alpha  }\|=
 \|Q^{[0]}; \mathfrak{C}^{l,\alpha  }(\bar{\mathfrak{D}}(2R_0))
\| \vee \|Q^{[\infty]};\mathfrak{C}_{(n)}^{l,\alpha }
(\overline{\mathfrak{D}^{\mathsf{c}}}(R_0))
\| .
$$
\end{Definition}

Let us prepare some calculus rules which will be used later about the functional spaces $\mathfrak{C}_{(n)}^{l}, \mathfrak{C}_{(n)}^{l,\alpha}, l=0,1$.

\begin{Proposition}\label{Prop.DExt}
Let $n=3,4,5$. There is a constant $C$ such that, for $j=1,3$ with
$\displaystyle \partial_1={\partial}/{\partial\varpi}, \partial_3={\partial}/{\partial z}$, 
\begin{align}
&\|\partial_jQ;\mathfrak{C}_{(n)}^0\|\leq C
\mathsf{a}^{-1}
\|Q;  \mathfrak{C}_{(n)}^{1}\|, \label{DExt.a}\\
&\|\partial_jQ; 
 \mathfrak{C}_{(n)}^{0,\alpha}\|\leq C
\mathsf{a}^{-1}
\|Q;  \mathfrak{C}_{(n)}^{1,\alpha}\|. \label{DExt.b}
\end{align}
\end{Proposition}

Proof. Since 
$$(\partial_jQ)^{[0]}=\partial_j(Q^{[0]})-(\partial_j\chi(r/R_0))\cdot Q, $$
we see
\begin{align*}
\|(\partial_jQ)^{[0]};\mathfrak{C}^0(\bar{\mathfrak{D}}(2R_0))\| &\leq 
 \sup_{r\leq R_0}|\partial_j(Q^{[0]})| +
\sup_{r\leq R_0}|(\partial_j\chi)\cdot Q| \\
&\leq\mathsf{a}^{-1}\Big( \|Q;\mathfrak{C}^1\|+
C\sup_{R_0\leq r\leq 2R_0}|Q| \Big) \\
&\leq \mathsf{a}^{-1}\Big(\|Q;\mathfrak{C}^1\|
+C\Big(\sup_{r\leq 2R_0}|Q^{[0]}|+\sup_{R_0\leq r}|Q^{[\infty]}|\Big) \Big) \\
&\leq\mathsf{a}^{-1}\Big( \|Q;\mathfrak{C}^1\| +2C
\|Q;\mathfrak{C}^0\| \Big),
\end{align*}
for
\begin{align*}
\sup_{R_0\leq r}|Q^{[\infty]}| &=\sup_{R_0\leq r}\Big|\Big(\frac{R_0}{r}\Big)^{n-2}
Q^{[\infty]\flat(n)}_{\star(n)}\Big| 
\leq \sup_{R_0\leq r}|Q^{[\infty]\flat(n)}_{\star(n)}| \\
&\leq \|Q^{[\infty]}; \mathfrak{C}^0(\overline{\mathfrak{D}^{\mathsf{c}}}(R_0))\| 
\leq \|Q; \mathfrak{C}^0\|.
\end{align*}

$\|(\partial_jQ)^{[\infty]}; \mathfrak{C}_{(n)}^0(\overline{\mathfrak{D}^{\mathsf{c}}}(R_0))\|$ can be estimated similarly, namely, it is sufficient to estimate
$$\|\partial_j(Q^{[\infty]});\mathfrak{C}_{(n)}^0(\overline{\mathfrak{D}^{\mathsf{c}}}(R_0))\|=
\sup_{R_0\leq r}
\Big|\Big(\frac{r}{R_0}\Big)^{n-2}\partial_j(Q^{[\infty]})\Big|.$$
But \eqref{D*} says, for $|\mbox{\boldmath$\xi$}^{\star}|\leq \Xi_0 (\Leftrightarrow |\mbox{\boldmath$\xi$}|\geq \Xi_0 )$, 
$$
\Big(\frac{|\mbox{\boldmath$\xi$}|}{\Xi_0}\Big)^{n-2}
|D_{\mbox{\boldmath$\xi$}}f(\mbox{\boldmath$\xi$})|\leq C \Bigg[
|f_{\star(n)}(\mbox{\boldmath$\xi$}^{\star})|+|D_{\mbox{\boldmath$\xi$}^{\star}}f_{\star(n)}(\mbox{\boldmath$\xi$}^{\star})|\Bigg] ,
$$
where 
the symbol $|D_XF(X)|$ stands for
$\displaystyle \max_k\Big|\frac{\partial}{\partial X^k}F(X^1,\cdots, X^n)\Big|$.
Here we have read \eqref{D*} transposing $\mbox{\boldmath$\xi$}, f $ and $\mbox{\boldmath$\xi$}^{\star}, f_{\star(n)}$, by keeping in mind that $(\mbox{\boldmath$\xi$}^{\star})^{\star}=\mbox{\boldmath$\xi$}, (f_{\star(n)})_{\star(n)}=f$, namely we read
$$
\frac{\partial}{\partial\xi_j}f(\mbox{\boldmath$\xi$})=
\frac{|\mbox{\boldmath$\xi$}^{\star}|^{n-2}}{(\Xi_0)^n}\Big[
-(n-2)\xi_j^{\star}f_{\star(n)}(\mbox{\boldmath$\xi$}^{\star})+
\sum_k(|\mbox{\boldmath$\xi$}^{\star}|^2\delta_{jk}-2\xi^{\star}_j\xi^{\star}_k)\frac{\partial}{\partial \xi^{\star}_k}
f_{\star(n)}(\mbox{\boldmath$\xi$}^{\star})\Big].
$$
Thus we see
\begin{align*}
\|\partial_j(Q^{[\infty]});\mathfrak{C}_{(n)}^0(\overline{\mathfrak{D}^{\mathsf{c}}}(R_0))\|&=
\sup_{R_0\leq r}
\Big|\Big(\frac{r}{R_0}\Big)^{n-2}\partial_j(Q^{[\infty]})\Big| \\
&\leq C \mathsf{a}^{-1} \|Q^{[\infty]}; \mathfrak{C}_{(n)}^1(\overline{\mathfrak{D}^{\mathsf{c}}}(R_0))\|.
\end{align*}
Summing up,  we have the estimate \eqref{DExt.a} of
$$
\|\partial_jQ;\mathfrak{C}_{(n)}^0\|= 
\|(\partial_jQ)^{[0]};\mathfrak{C}^0(\bar{\mathfrak{D}}(2R_0))\| \vee
\|(\partial_jQ)^{[\infty]};\mathfrak{C}_{(n)}^0(\overline{\mathfrak{D}^{\mathsf{c}}}(R_0))\|.
$$

In order to estimate $\|(\partial_jQ)^{[0]};\mathfrak{C}^{0,\alpha}(\bar{\mathfrak{D}}(2R_0))\|$ or
$$
|\Delta (\partial_jQ)^{[0]}|=
|(\partial_jQ)^{[0]}(\varpi', z')-(\partial_jQ)^{[0]}(\varpi,z)|,$$
it is sufficient to estimate
$$|\Delta \partial_j(Q^{[0]})|=
|\partial_j(Q^{[0]})(\varpi', z')-\partial_j(Q^{[0]})(\varpi, z)|$$
and
\begin{align*}
|\Delta Q|&=|Q(\varpi', z')-Q(\varpi, z)| \\
&\leq |Q^{[0]}(\varpi',z')-Q^{[0]}(\varpi, z)| + |Q^{[\infty]}(\varpi',z')-Q^{[\infty]}(\varpi, z)|.
\end{align*}
Of course
$$|\Delta\partial_jQ^{[0]}|\leq \|Q;\mathfrak{C}^{1,\alpha}(\bar{\mathfrak{D}}(2R_0))\|\cdot |\mbox{\boldmath$\xi$}'-\mbox{\boldmath$\xi$}|^{\alpha}
$$
and
$$
|\Delta Q^{[0]}|\leq \|Q;\mathfrak{C}^{0,\alpha}(\bar{\mathfrak{D}}(2R_0))\|\cdot|\mbox{\boldmath$\xi$}'-\mbox{\boldmath$\xi$}|^{\alpha}$$
by definition. As for the estimate of
$$
|\Delta Q^{[\infty]}|=\Big|\Delta \Big(\frac{R_0}{r}\Big)^{n-2}(Q^{[\infty]})_{\star(n)}^{\flat(n)}\Big| $$
for $\Xi_0 \leq |\mbox{\boldmath$\xi$}'|\wedge |\mbox{\boldmath$\xi$}|\leq
|\mbox{\boldmath$\xi$}'| \vee |\xi|\leq 2\Xi_0$, where 
$r'=\sqrt{(\varpi')^2+(z')^2}=|\mbox{\boldmath$\xi$}'|, r=\sqrt{\varpi^2+z^2}=|\mbox{\boldmath$\xi$}| $, we note that
$$ 
\Big|\Big(\frac{\Xi_0}{|\mbox{\boldmath$\xi$}'|}\Big)^{n-2}-\Big(\frac{\Xi_0}{|\mbox{\boldmath$\xi$}|}\Big)^{n-2}\Big|
\leq C |\mbox{\boldmath$\xi$}'-\mbox{\boldmath$\xi$}|\leq C'
|\mbox{\boldmath$\xi$}'-\mbox{\boldmath$\xi$}|^{\alpha}
$$
and
$ |(\mbox{\boldmath$\xi$}')^{\star}-\mbox{\boldmath$\xi$}^{\star}|^{\alpha} \leq C
|\mbox{\boldmath$\xi$}'-\mbox{\boldmath$\xi$}|^{\alpha}$. This observation deduces 
$$|\Delta Q^{[\infty]}| \leq C \|Q;\mathfrak{C}^{0,\alpha}(\bar{\mathfrak{D}}(2R_0))\|\cdot
|\mbox{\boldmath$\xi$}'-\mbox{\boldmath$\xi$}|^{\alpha},$$
and, therefore, 
$$
\|(\partial_jQ)^{[0]}; \mathfrak{C}^{0,\alpha}(\bar{\mathfrak{D}}(2R_0))\|
\leq C \mathsf{a}^{-1}
 \|Q;\mathfrak{C}^{1,\alpha}(\bar{\mathfrak{D}}(2R_0))\|.$$

In order to estimate
$\|(\partial_jQ)^{[\infty]};\mathfrak{C}^{0,\alpha}_{(n)}(\overline{\mathfrak{D}^{\mathsf{c}}}(R_0))\|$ or
$$\Big|\Delta 
\Big(\frac{r}{R_0}\Big)^{n-2}(\partial_jQ)^{[\infty]}\Big|=\Big|\Big(\frac{r'}{R_0}\Big)^{n-2}(\partial_jQ)^{[\infty]}(\varpi',z')-\Big(\frac{r}{R_0}\Big)^{n-2}(\partial_jQ)^{[\infty]}(\varpi,z)\Big|,$$
it is sufficient to estimate
$\displaystyle \Big|\Delta \Big(\frac{r}{R_0}\Big)^{n-2}\partial_j(Q^{[\infty]})\Big| $
on $R_0 \leq r'\wedge r$. But
 \eqref{D*} says, for 
$|\mbox{\boldmath$\xi$}'|\wedge |\mbox{\boldmath$\xi$}|\geq \Xi_0
 (\Leftrightarrow |(\mbox{\boldmath$\xi$}')^{\star}|\vee |\mbox{\boldmath$\xi$}^{\star}| \leq \Xi_0)$,
\begin{align*}
\Bigg|\Delta_{\mbox{\boldmath$\xi$}}\Bigg[\Big(\frac{|\mbox{\boldmath$\xi$}|}{\Xi_0}\Big)^{n-2}D_{\mbox{\boldmath$\xi$}}f(\mbox{\boldmath$\xi$})\Bigg]\Bigg| 
\leq C &\Bigg[\sup_{|\mbox{\boldmath$\xi$}_-^{\star}|\leq \Xi_0}\Big(|f_{\star(n)}(\mbox{\boldmath$\xi$}_-^{\star})|+
|D_{\mbox{\boldmath$\xi$}_-^{\star}}f_{\star(n)}(\mbox{\boldmath$\xi$}_-^{\star})|\Big)
|\Delta_{\mbox{\boldmath$\xi$}^{\star}}\mbox{\boldmath$\xi$}^{\star}| + \\
&+
\Big|\Delta_{\mbox{\boldmath$\xi$}^{\star}}f_{\star(n)}(\mbox{\boldmath$\xi$}^{\star})\Big|+
\Big|\Delta_{\mbox{\boldmath$\xi$}^{\star}}D_{\xi^{\star}}f_{\star(n)}(\mbox{\boldmath$\xi$}^{\star})\Big|
\Bigg],
\end{align*}
where
$\Delta_XF(X_1,\cdots,X_n)$ stands for
$F(X_1', \cdots. X_n')-F(X_1, \cdots, X_n)$.
Therefore we have
$$
\Big|\Delta\Big(\frac{r}{R_0}\Big)^{n-2}\partial_j(Q^{[\infty]})\Big|\leq C
\|Q^{[\infty]};\mathfrak{C}^{1,\alpha}_{(n)}(\overline{\mathfrak{D}^{\mathsf{c}}}(R_0))\|\cdot
|(\mbox{\boldmath$\xi$}')^{\star}-\mbox{\boldmath$\xi$}^{\star}|^{\alpha} 
$$
for $R_0 \leq r'\wedge r$.
 $\square$\\

\begin{Proposition}\label{Prop.Product}
Let $n, m \geq 3$.
 If $Q_1\in \mathfrak{C}_{(m)}^{0}$ and $ Q_2 \in \mathfrak{C}_{(n)}^{0}$, then $Q_1\cdot Q_2 \in \mathfrak{C}_{(n)}^{0}$ and $\|Q_1\cdot Q_2;\mathfrak{C}_{(n)}^{0}\|\leq C
\|Q_1; \mathfrak{C}_{(m)}^{0}\|
\cdot
\|Q_2; \mathfrak{C}_{(n)}^{0}\|$.
 If $Q_1\in \mathfrak{C}_{(m)}^{0,\alpha}$ and $ Q_2 \in \mathfrak{C}_{(n)}^{0,\alpha}$, then $Q_1\cdot Q_2 \in \mathfrak{C}_{(n)}^{0,\alpha}$ and $\|Q_1\cdot Q_2;\mathfrak{C}_{(n)}^{0,\alpha}\|\leq C
\|Q_1; \mathfrak{C}_{(m)}^{0,\alpha}\|
\cdot
\|Q_2; \mathfrak{C}_{(n)}^{0,\alpha}\|$.
\end{Proposition}

In order to prove this Proposition it is convenient to introduce the Kelvin transformation of the co-ordinate $(\varpi, z)$ and the Kelvin transformation of function of $(\varpi, z)$ as follows:
\begin{Definition}
Let us denote by $\mathbf{p}$ the vector $(\varpi, z)^{\top}$. The Kelvin transformation
$\mathbf{p} \mapsto \mathbf{p}^{\star}$ is defined by
\begin{equation}
\mathbf{p}^{\star}=\Big(\frac{R_0}{r}\Big)^2\mathbf{p},\quad\mbox{where}\quad
r=|\mathbf{p}|=\sqrt{\varpi^2+z^2}.
\end{equation}
The n-dimensional Kelvin transformation $Q_{{\star}(n)}$ of a function $Q$ of $(\varpi, z)$ is defined by
\begin{equation}
Q_{{\star}(n)}(\mathbf{p}^{\star})=\Big(\frac{r}{R_0}\Big)^{n-2}Q(\mathbf{p}).
\end{equation}
{\upshape Then we have 
$$(Q_{{\star}(n)})^{\flat(n)}=(Q^{\flat(n)})_{{\star}(n)}.$$
}
\end{Definition}

Using this notation, we can claim

\begin{Proposition}\label{Prop.leq}
If $Q \in \mathfrak{C}_{(m)}^{0}$ and $m \geq n$, then $Q \in \mathfrak{C}_{(n)}^{0}$.
If $Q \in \mathfrak{C}_{(m)}^{0,\alpha}$ and $m \geq n$, then $Q \in \mathfrak{C}_{(n)}^{0,\alpha}$.
\end{Proposition}

Proof. In fact we have
$$Q_{{\star}(n)}=\Big(\frac{r^{\star}}{R_0}\Big)^{m-n}Q_{{\star}(m)}, $$
where $\displaystyle r^{\star}=|\mathbf{p}^{\star}|=\frac{(R_0)^2}{r}$.
$\square$\\

Anyway, Proof of Proposition \ref{Prop.Product} is  given by the identity 
$$(Q_1Q_2)_{{\star}(n)}(\mathbf{p}^{\star})=\Big(\frac{|\mathbf{p}^{\star}|}{R_0}\Big)^{m-2}
(Q_1)_{{\star}(m)}(\mathbf{p}^{\star})(Q_2)_{{\star}(n)}(\mathbf{p}^{\star})$$
with $m-2 \geq 1$ as follows:

Proof of Proposition \ref{Prop.Product}.
 Let $m,n \geq 3$. Let $Q_1 \in \mathfrak{C}_{(m)}^0, Q_2 \in \mathfrak{C}_{(n)}^0$. We see
\begin{align*}
&\|(Q_1Q_2)^{[0]};\mathfrak{C}^0(\bar{\mathfrak{D}}(2R_0))\|=\sup_{r\leq 2R_0}|\chi(r/R_0)Q_1Q_2|=\sup_{r\leq 2R_0}|Q_1^{[0]}Q_2| \\
&\leq \sup_{r\leq 2R_0}|Q_1^{[0]}Q_2^{[0]}|+\sup_{r\leq 2R_0}|Q_1^{[0]}Q_2^{[\infty]}| \\
&\leq \sup_{r\leq 2R_0}|Q_1^{[0]}|\cdot\sup_{r\leq 2R_0}|Q_2^{[0]}|+\sup_{r\leq 2R_0}|Q_1^{[0]}|\cdot\sup_{R_0\leq r\leq 2R_0}|(Q_2^{[\infty]})_{\star(n)}| \\
&\leq 2\|Q_1;\mathfrak{C}_{(m)}^0\|\cdot\|Q_2;\mathfrak{C}_{(n)}^0\|.
\end{align*}
Here note that $$|Q^{[\infty]}|=
\Big|\Big(\frac{R_0}{r}\Big)^{n-2}Q^{[\infty]}_{\star(n)}\Big|\leq |Q^{[\infty]}_{\star(n))}|
\leq \|Q; \mathfrak{C}_{(n)}^0\|,$$
since either  $r \geq R_0$ or $Q^{[\infty]}=0$..

On the other hand we have
\begin{align*}
&\|(Q_1Q_2)^{[\infty]};\mathfrak{C}_{(n)}^0(\overline{\mathfrak{D}^{\mathsf{c}}}(R_0))\|=
\sup_{R_0\leq r}\Big|\Big(\frac{r}{R_0}\Big)^{n-2}(1-\chi)Q_1Q_2\Big| \\
&\leq \sup_{R_0\leq r}\Big|\Big(\frac{r}{R_0}\Big)^{n-2}Q_1^{[\infty]}Q_2^{[\infty]}\Big| 
+\sup_{R_0\leq r \leq 2R_0}
\Big|\Big(\frac{r}{R_0}\Big)^{n-2}Q_1^{\infty]}Q_2^{[0]}\Big| \\
&\leq \sup_{R_0\leq r}\Big|\Big(\frac{r}{R_0}\Big)^{m-2}Q_1^{[\infty]}
\Big(\frac{r}{R_0}\Big)^{n-2}Q_2^{[\infty]}\Big|+
\sup_{R_0\leq r \leq 2R_0}\Big|
\Big(\frac{r}{R_0}\Big)^{n-m}
\Big(\frac{r}{R_0}\Big)^{m-2}Q_1^{[\infty]}Q_2^{[0]}\Big| \\
&\leq (1+1\vee 2^{n-m})
\|Q_1;\mathfrak{C}_{(m)}^0\|\cdot\|Q_2;\mathfrak{C}_{(n)}^0\|.
\end{align*}
Summing up, we have
$$\|Q_1Q_2;\mathfrak{C}_{(n)}^0\|\leq C
\|Q_1;\mathfrak{C}_{(m)}^0\|\cdot\|Q_2;\mathfrak{C}_{(n)}^0\|. $$

Let $Q_1 \in \mathfrak{C}_{(m)}^{0,\alpha}, Q_2 \in \mathfrak{C}_{n)}^{0,\alpha}$. We shall use the notation
$$\Delta Q=Q(\varpi', z')-Q(\varpi, z).$$
Note that 
$$\Delta(S\cdot T) =(\Delta S)T'+S( \Delta T).$$

For $r'\vee r \leq 2R_0$, we have
\begin{align*}
&\Delta (Q_1Q_2)^{[0]} = \Delta(\chi Q_1Q_2) = \\
&=\Delta(\chi Q_1)Q_2'+Q_1\Delta(\chi Q_2)-(\Delta\chi)Q_1Q_2'= \\
&=(\Delta Q^{[0]})(Q_2^{[0]})' +
(\Delta Q^{[0]})(Q_2^{[\infty]})'+
Q_1^{[0]}\Delta Q_2^{[0]} +
Q_1^{[\infty]}\Delta Q_2^{[0]}
\\
&-(\Delta\chi)(Q_1^{[0]}+Q_1^{[\infty]})(Q_2^{[0]}+Q_2^{[\infty]})
\end{align*}
is estimated as
$$
|\Delta(Q_1Q_2)^{[0]}| \leq C \|Q_1;\mathfrak{C}_{(m)}^{0,\alpha}\|
\cdot\|Q_2;\mathfrak{C}_{(n)}^{0,\alpha}\|\cdot |\mathbf{p}'-\mathbf{p}|^{\alpha}.$$

For $R_0 \leq r'\wedge  r$, we have
\begin{align*}
&\Delta (Q_1Q_2)^{[\infty]}_{\star(n)}=\Delta\Bigg[\Big(\frac{r}{R_0}\Big)^{n-2}(1-\chi)Q_1Q_2 \Bigg]= \\
&=\Delta\Big(\frac{r}{R_0}\Big)^{n-2}\cdot (Q_1^{[\infty]})'Q_2' + \\
&+\Big(\frac{r}{R_0}\Big)^{n-2}\Bigg[\Delta Q_1^{[\infty]}\cdot Q_2'+
Q_1\cdot \Delta Q_2^{[\infty]}+\Delta\chi\cdot Q_1Q_2'\Bigg] \\
&=\Delta\Big(\frac{r}{R_0}\Big)^{n-2}\cdot\Big(\frac{R_0}{r}\Big)^{m-2}\cdot (Q_1^{[\infty]})_{\star(m)}\cdot Q_2 + \\
&+\Big(\frac{r}{R_0}\Big)^{n-2}\cdot\Delta\Big(\frac{R}{r}\Big)^{m-2}\cdot
(Q_1^{[\infty]})_{\star(m)}'\cdot Q_2' + \\
&+\Big(\frac{r}{R_0}\Big)^{n-m}\cdot\Delta (Q_1^{[\infty]})_{\star(m)}Q_2' +
\Big(\frac{r}{R_0}\Big)^{n-2}\Delta\Big(\frac{R_0}{r}\Big)^{n-2}\cdot Q_1\cdot (Q_2^{[\infty]})_{\star(n)} + \\
& +Q_1\cdot \Delta(Q_2^{[\infty]})_{\star(n)} +
\Big(\frac{r}{R_0}\Big)^{n-2}\Delta\chi\cdot Q_1Q_2'.
\end{align*}
Keeping in mind that $Q_1=Q_1^{[0]}+Q_1^{[\infty]}, Q_2=Q_2^{[0]}+Q_2^{[\infty]}$, we have the estimate
$$|\Delta((Q_1Q_2)^{[\infty]}_{\star(n)}|\leq C
\|Q_1;\mathfrak{C}_{(m)}^{0,\alpha}\|\cdot\|Q_2;\mathfrak{C}_{(n)}^{0,\alpha}\|\cdot
|(\mathbf{p}')^{\star}-\mathbf{p}^{\star}|^{\alpha}. $$
$\square$.\\

\subsection{Existence of the global metric}

We shall find a set of solutions $W,Y, X, V$ such that
\begin{align*}
&W \in \mathfrak{C}_{(3)}^{2,\alpha}, \quad
Y \in \mathfrak{C}_{(5)}^{2,\alpha}, \quad
X \in \mathfrak{C}_{(4)}^{2,\alpha}, \\
&V \upharpoonright \bar{\mathfrak{D}}(R_0) \in \mathfrak{C}^{0,\alpha}(\bar{\mathfrak{D}}(R_0)),\quad V^{\flat(3)} \in C^1(\mathbb{R}^3),\quad \lim_{r\rightarrow+\infty}V=0.
\end{align*}
If it will be done, then, thanks to Lemma \ref{removable} applied to $W, X, Y$ and the equations for $V$ (or $K$), there are constants $C_{\infty}^W, C_{\infty}^Y, C_{\infty}^X $ such that
\begin{align*}
& W=\frac{1}{r}\Big(C_{\infty}^w+O\Big(\frac{1}{r}\Big)\Big), \quad Y=\frac{1}{r^3}\Big(C_{\infty}^Y+O\Big(\frac{1}{r}\Big)\Big), \\
& X=\frac{1}{r^2}\Big(C_{\infty}^X+O\Big(\frac{1}{r}\Big)\Big),\\
&V= O\Big(\frac{1}{r^2}\Big)
\end{align*}
as $r \rightarrow +\infty$.
 This means that we have the expected flat asymptotic behavior \eqref{7.2a} $\sim$ \eqref{7.2d}, that is, 
\begin{align*}
&F=-\frac{\mathsf{G}M}{\mathsf{c}^2r}+O\Big(\frac{1}{r^2}\Big), \quad A=\frac{2\mathsf{G}J\varpi^2}{\mathsf{c}^3r^3}+O\Big(\frac{\varpi^2}{r^4}\Big), \\
&\Pi=\varpi\Big(1+O\Big(\frac{1}{r^2}\Big)\Big), \quad e^K=1+O\Big(\frac{1}{r^2}\Big)
\end{align*}
as $ r\rightarrow +\infty$. Here
$$M=M_{\mathsf{N}}+\frac{C_{\infty}^W}{\mathsf{G}\mathsf{c}^2},\qquad J=\frac{C_{\infty}^Y}{2\mathsf{G}}, $$
while $M_{\mathsf{N}}$ stands for the  total mass 
$\displaystyle \int\rho_{\mathsf{N}}^{\flat(3)}(\mbox{\boldmath$\xi$})d\mbox{\boldmath$\xi$}$ of the Newtonian limit
so that $\displaystyle \Phi_{\mathsf{N}}=-\frac{\mathsf{G}M_{\mathsf{N}}}{r}+O\Big(\frac{1}{r^2}\Big)$ as $ r\rightarrow +\infty$.

\begin{Remark}
We shall show later that $C_{\infty}^W=O(u_{\mathsf{O}}^2)$. Then 
$M \rightarrow M_{\mathsf{N}}$ as $u_{\mathsf{O}}/\mathsf{c}^2 \rightarrow 0$ with $u_{\mathsf{O}}=O(1)$. So $M>0$ for $u_{\mathsf{O}}/\mathsf{c}^2 \ll 1$, that is, the so called `positive mass' is realized. 
\end{Remark}

In order to construct such a set of solutions, we shall use the following

\begin{Definition}
 For any compactly supported continuous function $g$ on $\mathbb{R}^n$, we consider the function $\mathcal{K}^{(n)}g$ defined by
\begin{equation}
\mathcal{K}^{(n)}g(\mbox{\boldmath$\xi$})=
\frac{1}{(n-2)\Sigma_n}\int
\frac{g(\mbox{\boldmath$\xi$}')}{|\mbox{\boldmath$\xi$}-\mbox{\boldmath$\xi$}'|^{n-2}}d\mbox{\boldmath$\xi$}',
\end{equation}
where $\Sigma_n=2(n-2)\pi^{n/2}/\Gamma(n/2)$.
\end{Definition}

Under the situation of the above Definition suppose that $\mathrm{supp}[g] \subset \bar{B}^{(n)}(\Xi)$. 
Then
$$\|\mathcal{K}^{(n)}g\|_{C^1(\bar{B}^{(n)}(\Xi))}\leq C_0\|g\|_{C^0(\bar{B}^{(n)}(\Xi))},
$$
and if  $g \in C^{0,\alpha}(\bar{B}^{(n)}(\Xi))$, then 
$$\|\mathcal{K}^{(n)}g\|_{C^{2,\alpha}(\bar{B}^{(n)}(\Xi))} \leq C_0'\|g\|_{C^{0,\alpha}(\bar{B}^{(n)}(\Xi))}.$$
Here $C_0, C_0'$ depend only on $n, \alpha, \Xi$. We have
$$ \triangle^{(n)}.\mathcal{K}^{(n)}g=0
\quad\mbox{in}\quad \mathbb{R}^n \setminus \bar{B}^{(n)}(\Xi), 
$$
and 
$$\mathcal{K}^{(n)}g(\mbox{\boldmath$\xi$})=O\Big(\frac{1}{|\mbox{\boldmath$\xi$}|^{n-2}}\Big) \quad\mbox{as}\quad |\mbox{\boldmath$\xi$}| \rightarrow +\infty.$$

\begin{Definition}
 For $\mathfrak{g} \in \mathfrak{C}_{(n)}^0$, we define
$\mathfrak{K}^{(n)}\mathfrak{g}$ by
$$(\mathfrak{K}^{(n)}\mathfrak{g})^{\flat(n)}=f_0+f_{\infty},$$
where
\begin{align*}
&f_0=\mathcal{K}^{(n)}.\mathfrak{g}^{[0]\flat(n)}, \\
&f_{\infty}
=\Big(
\mathcal{K}^{(n)}[(\mathfrak{g}^{[\infty]\flat(n)})_{\blacklozenge(n)}] \Big)_{\star(n)}
\end{align*}
Here we define
\begin{equation}
g_{\blacklozenge(n)}(\mbox{\boldmath$\xi$}^{\star})=\Big(\frac{|\mbox{\boldmath$\xi$}^{\star}|}{\Xi_0}\Big)^4
g_{\star(n)}(\mbox{\boldmath$\xi$}^{\star})=
\Big(\frac{|\mbox{\boldmath$\xi$}|}{\Xi_0}\Big)^{n-6}g(\mbox{\boldmath$\xi$}).
\end{equation}
\end{Definition}

Under the situation of the above Definition we see that
$$\|g_{\blacklozenge(n)}; C^{l,\alpha}(B^{(n)}(\Xi_0))\|\leq C
\| g_{\star(n)}; C^{l,\alpha}(B^{(n)}(\Xi_0))\|,$$ for $(|\mbox{\boldmath$\xi$}^{\star}|/\Xi_0)^4$ is a quadratic polynomial of $\mbox{\boldmath$\xi$}^{\star}$. 
 We can consider that $(f_0)_{\star(n)}$ is harmonic in $ {B}^{(n)}(\Xi_0/2)$
and $f_{\infty}$ is harmonic in $B^{(n)}(\Xi_0)$ thanks to Lemma \ref{removable}, while the Poisson equations 
$$-\triangle f_0=\mathfrak{g}^{[0]\flat(n)},
\qquad -\triangle f_{\infty}=\mathfrak{g}^{[\infty]\flat(n)} $$
hold on the whole space $\mathbb{R}^n$.\\

Under these notations, we can claim

\begin{Proposition}
Let $n=3,4,5$.

1) The operator $\mathfrak{K}^{(n)}$ is continuous from
$\mathfrak{C}_{(n)}^{0}$ into 
$\mathfrak{C}_{(n)}^1$ and
from $\mathfrak{C}_{(n)}^{0,\alpha}$ into
$\mathfrak{C}_{(n)}^{2,\alpha}$, these operator norms
being independent of $\mathsf{a}$. 

2) If $\mathfrak{g} \in \mathfrak{C}_{(n)}^{0,\alpha}$, then
$Q=\mathsf{a}^2\mathfrak{K}^{(n)}\mathfrak{g}$ satisfies the equation
$$\Big[\frac{\partial^2}{\partial\varpi^2}+\frac{n-2}{\varpi}\frac{\partial}{\partial\varpi}
+\frac{\partial^2}{\partial z^2}\Big]Q+\mathfrak{g}=0 $$
on the whole space. 
\end{Proposition}

The proof can be done by dint of the Lemma \ref{removable}. 
Actually let $\mathfrak{g} \in \mathfrak{C}_{(n)}^0$ and consider $\mathfrak{K}^{(n)}\mathfrak{g}$ defined by
\begin{align*}
&\mathfrak{K}^{(n)}\mathfrak{g}=\mathfrak{f}_0+\mathfrak{f}_{\infty}, \\
&\mathfrak{f}_0^{\flat(n)}=f_0=\mathcal{K}^{(n)}g^{[0]\flat(n)},
\quad
\mathfrak{f}_{\infty}^{\flat(n)}=f_{\infty}=\Big(\mathcal{K}^{(n)}[(g^{[\infty]\flat(n)})_{\blacklozenge(n)}
]\Big)_{\star(n)}.
\end{align*}
Since $\mathfrak{g}^{[0]\flat(n)}\in C^0(\mathbb{R}^n)$ and $\mathrm{supp}[\mathfrak{g}^{[0]\flat(n)}]
\subset \{ |\mbox{\boldmath$\xi$}| < 2\Xi_0\} $, we have
$$\|f_0; C^1(\bar{B}^{(n)}(2\Xi_0))\|\leq
C\|\mathfrak{g}^{[0]\flat(n)};C^0(\bar{B}^{(n)}(2\Xi_0))\|\leq C'
\|\mathfrak{g}; \mathfrak{C}_{(n)}^0\|,$$
$f_0=O(|\mbox{\boldmath$\xi$}|^{-n+2})$ as $ |\mbox{\boldmath$\xi$}|\rightarrow +\infty$ and $(f_0)_{\star(n)}$ is harmonic on
$B^{(n)}(\Xi_0/2)$. Therefore we have
\begin{align*}
\|\mathfrak{f}_0^{[0]};\mathfrak{C}^1(\bar{\mathfrak{D}}(R_0))\|&=
\|f_0^{[0]};C^1(\bar{B}^{(n)}(\Xi_0))\| \\
&\leq C\|f_0;C^1(\bar{B}^{(n)}(\Xi_0))\| \leq C'\|\mathfrak{g};\mathfrak{C}_{(n)}^0\|, \\
\|\mathfrak{f}_0^{[\infty]};\mathfrak{C}_{(n)}^1
(\overline{\mathfrak{D}^{\mathsf{c}}}(R_0))\|&=
\|(f_0^{[\infty]})_{\star(n)};C^1(\bar{B}^{(n)}(\Xi_0))\| \\
&\leq C \|f_0;C^1(\bar{B}^{(n)}(2\Xi_0))\| \vee
\|(f_0)_{\star(n)};C^1(\bar{B}(\Xi_0/2))\| \\
&\leq C' \|f_0;C^1(\bar{B}^{(n)}(2\Xi_0))\| \vee
\|(f_0)_{\star(n)};C^1(\partial{B}(\Xi_0/2))\| \\
&\qquad\qquad \mbox{($\because$   by Lemma \ref{removable} \eqref{Pineq.1})}\\
&\leq C''\|f_0;C^1(\bar{B}^{(n)}(2\Xi_0))\| \vee
\|f_0; C^1(\partial{B}(2\Xi_0))\| \\
&\leq C'''\|\mathfrak{g};\mathfrak{C}_{(n)}^0\|.
\end{align*}
Thus we can claim 
$$ \|\mathfrak{f}_0;\mathfrak{C}_{(n)}^1\|\leq C \|\mathfrak{g};\mathfrak{C}_{(n)}^0\|. $$ 
Similarly we can estimate $\|\mathfrak{f}_{\infty};\mathfrak{C}_{(n)}^1\|$ by $\|\mathfrak{g};\mathfrak{C}_{(n)}^0\|$, therefore we can claim $$\|\mathfrak{K}^{(n)}\mathfrak{g}; \mathfrak{C}_{(n)}^1\|\leq C
\|\mathfrak{g}; \mathfrak{C}_{(n)}^0\|.$$
 Let $\mathfrak{g} \in \mathfrak{C}_{(n)}^{0,\alpha}$.
Since $\mathfrak{g}^{[0]\flat(n)}\in C^{0,\alpha}(\mathbb{R}^n)$ and $\mathrm{supp}[\mathfrak{g}^{[0]\flat(n)}]
\subset \{ |\mbox{\boldmath$\xi$}| < 2\Xi_0\} $, we have
$$\|f_0; C^{2,\alpha}(\bar{B}^{(n)}(2\Xi_0))\|\leq
C\|\mathfrak{g}^{[0]\flat(n)};C^{0,\alpha}(\bar{B}^{(n)}(2\Xi_0))\|\leq C'
\|\mathfrak{g}; \mathfrak{C}_{(n)}^{0,\alpha}\|,$$
$f_0=O(|\mbox{\boldmath$\xi$}|^{-n+2})$ as $ |\mbox{\boldmath$\xi$}|\rightarrow +\infty$ and $(f_0)_{\star(n)}$ is harmonic on
$B^{(n)}(\Xi_0/2)$. Therefore we have
\begin{align*}
\|\mathfrak{f}_0^{[0]};\mathfrak{C}^{2,\alpha}(\bar{\mathfrak{D}}(R_0))\|&=
\|f_0^{[0]};C^{2,\alpha}(\bar{B}^{(n)}(\Xi_0))\| \\
&\leq C\|f_0;C^{2,\alpha}(\bar{B}^{(n)}(\Xi_0))\| \leq C'\|\mathfrak{g};\mathfrak{C}_{(n)}^{0,\alpha}\|, \\
\|\mathfrak{f}_0^{[\infty]};\mathfrak{C}_{(n)}^{2,\alpha}
(\overline{\mathfrak{D}^{\mathsf{c}}}(R_0))\|&=
\|(f_0^{[\infty]})_{\star(n)};C^{2,\alpha}(\bar{B}^{(n)}(\Xi_0))\| \\
&\leq C \|f_0;C^{2,\alpha}(\bar{B}^{(n)}(2\Xi_0))\| \vee
\|(f_0)_{\star(n)};C^{2,\alpha}(\bar{B}(\Xi_0/2))\| \\
&\leq C' \|f_0;C^{2,\alpha}(\bar{B}^{(n)}(2\Xi_0))\| \vee
\|(f_0)_{\star(n)};C^{2,\alpha}(\partial{B}(\Xi_0/2))\| \\
&\qquad \qquad \mbox{($\because$   by Lemma \ref{removable} \eqref{Pineq.2})}\\
&\leq C''\|f_0;C^{2,\alpha}(\bar{B}^{(n)}(2\Xi_0))\| \vee
\|f_0; C^{2,\alpha}(\partial{B}(2\Xi_0))\| \\
&\leq C'''\|\mathfrak{g};\mathfrak{C}_{(n)}^{0,\alpha}\|.
\end{align*}
Thus we can claim 
$$ \|\mathfrak{f}_0;\mathfrak{C}_{(n)}^{2,\alpha}\|\leq C \|\mathfrak{g};\mathfrak{C}_{(n)}^{0,\alpha}\|. $$ 
Similarly we can estimate $\|\mathfrak{f}_{\infty};\mathfrak{C}_{(n)}^{2,\alpha}\|$ by $\|\mathfrak{g};\mathfrak{C}_{(n)}^{0,\alpha}\|$,  therefore we can claim $$\|\mathfrak{K}^{(n)}\mathfrak{g}; \mathfrak{C}_{(n)}^{2,\alpha}\|\leq C
\|\mathfrak{g}; \mathfrak{C}_{(n)}^{0,\alpha}\|.$$\\

Moreover, thanks to the same reasoning as that of \cite[Proposition 10]{asEE}, we can claim the following
\begin{Proposition}
There is a bounded linear operator $\mathfrak{L}$ which enjoys the following properties:

1) $\mathfrak{L}$ is continuous from $\mathfrak{C}_{(3)}^{0}$ into
$\mathfrak{C}_{(3)}^{1}$ and
from $\mathfrak{C}_{(3)}^{0,\alpha}$ into
$\mathfrak{C}_{(3)}^{2,\alpha}$, the operator norm being independent of $\mathsf{a}$;

2) For $\mathfrak{g}\in \mathfrak{C}_{(3)}^{0,\alpha}$, the function $Q=\mathsf{a}^2\mathfrak{Lg}$ satisfies the equation
$$\Big[
\frac{\partial^2}{\partial \varpi^2}+\frac{1}{\varpi}\frac{\partial}{\partial\varpi}+\frac{\partial^2}{\partial z^2}+4\pi\mathsf{G}\frac{1}{\gamma-1}\frac{\rho_{\mathsf{N}}}{u_{\mathsf{N}}}\Big] Q+\mathfrak{g}=0
$$ in the whole space and satisfies $Q(O)=0$.
\end{Proposition}

In fact the equation to be solved for $W=Q^{\flat(3)}, G=\mathsf{a}^2\mathfrak{g}^{\flat(3)}$ is
$$\triangle W +
\frac{1}{\gamma-1}(\Theta\vee 0)^{\frac{2-\gamma}{\gamma-1}}W
+G=0.$$
The solution $W$ is given as $W_0+W_{\infty}$ by the solutions  $W_0, W_{\infty}$ of the equations
\begin{align*}
&\triangle W_0+
\frac{1}{\gamma-1}(\Theta\vee 0)^{\frac{2-\gamma}{\gamma-1}}W_0+G^{[0]}=0, \\
&\triangle W_{\infty}+G^{[\infty]}=0.
\end{align*}
Here we recall $\Theta\vee 0 =0$ on the domain 
$B(2\xi_1(\frac{1}{\gamma-1}))^{\mathsf{c}} \supset  B(\Xi_0)^{\mathsf{c}}$. The former equation can be treated as \cite[Proposition 10]{asEE}, say, as the integral equation
$$W_0=\mathcal{K}\Big[
\frac{1}{\gamma-1}(\Theta\vee 0)^{\frac{2-\gamma}{\gamma-1}}W_0+G^{[0]}\Big]$$
with
$\mathcal{K} : f \mapsto \mathcal{K}^{3)}f-\mathcal{K}^{(3)}f(O),
$
and the solution $W_{\infty}$ is given by
$$W_{\infty}=\Big(\mathcal{K}^{(3)}[(G^{[\infty]})_{\blacklozenge(3)}]\Big)_{\star(3)}.
$$\\

We are going to find solutions $W,Y, X, V$ in
$$\mathfrak{C}_{(3)}^{2,\alpha},
\mathfrak{C}_{(5)}^{2,\alpha},
\mathfrak{C}_{(4)}^{2,\alpha},
\mathfrak{C}^{0,\alpha}(\bar{\mathfrak{D}}(R_0)).
$$\\

In order to fix the idea, hereafter we suppose the following\\

{\bf (D2)}: {\it It holds that
\begin{equation}
u_{\mathsf{O}}\leq C^0, \quad \frac{1}{\mathsf{c}^2}\leq C^0,\quad
\frac{1}{\mathsf{a}^2}\leq C^0, \quad
\frac{u_{\mathsf{O}}}{\mathsf{c}^2}\leq\delta^0
\end{equation}
with a fixed constant $C^0$ and a sufficiently small positive number $\delta^0$}.\\

Here let us recall Definition \ref{Def.R0} which gives
\begin{equation}
\mathsf{a}=\frac{1}{\sqrt{4\pi\mathsf{G}}}\Big(\frac{\mathsf{A}\gamma}{\gamma-1}\Big)^{\frac{1}{2(\gamma-1)}}
u_{\mathsf{O}}^{-\frac{2-\gamma}{2(\gamma-1)}}, \quad
\mathsf{b}=\frac{1}{4\pi\mathsf{G}}\Big(\frac{\mathsf{A}\gamma}{\gamma-1}\Big)^{\frac{1}{\gamma-1}}
\Omega_{\mathsf{O}}^2u_{\mathsf{O}}^{-\frac{1}{\gamma-1}}.
\end{equation}
Hence
\begin{equation}
\Omega_{\mathsf{O}}^2=\mathsf{a}^{-2}\mathsf{b}u_{\mathsf{O}}
\end{equation}
so that
$$\frac{\Omega_{\mathsf{O}}^2}{\mathsf{c}^2}\leq C^0\beta^0\cdot 
\frac{ u_{\mathsf{O}} }{\mathsf{c}^2} \leq C^0\beta^0\delta^0 $$
provided {\bf (D1)(D2) }.

Note that $u_{\mathsf{N}}\leq u_{\mathsf{O}}$ everywhere. We see
\begin{align}
&\|u_{\mathsf{N}}: \mathfrak{C}_{(3)}^{2,\alpha}\|\leq Cu_{\mathsf{O}}, \quad \|\Phi_{\mathsf{N}}:\mathfrak{C}_{(3)}^{2,\alpha}\| \leq Cu_{\mathsf{O}}, \nonumber \\
&\|\rho_{\mathsf{N}};\mathfrak{C}_{(3)}^{0,\alpha}\|\leq C\mathsf{G}^{-1}\mathsf{a}^{-2}u_{\mathsf{O}},
\quad \|P_{\mathsf{N}}:\mathfrak{C}_{(3)}^{1,\alpha}\|\leq C\mathsf{G}^{-1}\mathsf{a}^{-2}u_{\mathsf{O}}^2, \nonumber \\
& \Big\|\frac{\rho_{\mathsf{N}}}{u_{\mathsf{N}}}; \mathfrak{C}_{(3)}^{0,\alpha}\Big\| \leq C \mathsf{G}^{-1}\mathsf{a}^{-2}.
\end{align}\\

{\bf (D3)  Specification of $\delta_0$}:
{\it We take $\delta_0>0$ such that
\begin{equation}
\frac{|w|}{\mathsf{c}^2}\leq\delta_0u_{\mathsf{O}} \qquad \Rightarrow \quad
u_{\mathsf{N}}+\frac{w}{\mathsf{c}^2}<0\quad\mbox{for}\quad 3r_1\leq r.
\end{equation}
}\\

It is possible, since we have extended $\Theta$ so that 
$\sup\{\Theta | 2\xi_1\leq  |\mbox{\boldmath$\xi$}|<+\infty \}<0$. Then $\{u >0\}=\{\rho >0\}
\subset \mathfrak{D}(3r_1)$, provided that $|w|/\mathsf{c}^2 \leq\delta_0u_{\mathsf{O}}$. \\

{ \bf (D4)  Specification of $\delta_1$}: {\it
We take 
 $\delta_1>0$ such that,
if
\begin{equation}
\frac{1}{\mathsf{c}^2}|W|\leq\delta_1u_{\mathsf{O}},\quad
\frac{1}{\mathsf{c}^2}|Y|\leq\delta_1|\Omega_{\mathsf{O}}|,\quad
\frac{1}{\mathsf{c}^2}|X|\leq\delta_1u_{\mathsf{O}},
\end{equation}
 then it holds that
 $|w|/\mathsf{c}^2\leq\delta_0u_{\mathsf{O}}$,
while $|Z|/\mathsf{c}^2 \ll 1$. Here we take $\beta^0$ of the assumption {\bf (D1)}
and $\delta^0$ of {\bf (D2)}  smaller if necessary. }\\

It is possible, since  \eqref{Ww} with \eqref{WwZ} implies the existence of such $\delta_1$. Actually \eqref{Ww}, \eqref{WwZ} imply the following estimates on $\bar{\mathfrak{D}}(2R_0)$:
\begin{align*}
\frac{1}{\mathsf{c}^2}|w-W|&\leq C\Bigg[
\mathsf{a}^2|\Omega_{\mathsf{O}}|^2\frac{|\Phi_{\mathsf{N}}|}{\mathsf{c}^2}+\mathsf{a}^2|\Omega_{\mathsf{O}}|\frac{|Y|}{\mathsf{c}^2} +\mathsf{a}^4|\Omega_{\mathsf{O}}|^4\frac{1}{\mathsf{c}^2}+\frac{1}{\mathsf{c}^4}|\mathfrak{Q}_0|
\Bigg], \\
\frac{1}{\mathsf{c}^2}|\mathfrak{Q}_0|&\ll 
\mathsf{a}^2|\Omega_{\mathsf{O}}|\frac{|Y|}{\mathsf{c}^2}+\frac{|\Phi_{\mathsf{N}}|}{\mathsf{c}^2}+
\frac{1}{\mathsf{c}^4}(|W|+|X|),
\end{align*}
provided that
$$\frac{|\Phi_{\mathsf{N}}|}{\mathsf{c}^2}+\frac{|W|}{\mathsf{c}^4}+\frac{|X|}{\mathsf{c}^4}
\ll 1
\quad\mbox{and}\quad \mathsf{a}^2
|\Omega_{\mathsf{O}}|\frac{|Y|}{\mathsf{c}^2} \ll 1,
$$
and we have $w-W=0$ on $\bar{\mathfrak{D}}(2R_0)^{\mathsf{c}}$ where $\Omega=0$. Recall
$\mathsf{a}^2\Omega_{\mathsf{O}}^2/\mathsf{c}^2 \leq \beta^0u_{\mathsf{O}}/\mathsf{c}^2$.\\

Supposing that $V \in \mathfrak{C}^{0,\alpha}(\bar{\mathfrak{D}}(R_0))$ is given, we solve the equations
\eqref{PNEq.a}, \eqref{PNEq.b}, \eqref{PNEq.c} for unknowns $W, Y, X$
by solving the integral equations
\begin{subequations}
\begin{align}
&W=\mathsf{a}^2\mathfrak{L}(\mathfrak{g}_a-4\pi\mathsf{G}\frac{1}{\gamma-1}\frac{\rho_{\mathsf{N}}}{u_{\mathsf{N}}}
\Omega\varpi^2Y+\mathfrak{R}_a), \label{Inteq.a} \\
&Y=\mathsf{a}^2\mathfrak{K}^{(5)}(\mathfrak{g}_b+\mathfrak{R}_b), \label{Inteq.b} \\
&X=\mathsf{a}^2\mathfrak{K}^{(4)}(\mathfrak{g}_c+\mathfrak{R}_c), \label{Inteq.c}
\end{align}
\end{subequations}
where
\begin{subequations}
\begin{align}
\mathfrak{g}_a&=-4\pi\mathsf{G}\frac{1}{\gamma-1}\frac{\rho_{\mathsf{N}}}{u_{\mathsf{N}}}
\Omega^2\Big(2\Phi_{\mathsf{N}}-\frac{1}{4}\Omega^2\varpi^2\Big)+ \nonumber \\
&+4\pi\mathsf{G}\Upsilon_1\rho_{\mathsf{N}}u_{\mathsf{N}}
-8\pi\mathsf{G}\rho_{\mathsf{N}}(\Phi_{\mathsf{N}}+2\Omega^2\varpi^2)+
12\pi\mathsf{G}P_{\mathsf{N}}, \\
\mathfrak{g}_b&=16\pi\mathsf{G}\Omega\rho_{\mathsf{N}}, \\
\mathfrak{g}_c&=-16\pi\mathsf{G}P_{\mathsf{N}}.
\end{align}
\end{subequations}
Here we suppose that $\{\rho>0\} \subset \mathfrak{D}(R_0)$. Then, since the evaluation  of $\mathfrak{R}_a, \mathfrak{R}_b, \mathfrak{R}_c$ can be done by using the values of $V$ (or $K$) only on the domain $\mathfrak{D}(R_0)$, it is sufficient that $V$ on $\mathfrak{D}(R_0)$ is given. \\

We see
\begin{subequations}
\begin{align}
&\|\mathfrak{g}_a;\mathfrak{C}_{(3)}^{0,\alpha}\|\leq C_1 \mathsf{a}^{-2}u_{\mathsf{O}}^2, \label{ga}\\
&\|\mathfrak{g}_b;\mathfrak{C}_{(5)}^{0,\alpha}\|\leq C_1\mathsf{a}^{-2}|\Omega_{\mathsf{O}}|u_{\mathsf{O}}, \label{gb}\\
&\|\mathfrak{g}_c;\mathfrak{C}_{(5)}^{0,\alpha}\|\leq C_1 \mathsf{a}^{-2}u_{\mathsf{O}}^2.\label{gc}
\end{align}
\end{subequations}

Given $W \in \mathfrak{C}_{(3)}^{2,\alpha}\subset \mathfrak{C}_{(3)}^1,
Y \in \mathfrak{C}_{(5)}^{2,\alpha} \subset \mathfrak{C}_{(5}^1,
X \in \mathfrak{C}_{(4)}^{2,\alpha}\subset \mathfrak{C}_{(4)}^1$, we evaluate $\mathfrak{g}_a, \mathfrak{g}_b, \mathfrak{g}_c, \mathfrak{R}_a, \mathfrak{R}_b, \mathfrak{R}_c$ by them, and we put
\begin{subequations}
\begin{align}
&\tilde{W}=\mathsf{a}^2\mathfrak{L}(\mathfrak{g}_a
-4\pi\mathsf{G}\frac{1}{\gamma-1}\frac{\rho_{\mathsf{N}}}{u_{\mathsf{N}}}
\Omega\varpi^2\tilde{Y}+\mathfrak{R}_a), \label{FPeq.a} \\
&\tilde{Y}=\mathsf{a}^2\mathfrak{K}^{(5)}(\mathfrak{g}_b+\mathfrak{R}_b), \label{FPeq.b} \\
&\tilde{X}=\mathsf{a}^2\mathfrak{K}^{(4)}(\mathfrak{g}_c+\mathfrak{R}_c). \label{FPeq.c}
\end{align}
\end{subequations}
Here $\tilde{Y}$ in the right hand side of \eqref{FPeq.a} means $\tilde{Y}$ determined by
\eqref{FPeq.b}.

Note that 
$$|\Omega_{\mathsf{O}}|\|4\pi\mathsf{G}\frac{1}{\gamma-1}\frac{\rho_{\mathsf{N}}}{u_{\mathsf{N}}}\Omega\varpi^2 ;
\mathfrak{C}_{(3)}^{0,\alpha}\|
\leq C'\mathsf{b}\mathsf{a}^{-2}u_{\mathsf{O}}\leq C\mathsf{a}^{-2}u_{\mathsf{O}}, $$
for $\displaystyle \frac{\rho_{\mathsf{N}}}{u_{\mathsf{N}}}\Omega\varpi^2=0$ on $\bar{\mathfrak{D}}(2R_0)^{\mathsf{c}}$ where $\Omega=0$.\\

Then, after these preparations, by a similar argument as that of \cite{asEE}, we can claim the following assertion:\\

 Let us fix $V$ given on $\bar{\mathfrak{D}}(R_0)$ so that $\|V;
\mathfrak{C}^{0,\alpha}(\bar{\mathfrak{D}}(R_0))
\|\leq u_{\mathsf{O}}^2M$ with
$u_{\mathsf{O}}M/\mathsf{c}^2 \leq C_0$.

Let 
\begin{equation}
\|W;\mathfrak{C}_{(3)}^1\| \leq u_{\mathsf{O}}^2B, \quad 
\|Y;\mathfrak{C}_{(5)}^1\| \leq |\Omega_{\mathsf{O}}|u_{\mathsf{O}}B,
\quad \|X:\mathfrak{C}_{(4)}^1\| \leq u_{\mathsf{O}}^2B, \label{Bounce1}
\end{equation}
with
\begin{equation}
 u_{\mathsf{O}}B/\mathsf{c}^2 \leq \delta_1
\end{equation}
and
\begin{equation}
\|W;\mathfrak{C}_{(3)}^{2,\alpha}\|\leq u_{\mathsf{O}}^2\hat{B},\quad 
\|Y;\mathfrak{C}_{(5)}^{2,\alpha}\|\leq |\Omega_{\mathsf{O}}|u_{\mathsf{O}}\hat{B},
\quad \|X:\mathfrak{C}_{(4)}^{2,\alpha}\|\leq u_{\mathsf{O}}^2\hat{B} \label{Bounce2}
\end{equation}
with
\begin{equation}
B\leq \hat{B}.
\end{equation}
Then we can find $B, \hat{B}$ such that
$\tilde{W},\tilde{Y},\tilde{X}$ ensure the same estimates, provided that we take $\beta^0$ in 
the assumption {\bf (D1)} and $\delta^0$ 
in the assumption {\bf (D2)}
 smaller if necessary. 
Since \eqref{Bounce1} is supposed with $u_{\mathsf{O}}B/\mathsf{c}^2 \leq \delta_1$,
it is guaranteed that $\{\rho >0\} \subset \mathfrak{D}(R_0)$. 

However, when we repeat the argument of \cite{asEE}, which was done on the bounded domain, in our situation on the whole space, we must deal 
with estimates of 
$\|\mathfrak{R}_a;\mathfrak{C}_{(3)}^{0}\|, \|\mathfrak{R}_a;\mathfrak{C}_{(3)}^{0,\alpha}\|$ 
and $\|\mathfrak{R}_b; \mathfrak{C}_{(5)}^{0,\alpha}\|$, while $\mathfrak{g}_a$, $\mathfrak{g}_b$, $
\mathfrak{g}_c$, $\mathfrak{R}_c$ vanish on 
$\mathfrak{D}(3r_1)^{\mathsf{c}} \supset \mathfrak{D}(R_0)^{\mathsf{c}}=\mathfrak{D}(4r_1)^{\mathsf{c}}$ and can be neglected in the exterior domain. To do so,
we recall Propositions \ref{Prop.DExt}, \ref{Prop.Product}, \ref{Prop.leq} prepared in advance..

As for the estimates of $H_{\rho}(w)$ which are necessary to estimate $\mathfrak{R}_a$ the Proposition \cite[Proposition 11]{asEE} should read

\begin{Proposition}\label{Prop.11}
Suppose $w, w_1, w_2 \in \mathfrak{C}_{(3)}^0$ and 
$\displaystyle u_{\mathsf{N}}+\frac{w}{\mathsf{c}^2} <0$ for $r \geq 3r_1$. Then, since $R_0 = 4r_1$, we have
$H_{\rho}(w)=H_{\rho}(w_1)=H_{\rho}(w_2)=0$ on $\mathfrak{D}(R_0)^{\mathsf{c}}$ and we have:

i)  
\begin{equation}
\|H_{\rho}(w); \mathfrak{C}^{0}_{(3)}\|\leq C
\Big(\frac{\|w;\mathfrak{C}_{(3)}^0\|}{\mathsf{c}^2}\Big)^{\frac{1}{\gamma-1}},
\end{equation}

ii) 
\begin{equation}
\|H_{\rho}(w_2)-H_{\rho}(w_1); \mathfrak{C}_{(3)}^0\|
\leq C \Big(\frac{1}{\mathsf{c}^2}(\|w_1;\mathfrak{C}_{(3)}^0\|+
\|w_2;\mathfrak{C}_{(3)}^0\|)\Big)^{\frac{1}{\gamma-1}-1}
\cdot
\frac{\|w_2-w_1;\mathfrak{C}_{(3)}^0\|}{\mathsf{c}^2}, \label{HW1W2.1}
\end{equation}

iii)  If $w \in \mathfrak{C}^1_{(3)}$, then
$H_{\rho}(w)\in \mathfrak{C}^{0, \alpha}_{(3)}$ and
\begin{align}
\|H_{\rho}(w);\mathfrak{C}^{0,\alpha}_{(3)}\|&\leq C
u_{\mathsf{O}}^{\frac{1}{\gamma-1}-1}\Big[
\Big(1+\frac{\|w;\mathfrak{C}_{(3)}^1\|}{\mathsf{c}^2u_{\mathsf{O}}}
\Big)^{\frac{1}{\gamma-1}-1}
\frac{\|w;\mathfrak{C}_{(3)}^0\|}{\mathsf{c}^2} + \nonumber \\
&+\frac{\|w;\mathfrak{C}_{(3)}^0\|}{\mathsf{c}^2u_{\mathsf{O}}}
\frac{\|w;\mathfrak{C}_{(3)}^{0,\alpha}\|}{\mathsf{c}^2}\Big],
\end{align}

and

iv) If $w_1,w_2 \in \mathfrak{C}^1_{(3)}$, then
\begin{align}
&\|H_{\rho}(w_2)-H_{\rho}(w_1);\mathfrak{C}_{(3)}^{0,\alpha}\| \leq \nonumber \\
& \leq C
u_{\mathsf{O}}^{\frac{1}{\gamma-1}-1}\Bigg[\Big(1+\frac{1}{\mathsf{c}^2u_{\mathsf{O}}}
(\|w_1;\mathfrak{C}_{(3)}^1\|+\|w_2;\mathfrak{C}_{(3)}^1\|\Big)^{\frac{1}{\gamma-1}-1}
\frac{\|w_2-w_1;\mathfrak{C}_{(3)}^0\|}{\mathsf{c}^2} + \nonumber \\
&+\frac{1}{\mathsf{c}^2u_{\mathsf{O}}}(\|w_1;\mathfrak{C}_{(3)}^0\|+\|w_2;\mathfrak{C}_{(3)}^0\|)
\frac{\|w_2-w_1;\mathfrak{C}_{(3)}^{0,\alpha}\|}{\mathsf{c}^2}\Bigg]. \label{HW1W2.2}
\end{align}

\end{Proposition}

Using Propositions \ref{Prop.DExt}, \ref{Prop.Product}, \ref{Prop.leq},
we can claim that, thanks to Proposition \ref{Prop.11}, i), 
$$
\|\mathfrak{R}_a; \mathfrak{C}_{(3)}^0\| \leq C
\mathsf{a}^{-2}\Bigg[\frac{u_{\mathsf{O}}^3}{\mathsf{c}^2}(B(1+B+\hat{B}+M)+M)+ 
\Big(\frac{u_{\mathsf{O}}}{\mathsf{c}^2}\Big)^{\frac{2-\gamma}{\gamma-1}}u_{\mathsf{O}}^2B\Bigg],
$$
and, thanks to Proposition \ref{Prop.11}, iii), 
$$
\|\mathfrak{R}_a; \mathfrak{C}_{(3)}^{0,\alpha}\| \leq C
\mathsf{a}^{-2}\Bigg[\frac{u_{\mathsf{O}}^3}{\mathsf{c}^2}(\hat{B}(1+\hat{B}+M)+M)
+ u_{\mathsf{O}}^2B\Bigg].
$$
Here we note that, looking at \eqref{Ww}, \eqref{WwZ}, we see
$$
\|w-W;\mathfrak{C}_{(3)}^1\| \leq C \beta^0u_{\mathsf{O}}^2(1+B), 
\qquad \|w-W;\mathfrak{C}_{(3)}^{1,\alpha}\|\leq C \beta^0u_{\mathsf{O}}^2(1+\hat{B}).
$$
As for $\mathfrak{R}_b, \mathfrak{R}_c$, it is relatively easy to see 
\begin{align}
&\|\mathfrak{R}_b; \mathfrak{C}_{(5)}^{0,\alpha}\| \leq C 
\mathsf{a}^{-2}|\Omega_{\mathsf{O}}|\frac{u_{\mathsf{O}}^2}{\mathsf{c}^2}(\hat{B}+M), \\
&\|\mathfrak{R}_c;\mathfrak{C}_{(4)}^{0,\alpha}\|\leq C
\mathsf{a}^{-2} \frac{u_{\mathsf{O}}^3}{\mathsf{c}^2}(B+M).
\end{align}\\

Summing up, if $W, Y, X$ satisfy \eqref{Bounce1}, \eqref{Bounce2}, then $\tilde{W}, \tilde{Y}, \tilde{X}$ satisfy
\begin{align*}
\|\tilde{W}; \mathfrak{C}_{(3)}^1\|&\leq C_1u_{\mathsf{O}}^2\Big[1+\frac{u_{\mathsf{O}}}{\mathsf{c}^2}\Big(B(1+B+\hat{B}+M)+M\Big)+
\Big(\frac{u_{\mathsf{O}}}{\mathsf{c}^2}\Big)^{\frac{2-\gamma}{\gamma-1}}(1+B)\Big]+ \\
&+C_2|\Omega_{\mathsf{O}}|\mathsf{a}^{-2}\|\tilde{Y};\mathfrak{C}_{(5)}^0\|, \\
\|\tilde{W}; \mathfrak{C}_{(3)}^{2,\alpha}\|&\leq C_1u_{\mathsf{O}}^2\Big[1+\frac{u_{\mathsf{O}}}{\mathsf{c}^2}\Big(\hat{B}(1+\hat{B}+M)+M\Big)+(1+B)
\Big]+ \\
&+C_2|\Omega_{\mathsf{O}}|\mathsf{a}^{-2}\|\tilde{Y};\mathfrak{C}_{(5)}^{0,\alpha}\|, \\
\|\tilde{Y};\mathfrak{C}_{(5)}^{2,\alpha}\|&\leq C_1|\Omega_{\mathsf{O}}|u_{\mathsf{O}}\Big[1+\frac{u_{\mathsf{O}}}{\mathsf{c}^2}(\hat{B}+M)\Big], \\
\|\tilde{X};\mathfrak{C}_{(4)}^{2,\alpha}\|&\leq C_1u_{\mathsf{O}}^2\Big[1+\frac{u_{\mathsf{O}}}{\mathsf{c}^2}(B+M)\Big].
\end{align*}
Here $C_1$ is replaced by a larger one than that in \eqref{ga} $\sim$ \eqref{gc} if necessary. 
Hence we can take $B, \hat{B}$ so that 
\begin{align*}
&C_1\Big[1+\delta^0\Big(B(1+B+\hat{B}+M)+M\Big)+(\delta^0)^{\frac{2-\gamma}{\gamma-1}}(1+B)\Big]+C_2\beta^0\hat{B} \leq B, \\
&C_1\Big[1+\delta^0\Big(\hat{B}(1+\hat{B}+M)+M\Big)+(1+B)\Big]+C_2\beta^0\hat{B} \leq \hat{B}, \\
&C_1\Big[1+\delta^0(\hat{B}+M)\Big] \leq \hat{B}, \\
&C_1\Big[1+\delta^0(B+M)\Big]\leq \hat{B},\\
&\delta^0B \leq \delta_1, \qquad B \leq \hat{B}.
\end{align*}
Then $\tilde{W}, \tilde{X}, \tilde{Y}$ enjoy the estimates \eqref{Bounce1}, \eqref{Bounce2}. Here we take $\delta^0, \beta^0$ smaller if necessary.
More concretely speaking, it is sufficient to take 
\begin{align}
&B:=4C_1+C_2,\nonumber \\
&\hat{B}:=C_1(4+4C_1+C_2)+C_2\quad  (=(1+C_1)B),
\end{align}
provided that
\begin{align}
&\delta^0B(1+B+\hat{B})\leq 1, \quad
\delta^0(1+B)M \leq 1,
\quad (\delta^0)^{\frac{2-\gamma}{\gamma-1}}(1+B)\leq 1, \nonumber \\
& \delta^0\hat{B}(1+\hat{B})\leq 1, \quad \delta^0(1+\hat{B})M \leq 1, \nonumber \\
&\beta^0\hat{B} \leq 1.
\end{align}
Note that $B, \hat{B}$ and the smallness of $\beta^0$ can be specified  independently of $M$, although the smallness of $\delta^0$ may depend on $M$.

Thus by the same arguments as that of \cite{asEE}, we can claim that the mapping $(W, Y, X) \mapsto (\tilde{W},\tilde{Y}, \tilde{X})$ turns out to be a contraction mapping from the functional set 
\begin{equation}
\mathfrak{X}:=\{(W,X,Y) \  |\   \mbox{\eqref{Bounce1} and\eqref{Bounce2}  hold} \}
\end{equation}
into itself with respect to a suitable distance. Actually the distance
$\mathrm{dist.}(U_2,U_1):=\mathfrak{N}(U_2-U_1)$ works, where
\begin{align}
&\mathcal{N}(U):=\|W;\mathfrak{C}_{(3)}^1\|\vee\Big(|\Omega_{\mathsf{O}}|^{-1}u_{\mathsf{O}}\|Y;\mathfrak{C}_{(5)}^1\|\Big)
\vee\|X;\mathfrak{C}_{(4)}^1\|, \\
&\hat{\mathcal{N}}(U):=\|W;\mathfrak{C}_{(3)}^{2,\alpha}\|\vee\Big(|\Omega_{\mathsf{O}}|^{-1}u_{\mathsf{O}}\|Y;\mathfrak{C}_{(5)}^{2,\alpha}\|\Big)
\vee\|X;\mathfrak{C}_{(4)}^{2,\alpha}\|, \\
&\mathfrak{N}(U):=\mathcal{N}(U)+\kappa\hat{\mathcal{N}}(U),\qquad \kappa:=2(\delta^0)^{1-\alpha},
\end{align}
provided that $\Omega_{\mathsf{O}}\not=0$. (Of course, when $\Omega_{\mathsf{O}}=0$, then we neglect $Y\equiv 0$.)
We take $\delta^0$ sufficiently small. Thanks to   Proposition \ref{Prop.11}, ii), \eqref{HW1W2.1}, we have
\begin{equation}
\mathcal{N}(\tilde{U}_2-\tilde{U}_1)\leq C_N\Big[
(\delta^0)^{\alpha}\mathcal{N}(U_2-U_1)+\delta^0\hat{\mathcal{N}}(U_2-U_1)\Big],
\end{equation}
and, thanks to   Proposition \ref{Prop.11}, iv), \eqref{HW1W2.2}, we have
\begin{equation}
\hat{\mathcal{N}}(\tilde{U}_2-\tilde{U}_1)\leq C_N\Big[
\mathcal{N}(U_2-U_1)+\delta^0\hat{\mathcal{N}}(U_2-U_1)\Big],
\end{equation}
therefore
\begin{equation}
\mathfrak{N}(U_2-U_1)\leq K_N\Big[
\mathcal{N}(U_2-U_1)+\kappa'\hat{\mathcal{N}}(U_2-U_1)\Big].
\end{equation}
with $K_N:=C_N((\delta^0)^{\alpha}+\kappa), \kappa'=\delta^0(1+\kappa)((\delta^0)^{\alpha}+\kappa)^{-1}$. Assuming $\kappa =2(\delta^0)^{1-\alpha}\leq 1$, we have $\kappa'\leq \kappa$ so that
$$\mathfrak{N}(\tilde{U}_2-\tilde{U}_1)\leq K_N\mathfrak{N}(U_2-U_1).$$
Supposing $K_N=C_N((\delta^0)^{\alpha}+2(\delta^0)^{1-\alpha})<1$, we have a contraction with respect to $\mathfrak{N}$.

In the previous study \cite{asEE} we considered equations on a bounded domain on which $\Omega$ is keeping to be a constant. Now in this study we consider equations on the whole space, but we are setting  that $\Omega=\Omega_{\mathsf{O}}$, a constant, on the domain $\bar{\mathfrak{D}}(R_0)$, which includes the support of $\rho$ as result, and $\Omega=0$ on the exterior domain $\bar{\mathfrak{D}}(2R_0)^{\mathsf{c}}$; In the transit band $\bar{\mathfrak{D}}(2R_0)\setminus \bar{\mathfrak{D}}(R_0)$, we have
$\|\Omega; \mathfrak{C}^{2,\alpha}\|\leq C|\Omega_{\mathsf{O}}|$,
thanks to the fixed cut off function $\chi$.\\

Thus we can claim

\begin{Theorem} \label{Th.WYX}
There is a unique set of solutions $(W, Y, X)$ of \eqref{PNEq.a}\eqref{PNEq.b}\eqref{PNEq.c} in $\mathfrak{X}$ for any given $V \in \mathfrak{C}^{0,\alpha}(\bar{\mathfrak{D}}(R_0))$
with $\|V;\mathfrak{C}^{0,\alpha}(\bar{\mathfrak{D}}(R_0))\|\leq u_{\mathsf{O}}^2M$.
\end{Theorem}

Let us denote the solution of Theorem \ref{Th.WYX} by $U=(W,Y,X) = \mathcal{S}(V)$.\\

We are going to solve \eqref{Veq.a}, \eqref{Veq.b}, that is, we are looking for $V$ such that $V$ and $U=\mathcal{S}(V)$ satisfy \eqref{Veq.a},\eqref{Veq.b}. Note that $V$ does not appear explicitly in the right-hand sides of \eqref{Veq.a}, \eqref{Veq.b}.
We have
\begin{align}
&\hat{\mathcal{N}}(\mathcal{S}(V))\leq
C_3u_{\mathsf{O}}^2, \label{NV} \\
&\hat{\mathcal{N}}(\mathcal{S}(V_2)-\mathcal{S}(V_1))\leq
C_3\delta^0\|V_2-V_1;\mathfrak{C}^{0,\alpha}(\bar{\mathfrak{D}}(R_0))\|, \label{NV1V2}
\end{align}
while $\hat{\mathcal{N}}$ is equivalent to $\mathfrak{N}$, since
$\kappa\hat{\mathcal{N}}\leq\mathfrak{N}\leq (1+\kappa)\hat{\mathcal{N}}$. We can specify $C_3$ independently of $M$, while $\delta^0$ is supposed to be small depending on $M$.

Let us denote by $\mathcal{R}_1(V), \mathcal{R}_3(V)$ the right-hand sides of
the equations \eqref{Veq.a}, \eqref{Veq.b}, respectively, evaluated by $U=\mathcal{S}(V)$, and define
$\hat{V} $  by
\begin{equation}
\hat{V}(\varpi, z)=
\int_0^z\mathcal{R}_3(V)(0,z')dz'+
\int_0^{\varpi}\mathcal{R}_1(V)(\varpi', z)d\varpi'. \label{InteqV}
\end{equation}
The function $\hat{V}$ is considered to be defined on $\mathfrak{D}(+\infty):=
\{ (\varpi, z)| 0\leq \varpi, |z| <+\infty \}$, since $\mathcal{R}_1(V), \mathcal{R}_3(V) $ are defined on $\mathfrak{D}(+\infty)$.

Estimating the right-hand sides of the equations \eqref{N.EQd}, \eqref{N.EQe}, we see 
$$\Big|\frac{\partial \hat{V}}{\partial\varpi}\Big| \leq C\mathsf{a}^2u_{\mathsf{O}}^2\frac{1}{r^3}, \quad\Big|\frac{\partial \hat{V}}{\partial z}\Big| \leq C\mathsf{a}^2u_{\mathsf{O}}^2\frac{1}{r^3}, \quad\mbox{on}\quad r \geq R_0. $$ Integrating this, we see that 
$\hat{V}=C_{\infty}(V)+O(1/r^2)$ as $r \rightarrow +\infty$ with a constant $C_{\infty}(V)=O(u_{\mathsf{O}}^2)$. We put $\tilde{V}$ by
\begin{align}
\tilde{V}(\varpi,z)&=-C_{\infty}(V)+\hat{V}(\varpi,z) \nonumber \\
&=-C_{\infty}(V)+\int_0^z\mathcal{R}_3(V)(0,z')dz'+
\int_0^{\varpi}\mathcal{R}_1(V)(\varpi', z)d\varpi'. 
\end{align}The function $\tilde{V}$ is considered to be defined on $\mathfrak{D}(+\infty)$. Let us consider the operator
$$\mathcal{T}: \mathfrak{C}^{0,\alpha}(\bar{\mathfrak{D}}(R_0)) \rightarrow 
\mathfrak{C}^{0,\alpha}(\bar{\mathfrak{D}}(R_0)) : V \mapsto \tilde{V} \upharpoonright
\bar{\mathfrak{D}}(R_0).
$$
Thanks to Theorem \ref{ThN2}, if $V$ is a fixed point of $\mathcal{T}$, it is a solution.
Fixing  $M$ sufficiently large, we can take $\delta^0$ sufficiently small  such that $\mathcal{T}$ turns out to be a contraction mapping from
\begin{equation}
\mathfrak{V}=\{ V\  |\  \|V;
\mathfrak{C}^{0,\alpha}(\bar{\mathfrak{D}}(R_0))
\| \leq u_{\mathsf{O}}^2M\}
\end{equation}
into itself with respect to the distance given by $\|\cdot; 
\mathfrak{C}^{0,\alpha}(\bar{\mathfrak{D}}(R_0))
\|$ thanks to \eqref{NV}, \eqref{NV1V2}. Actually, looking at \eqref{Veq.a},\eqref{Veq.b}, we have,
 thanks to \eqref{NV}, \eqref{NV1V2},
\begin{equation}
\|\tilde{V};\mathfrak{C}^{0,\alpha}(\bar{\mathfrak{D}}(R_0))\|
\leq C_4\hat{\mathcal{N}}(\mathcal{S}(V)) 
\leq C_4 C_3  u_{\mathsf{O}}^2
\end{equation}
and
\begin{align}
\|\tilde{V}_2-\tilde{V}_1;\mathfrak{C}^{0,\alpha}(\bar{\mathfrak{D}}(R_0))\|
&\leq C_4 \hat{\mathcal{N}}(\mathcal{S}(V_2)-\mathcal{S}(V_1)) \nonumber \\
&\leq C_4C_3\delta^0\|V_2-V_1; \mathfrak{C}^{0,\alpha}(\bar{\mathfrak{D}}(R_0))\|.
\end{align}
In fact these estimates come from
\begin{align*}
&\|\tilde{V};\mathfrak{C}^{0,\alpha}(\bar{\mathfrak{D}}(R_0))\|
\leq C_4\hat{\mathcal{N}}_0(\mathcal{S}(V)) \\
&\|\tilde{V}_2-\tilde{V}_1;\mathfrak{C}^{0,\alpha}(\bar{\mathfrak{D}}(R_0))\|
\leq C_4 \hat{\mathcal{N}}_0(\mathcal{S}(V_2)-\mathcal{S}(V_1)) ,
\end{align*}
where 
\begin{align*}
\hat{\mathcal{N}}_0(U)&:= \|W\restriction \bar{\mathfrak{D}}(R_0));\mathfrak{C}^{2,\alpha}(\bar{\mathfrak{D}}(R_0))\|\vee \\
&\vee \Big(|\Omega_{\mathsf{O}}|^{-1}u_{\mathsf{O}}\|Y\restriction\bar{\mathfrak{D}}(R_0);\mathfrak{C}^{2,\alpha}(\bar{\mathfrak{D}}(R_0))\|\Big)\vee
\|X\restriction\bar{\mathfrak{D}}(R_0);\mathfrak{C}^{2,\alpha}(\bar{\mathfrak{D}}(R_0))\| \\
&\leq \hat{\mathcal{N}}(U).
\end{align*}
That is, we need not take care of values of $\mathcal{S}(V), \mathcal{S}(V_1), \mathcal{S}(V_2)$ on 
$\bar{\mathfrak{D}}(R_0)^{\mathsf{c}}=\{ (\varpi, z) | r >R_0\}$.
Note that $ C_4$ can be specified independently of $M$. 
We fix $M$ such that $M \geq C_4 C_3$ and replace $\delta^0$, which is supposed to be small depending on $M$,  by a smaller one so that
$C_4C_3\delta^0 <1$ if necessary.

  Therefore we have

\begin{Proposition}
There is a solution $V \in 
\mathfrak{C}^{0,\alpha}(\bar{\mathfrak{D}}(R_0))
$ of $V=\mathcal{T}(V)$ with
$U=\mathcal{S}(V)$. This is the unique solution in $\mathfrak{V}$ and the equations \eqref{Veq.a}, \eqref{Veq.b} are satisfied.
\end{Proposition}
It is easy to extend $V=\tilde{V}\restriction \bar{\mathfrak{D}}(R_0)$ to the whole space. Actually it can be done as nothing but $\tilde{V}$ by \eqref{InteqV}, since $\mathcal{R}_1(V), \mathcal{R}_3(V)$ are determined by the values of $V$ only on
the domain $\{ \rho >0\}\subset \mathfrak{D}(3r_1) \subset \mathfrak{D}(R_0)$. Thanks to Theorem \ref{ThN2} applied for arbitrarily large $R$, thus extended $V$ gives $K=V/\mathsf{c}^4$ which satisfies the equations \eqref{N.EQd}, \eqref{N.EQe} on the whole space.
We take $M$ larger if necessary. Thus we have
\begin{Theorem}\label{Th.K}
There is a solution $V$ of \eqref{Veq.a}\eqref{Veq.b} such that $V(O)=0$, $V^{\flat(3)}\in C^1(\mathbb{R}^3),
V\upharpoonright \bar{\mathfrak{D}}(R_0) \in
\mathfrak{C}^{0,\alpha}(\bar{\mathfrak{D}}(R_0))$,
$\|V \upharpoonright \bar{\mathfrak{D}}(R_0); \mathfrak{C}^{0,\alpha}(\bar{\mathfrak{D}}(R_0))\|\leq u_{\mathsf{O}}^2M$, and
$$V=O\Big(\frac{1}{r^2}\Big)\quad\mbox{as}\quad r\rightarrow \infty.
$$ The solution is unique in $\mathfrak{V}$. 
\end{Theorem}

Summing up, we can claim the following final :

\begin{Theorem}\label{Th.5}
Suppose {\bf (D0)}, {\bf (D1)}, and {\bf (D2)}. There is a set of potentials $F,A,\Pi, K$ and a function $u$, which enjoys $u(O)=u_{\mathsf{O}}$ and gives the density distribution $\rho=f^{\rho}(u)$ whose support is included in $\{ r < 3r_1=3\mathsf{a}\xi_1(\frac{1}{\gamma-1}) \}$, such that $F, A, \Pi/\varpi, K$ are $C^2$-functions on the whole space $\mathbb{R}^3$ and give the axisymmetric metric which satisfies the Einstein equations with the energy-momentum tensor of $\rho, \Omega_{\mathsf{O}} $ and is asymptotically flat. Moreover $F=\frac{1}{\mathsf{c}^2}\Phi_{\mathsf{N}}+O(u_{\mathsf{O}}^2/\mathsf{c}^4)$, $A=O(|\Omega_{\mathsf{O}}|u_{\mathsf{O}})\varpi^2/\mathsf{c}^3$, $\Pi=\varpi(1+O(u_{\mathsf{O}}^2/\mathsf{c}^4))$,
$K=O(u_{\mathsf{O}}^2/\mathsf{c}^4)$,
$u=u_{\mathsf{N}}+O(u_{\mathsf{O}}^2/\mathsf{c}^2) $. Here $u_{\mathsf{O}}=u_{\mathsf{N}}(O)$
and $\Phi_{\mathsf{N}}$ is the Newtonian gravitational potential generated by the density distribution $\rho_{\mathsf{N}}=f^{\rho}_{\mathsf{N}}(u_{\mathsf{N}})$.
\end{Theorem}

\section{Concluding Remark}

We have constructed an asymptotically flat metric with a rotating compactly supported 
perfect fluid. However, this work is done under the assumption that the angular velocity is constant and small. Rapidly or differentially rotating models are out of the scope of this work. Moreover the case with strong gravitational effect, say, the case with too big central density or $u_{\mathsf{O}}/\mathsf{c}^2 \gg 1$,  is out of the scope. The treatises of these cases are important and open problems. 

Actually as for the rapidly rotating models, we have much work to do before discussing such relativistic, post-Newtonian analysis as done in this article, since even in the non-relativistic Newtonian frame work the mathematical analysis of rapidly rotating gaseous star models has not yet been well developed. Imagine that the angular velocity becomes larger and larger. The stationary figure of the star might deviate from nearly ellipsoidal one to some very irregular configuration and finally could split into many non-connected pieces, say, a pair of binary stars or something like that could appear. The analysis of such bifurcations of rotating star models along the increase of the angular velocity even in the study of incompressible (homogeneous, say, not gaseous but liquid ) celestial bodies has a long history of innumerable many researches, by Maclaurin, Jacobi, Poincar\'{e}, Darwin, Tchebycheff, Liapounoff, Lichtenstein and so on. (See the monographs
\cite{Jardetzky}, \cite{ChandraEFE}, \cite{Hagihara}.) The serious disputes among these prominent mathematicians seem not to be completely solved, and there remain open problems for the study of compressible gaseous star models a fortiori, even in the non-relativistic problem setting. Direct treatise of the relativistic problem seems impossible without the development of the study of  rapidly rotating star models in the non-relativistic framework.

As for the differentially rotating models, in order to discuss them, first of all we should reconsider Proposition \ref{Prop4}. The proof of Proposition \ref{Prop4} presented in this article needs the assumption that $\Omega$ is constant where $\rho >0$, that is, that the rotation is rigid. But Proposition \ref{Prop4} is inevitable for Theorem \ref{Th.1} which guarantees the equivalence of the system of equations \eqref{N.EQa} $\sim$ \eqref{N.EQf} for $F, A, \Pi,  K, u$ and the full system of the Einstein equations. Completely different elegant proof than the direct but messy proof presented in this article could be provided by experts of the theory of differential geometry, when it might be allowed that $\Omega$ is not constant where $\rho >0$. . 

As for the treatise of the case with very strong gravitational effect by the relatively large central density, the author has no idea. Actually the discussion presented in this article is too heavily relying upon the post-Newtonian approximation. This situation is different from that of spherically symmetric problem, for which we can construct the interior metric by solving the Tolman-Oppenheimer-Volkoff equation without worrying about how large the central density be and extend it to the exterior vacuum metric by the Schwarzschild metric.

Moreover the so called `matter-vacuum matching problem' introduced in Section 1 cannot be said to be solved by this work. Namely, if the angular velocity is zero and spherically  symmetric metric is concerned, the metric in the exterior domain with the vacuum can be identified with the exterior part of the Schwarzschild metric by a suitable change of co-ordinates. However,
if $\Omega_{\mathsf{O}}\not=0$, it is not yet clarified 
whether there is a co-ordinate transformation which reduces the metric constructed in our wok to the exterior part of the Kerr metric in an exterior domain $\{ \rho =0\}$, or not. At least we can hardly expect that the vacuum boundary
$\partial\{ \rho >0\}$ of the density distribution constructed in this wok would be  an exact ellipsoid. In this sense, we cannot claim that we have solved the so called `matter-vacuum matching problem'. The problem is still open, although someone might be optimistic.

\vspace{20mm}
{\bf\Large Acknowledgment}\\

This work was partially done on the occasion of the BIRS-CMO Workshop `Time-like Boundaries in General Relativistic Evolution Problems (19w5140)' held at Oaxaca, Mexico on July 28 - August 2, 2019. The author would like to express his thanks to Professors Helmut Friedrich, Olivier Sarbach, Oscar Reula for organizing this workshop, to the Casa Matematica Oaxaca of the Banff International Research Station for Mathematical Innovation and Discovery for the hospitality, and to the participants for discussions. Especially the author deeply thanks to Professor O. Sarbach who constantly encouraged this work and gave helpful suggestions to ameliorate the fashion of presentation. 
 The author expresses his thanks to the anonymous referees who read the manuscript carefully and gave suggestions for amelioration of the text. Particularly, the acquaintance of important preceding achievements owes to their kind instructions.
This work is supported by JSPS KAKENHI Grant Number JP18K03371.

\vspace{20mm}

{\bf\Large Appendix:  Glossary of constants, parameters and variables}\\

Notation: $$Q\wedge Q':=\min(Q, Q'), \quad Q\vee Q':=\max(Q, Q'),\qquad\partial_j:=\frac{\partial}{\partial x^j}$$

{\bf 1. Positive constants and parameters}\\

$c$\  \  : the speed of light, \  \   $\mathsf{G}$\  \  : the gravitation constant, \  \  $\mathsf{A}$\  \  :the constant of isentropy, \  \  $\gamma$\  \  : the adiabatic exponent, \  \  $\alpha$\  \  : the Hl\"{o}lder continuity exponent:
$$ 1 <\gamma < 2,\quad 0<\alpha <\Big(\frac{1}{\gamma-1}-1\Big)\wedge 1. $$

$\rho_{\mathsf{NO}}$\  \  : the central value of $\rho_{\mathsf{N}}$, \  \  $u_{\mathsf{O}}$\  \  : the central value of $u_{\mathsf{N}}$ and $u$

$$\mathsf{a}=\sqrt{\frac{\mathsf{A}\gamma}{4\pi\mathsf{G}(\gamma-1)}}
\rho_{\mathsf{N}\mathsf{O}}^{-\frac{2-\gamma}{2}}
=
\frac{1}{\sqrt{4\pi\mathsf{G}}}\Big(\frac{\mathsf{A}\gamma}{\gamma-1}\Big)^{\frac{1}{2(\gamma-1)}}
u_{\mathsf{O}}^{-\frac{2-\gamma}{2(\gamma-1)}}
$$

$\xi_1\Big(\frac{1}{\gamma-1}\Big)$\  \  : the zero of the Lane-Emden function $\theta\Big(\xi,; \frac{1}{\gamma-1}\Big)$

$$r_1=\mathsf{a}\xi_1\Big(\frac{1}{\gamma-1}\Big),\quad \Xi_0=\xi_1\Big(\frac{1}{\gamma-1}\Big),
\quad R_0=4r_1=\mathsf{a}\Xi_0$$

{\bf 2. Nonnegative constants and parameters}\\

$\Omega_{\mathsf{O}}$\  \  : the constant angular velocity on a neighborhood of the support of $\rho$

$$\mathsf{b}=
\frac{\Omega_{\mathsf{O}}^2}{4\pi\mathsf{G}\rho_{\mathsf{N}\mathsf{O}}}
=\frac{1}{4\pi\mathsf{G}}\Big(\frac{\mathsf{A}\gamma}{\gamma-1}\Big)^{\frac{1}{\gamma-1}}
\Omega_{\mathsf{O}}^2u_{\mathsf{O}}^{-\frac{1}{\gamma-1}}
$$

{\bf 3. Variables and functions }\\

Greek indexes $\mu, \nu, \cdots $ stand for $0,1,2,3$

Latin indexes $k,j,\cdots$ stand for $1,2,3$

$x^0=\mathsf{c}t, \quad x^1=\varpi, \quad x^2=\phi, \quad x^3=z$\  \  \  : co-ordinates

$g_{\mu\nu}$\  \   : generally described coefficients of the metric
$$ds^2=g_{\mu\nu}dx^{\mu}dx^{\nu} $$
$R_{\mu\nu}$\  \  : the Ricci tensor

$$R=g^{\mu\nu}R_{\mu\nu}$$

$T^{\mu\nu}$\  \  : the energy-momentum tensor

$$T:=g^{\mu\nu}T_{\mu\nu}, \quad S_{\mu\nu}:=T_{\mu\nu}-\frac{1}{2}g_{\mu\nu}T $$
$F, A, K, \Pi$\  \  : coefficients of the axially symmetric metric (Lanczos form)
$$
ds^2=e^{2F}(\mathsf{c}dt+Ad\phi)^2
-e^{-2F}[e^{2K}(d\varpi^2+dz^2)+\Pi^2d\phi^2] $$

$U^{\mu}$\  \  : components of the 4-velocity of the fluid

$$U^0=e^{-G}, \quad U^1=U^3=0,\quad U^2=e^{-G}\frac{\Omega}{\mathsf{c}},
\quad \Omega :\mbox{angular velocity} $$

$\rho$\ \   : mass density, \  \  $P$: pressure, \  \  $u$: the enthalpy density:

$$u=\int_0^{\rho}\frac{dP}{\rho+P/\mathsf{c}^2}\quad\mbox{for}\quad \rho >0$$
but $u$ can take negative values and
$$\rho=f^{\rho}(u)=
\Big(\frac{\gamma-1}{\mathsf{A}\gamma}\Big)^{\frac{1}{\gamma-1}}
(u\vee 0)^{\frac{1}{\gamma-1}}\Big[1+\Upsilon_{\rho}(u/\mathsf{c}^2)\Big]$$
$$\Upsilon_{\rho}(X)=\sum_{k\geq 1}\Upsilon_kX^k : \quad\mbox{power series expansion}$$

$\rho_{\mathsf{N}}, P_{\mathsf{N}}, u_{\mathsf{N}}$\  \  : the Newtonian limits of the state variables $\rho, P, u$

$$\rho_{\mathsf{N}}=f_{\mathsf{N}}^{\rho}(u_{\mathsf{N}})=
\Big(\frac{\gamma-1}{\mathsf{A}\gamma}\Big)^{\frac{1}{\gamma-1}}(u_{\mathsf{N}}\vee 0)^{\frac{1}{\gamma-1}}$$

$\Phi_{\mathsf{N}}$\  \  : the Newtonian gravitational potential generated by $\rho_{\mathsf{N}}$

\begin{align*}
H_{\rho}(w)&:=f_{\mathsf{N}}^{\rho}(u_{\mathsf{N}}+\frac{w}{\mathsf{c}^2})
-f_{\mathsf{N}}^{\rho}(u_{\mathsf{N}})
-Df_{\mathsf{N}}^{\rho}(u_{\mathsf{N}})\frac{w}{\mathsf{c}^2} \nonumber \\
&=f_{\mathsf{N}}^{\rho}(u_{\mathsf{N}}+\frac{w}{\mathsf{c}^2})
-f_{\mathsf{N}}^{\rho}(u_{\mathsf{N}})
-
\frac{1}{\gamma-1}\frac{\rho_{\mathsf{N}}}{u_{\mathsf{N}}}
\frac{w}{\mathsf{c}^2}. 
\end{align*}

$W, Y, X, V, w$\  \  \  : the alternative variables for $F, A, \Pi, K, u$ defined by

\begin{align*}
&F=\frac{1}{\mathsf{c}^2}\Phi_{\mathsf{N}}-\frac{1}{\mathsf{c}^4}W,\quad A=\frac{1}{\mathsf{c}^3}\varpi^2Y, \\
&\Pi=\varpi\Big(1+\frac{1}{\mathsf{c}^4}X\Big), \quad K=\frac{1}{\mathsf{c}^4}V, \\
&u=u_{\mathsf{N}}+\frac{1}{\mathsf{c}^2}w 
\end{align*}

Note that the variables $W$ and $w$ are related through \eqref{0313} with \eqref{0314}
so that $w=W+O(|\Omega|)$, or $w=W$ where $\Omega=0$. In fact, $F=G$ where $\Omega=0$, $\displaystyle G+\frac{u}{\mathsf{c}^2}=\mbox{Const.}$, and
$\displaystyle u_{\mathsf{N}}=\Phi_{\mathsf{N}}+\frac{\Omega^2}{2}\varpi^2+\mbox{Const.}$.\\

$f, k, l,m$ \  \  \  :the alternative coefficients of the metric (Lewis form):
$$ds^2=f(\mathsf{c}dt)^2-2k\mathsf{c}dtd\phi
-ld\phi^2-e^m(d\varpi^2+dz^2)$$
so that
$$f=e^{2F},\quad k=-e^{2F}A,\quad m=2(-F+K),\quad
l=-e^{2F}A^2+e^{-2F}\Pi^2 $$

$\theta\Big(\xi; \frac{1}{\gamma-1}\Big)$\  \  \  : the Lane-Emden function of index
 $\frac{1}{\gamma-1}$

$\Theta\Big(\xi,\zeta; \frac{1}{\gamma-1},\mathsf{b}\Big)$\  \  \  : the distorted Lane-Emden function

$$ r=\sqrt{\varpi^2+z^2},\quad \zeta=\frac{z}{r} $$
$$\xi=|\mbox{\boldmath$\xi$}|,\qquad \mbox{\boldmath$\xi$}=(\xi_1,,\cdots,\xi_{n-1},\xi_n)^{\top},\quad (n=3,4,5) $$
$$\varpi=\mathsf{a}\sqrt{(\xi_1)^2+\cdots (\xi_{n-1})^2},\quad z=\mathsf{a}\xi_n,
\quad 
r=\mathsf{a}\xi $$

$$u_{\mathsf{N}}(\varpi, z)=u_{\mathsf{O}}\Theta\Big(\frac{r}{\mathsf{a}}, \zeta; \frac{1}{\gamma-1},\mathsf{b}\Big) $$

$$\Omega=\Omega(\varpi, z)=\Omega_{\mathsf{O}}\chi\Big(\frac{r}{R_0}\Big) \quad \chi: \mbox{the cut-off function}: $$
$$\chi(t)=1\quad\mbox{for}\quad t\leq 1, \quad
0<\chi(t) <1\quad\mbox{for}\quad 1<t<2,\quad \chi(t)=0\quad\mbox{for}\quad
2\leq t.$$


\vspace{20mm}









\begin{thebibliography}{99}
\bibitem{Andersson2008} L. Andersson, R. Beig and B. G. Schmidt, Static self-gravitating elastic
bodies in Einstein gravity, Commun. Pure Appl. Math., LXI(2008), 0988-1023.
\bibitem{Andersson2016} L. Andersson, T. Oliynyk and B. G. Schmidt, Dynamical compact
elastic bodies in general relativity, Arch. Rational Mech. Anal., 220(2016), 849-887.
\bibitem{Beig} R. Beig and B. G. Schmidt, Relativistic elasticity, Classical and Quantum Gravity, 20(2003), 889-904.
\bibitem{BoyerL} R. H. Boyer and R. W. Lindquist, Maximal analytic extension of the Kerr metric, J. Mathematical Physics, 8(1967), 265-281.
\bibitem{ButterworthI} E. M. Butterworth and J. R. Ipser, On the structure and stability of rapidly rotating fluid bodies in general relativity. I. The numerical method for computing
structure and its application to uniformly rotating homogeneous bodies,
Astrophysical J., 294(1976), 200-223.
\bibitem{ChandraSS} S. Chandrasekhar, An Introduction to the Study of Stellar Structure, Univ. Chicago Press, Chicago, 1936.
\bibitem{ChandraEFE} S. Chandrasekhar, Ellipsoidal Figures of Equilibrium, Yale UP, New Haven, Conn., 1969; Dover, 1987.
\bibitem{ChandraBH} S. Chandrasekhar, The Mathematical Theory of Black Holes, Oxford UP,
Oxford, 1992.
\bibitem{GilbargT} D. Gilbarg and N. S. Trudinger, Elliptic Partial Differential Equations of Second Order, Springer, 1998.
\bibitem{Hagihara} Y. Hagihara, Theories of Equilibrium Figures of a Rotating Homogeneous Fluid Mass, NASA, Washington DC, 1970.
\bibitem{Hernandez1967} W. C. Hernandez, Jr., Static, axially symmetric, interior solution in general relativity, Physical Review, 153(1967), 1359-1363.
\bibitem{Hernandez1967b} W. C. Hernandez, Jr., Material sources for the Kerr metric, Physcal Review, 159(1967), 1070-1072.
\bibitem{Hernandez2016} J. L. Hernamdez-Pastora, L. Herrera and J. Mart\'{i}n, Axially symmetric static sources of gravitational field, Classical Quantum Gravity, 33(2016), 235005.
\bibitem{Hernandez2017} J. L. Hernandez-Pastora and L. Herrera, Interior solution for the Kerr metric,
Physical Review D, 95(2017), 024003.
\bibitem{Herrera2013} L. Herrera, A. Di Prisco, J. Ib\'{a}\~{n}ez and J. Ospio, Axially symmetric sources: A general framework and some analytical solutions,
Physical Review D, 87(2013), 024014.
\bibitem{Islam} J. N. Islam, {\it Rotating Fields in General Relativity}, Cambridge UP., 1985.
\bibitem{Jardetzky} W. S. Jardetzky, Theories of Figures of Celestial Bodies, Interscience Publ., New York, 1958; Dover, 2005.
\bibitem{Israel} W. Israel, Sources of the Kerr metric, Physical Review D, 2(1970), 641-646.
\bibitem{JJTM1} Juhi Jang and T. Makino, On slowly rotating axisymmetric solutions of the Euler-Poisson equations,
Arch. Rational Mech. Anal., 225(2017), 873-900.
\bibitem{JJTM2} Juhi Jang and T. Makino, On rotating axisymmetric solutions of the Euler-Poisson equations, J. Differential Equations, 266(2019),3942-3972.
\bibitem{Kellog} O. D. Kellog, Foundations of Potential Theory, Springer, Berlin, 1929.
\bibitem{Krasinski1978} A. Krasi\'{n}ski, Ellipsoidal space-time, sources for the Kerr metric, Annals of Physics, 112(1978), 22-40.
\bibitem{Lanczos} K. Lanzcos, \"{U}ber eine station\"{a}re Kosmologie
im Sinne der Einsteinschen Gravitationstheorie,
Zeitschrift f\"{u}r Physik, 21 (1924), 73-110.
\bibitem{Lewis} T. Lewis, Some special solutions of the equations of axially symmetric gravitational fields, Proc. Roy. Soc. London, Ser. A, 136(1931), 176 -192.
\bibitem{TM1998} T. Makino, On spherically symmetric stellar models in general relativity, Kyoto J. Math., 38(1998), 55-69.
\bibitem{ssEE} T. Makino, On spherically symmetric solutions of the Einstein-Euler equations, Kyoto J. Math., 56(2016), 243-282.
\bibitem{asEE} T. Makino, On slowly rotating axisymmetric solutions of the Einstein-Euler equations, J. Math. Physics, 59(2018), 102502.
\bibitem{Meinel} R. Meinel, M. Ansorg, A. Kleinw\"{a}rter,
G. Neugebauer and D. Petroff, {\it Relativistic Figures of Equilibrium}, Cambridge UP., 
2008.
\bibitem{MisnerTW} Ch. W. Misner, K. S. Thorne and J. A. Wheeler, Gravitation, Freeman, New York, 1970.
\bibitem{Moritz} H. Moritz, The Figure of the Earth, Wichmann, Karlsruhe, 1990.
\bibitem{Papapetrou} A. Papapetrou, Champs gravitationnels
stationnaires \`{a} symm\'{e}trie axiale, Ann. Inst. H. Poincar\'{e}, 4(1966), 83-105.
\bibitem{Pizzetti} P. Pizzetti, Principi della Teoria Meccanica della Figura dei Pianeti, Enrico Spoerri, Pisa, 1913.
\bibitem{PlebanskiK} J. Pleba\'{n}ski and A. Krasi\'{n}ski, An Introduction to General Relativity and Cosmology, Cambridge U. P., Cambridge, NY, et al., 2006.
\bibitem{RendallS} A. D. Rendall and B. G. Schmidt, Existence and properties of spherically symmetric static fluid bodies with a given equation of state, Classical Quantum Gravity, 8(1991), 985-1000.
\bibitem{Roos} W. Roos, On the existence of interior solutions in general relativity, Gen. Rel. Grav., 7(1976), 431-444.
\bibitem{Synge} J. L. Synge, Relativity: The General Theory, North-Holland Publ. Co., Amsterdam-New York-Oxford, 1960.
\bibitem{Wavre} R. Wavre, Figures Plan\'{e}taires et G\'{e}od\'{e}sie, Gauthier-Villars, Paris, 1932.
\bibitem{Weyl} H. Weyl, Zur Gravitationstheorie, Ann. Physik, 54(1917), 117-145.
\bibitem{Zeldovich} Ya. B. Zeldvich and I. D. Novikov, Relativistic Astrophysics, 1: Stars and Relativity, Univ. Chicago Press, Chicago, 1971.
\end{thebibliography}
\end{document}